\documentclass[final]{amsart}

\usepackage{amsmath}
\usepackage{amssymb,amsfonts,bm}
\usepackage{color,graphicx,caption,subfig}
\usepackage{xspace}
\usepackage[usenames,dvipsnames]{xcolor}

\begin{document}

\newtheorem{remark}{Remark}[section]
\makeatletter
\define@key{Gin}{Trim}
           {\let\Gin@viewport@code\Gin@trim\expandafter\Gread@parse@vp#1 \\}
\makeatother

\newcommand{\trimval}{20 15 40 35}
\newcommand{\trimvalr}{20 15 40 35}
\newcommand{\trimvalc}{10 10 0 10}

\definecolor{lgray}{gray}{0.95}

\newcommand{\tab}{\hspace{20pt}}
\newcommand{\step}[1]{\noindent\raisebox{1.5pt}[10pt][0pt]{\tiny\framebox{$#1$}}\xspace}
\newcommand{\two}{0.4\textwidth}
\newcommand{\three}{0.3\textwidth}
\newcommand{\threep}{0.3\textwidth}
\newcommand{\four}{0.22\textwidth}

\newcommand{\bmu}{{\bm\mu}}
\newcommand{\x}{{\bm x}}
\newcommand{\dom}{\mathbb P}
\newcommand{\Real}{\mathbb R}

\newcommand{\Sbb}{\mathbb S}
\newcommand{\Wbb}{\mathbb W}
\newcommand{\myspan}[1]{\text{span}\{{#1}\}}

\newcommand{\tN}{{\tilde N}}

\newcommand{\vrb}{v_{\tt rb}}
\newcommand{\Xitrial}{\Xi_{\tt train}}
\newcommand{\Xitest}{\Xi_{\tt test}}
\newcommand{\Xitrialper}{\Xi^{\%}_{\tt trial}}
\newcommand{\XiH}{\Xitrial}

\newcommand{\dbar}[1]{\overline{\overline{#1}}}

\newcommand{\ben}[1]{\color{red} {#1} \color{black}}
\newcommand{\new}[1]{\color{blue} {#1} \color{black}}

\newcommand{\mybox}[3]{
\begin{center}
\begin{minipage}[t]{\dimexpr\textwidth}
\begin{center}
\fcolorbox{black}{lgray}{%
\begin{minipage}[t]{\dimexpr\textwidth-2\fboxsep-2\fboxrule}
#1
\end{minipage}
}
\fcolorbox{black}{white}{%
\begin{minipage}[t]{\dimexpr\textwidth-2\fboxsep-2\fboxrule}
#2
\vspace{4pt}
\begin{center}
#3
\end{center}
\end{minipage}
}
\end{center}
\end{minipage}
\end{center}
}

\newcommand{\mysymb}[1]{
\vspace{5pt}
\noindent
\step{\tt #1}
}

\title{Locally Adaptive Greedy Approximations for Anisotropic Parameter Reduced Basis Spaces}


\author{Yvon Maday}
\address{
UPMC University Paris 06, UMR 7598, Laboratoire Jacques-Louis Lions, F-75005, Paris, France, and
 Division of Applied Mathematics, Brown University, Providence RI, USA.}
\email{\tt maday@ann.jussieu.fr}

\author{Benjamin Stamm}
\address{UPMC University Paris 06, UMR 7598, Laboratoire Jacques-Louis Lions, F-75005, Paris, France and CNRS, UMR 7598, Laboratoire Jacques-Louis Lions, F-75005, Paris, France.}
\email{\tt stamm@ann.jussieu.fr}

\keywords{
model reduction, reduced basis method, greedy algorithm
}

\begin{abstract}
Reduced order models, in particular the reduced basis method, rely on empirically built and problem dependent basis functions that are constructed during an off-line stage. 
In the on-line stage, the precomputed problem-dependent solution space, that is spanned by the basis functions, can then be used in order to reduce the size of the computational problem.
For complex problems, the number of basis functions required to guarantee a certain error tolerance can become too large in order to benefit computationally from the model reduction. 
To overcome this, the present work introduces a framework where local approximation spaces (in parameter space) are used to define the reduced order approximation in order to have explicit control over the on-line cost. This approach also adapts the local approximation spaces to local anisotropic behavior in the parameter space. We present the algorithm and numerous numerical tests.
\end{abstract}

\maketitle


\pagestyle{myheadings}
\thispagestyle{plain}
\markboth{Yvon Maday and Benjamin Stamm}{Local Greedy Approximations}

\section{Introduction}
\label{sec:introduction}

The recent progresses in the numerical simulation of physical phenomena obtained through the combination of robust and accurate approximation methods and faster and larger computational platforms have permitted to investigate more and more complex problems with improved reliability. These progresses have, in turn, lead to new demands for the numerical simulations that are not only used to understand a given state but investigate control and optimization problems. 
In such problems, a large number of similar computations of a parameter-dependent model have to be performed associated with 
different values of the parameter chosen either randomly of guided by some well suited recursive algorithm.

Faster solution algorithms are often not sufficient to achieve these new demands in many engineering applications and reduced numerical methods have been proposed as surrogates to standard numerical high fidelity approximations of given mathematical models. 
Among these methods, the reduced basis (RB) methods \cite{Rozza:2008p921,Prudhomme:2002p183,Quarteroni:2011p6628} is a class of higher-order mathematical technique that has been widely developed and gained in generality and reliability during the last decade. 
The greedy algorithm \cite{Veroy:2003p56} plays a crucial ingredient and will be the main focus of the present paper.
For parametrized time-dependent problems, the RB (greedy-based) strategy can be successfully combined with the proper orthogonal decomposition (POD) \cite{book_POD_Holmes} as is illustrated in \cite{HO_2008,NRP_2009}.

The basic idea behind these reduced order methods is the notion of small Kolmogorov $n$-width of the set of all solutions under variation of the parameter, i.e. that the solution manifold, embedded in a Hilbert space $\mathbb W$, can be well approximated with a well chosen $n$-dimensional subspace of $\mathbb W$.
In essence, RB methods are based on a two steps strategy : the first step (off-line stage) allows to select particular instances of the parameters, for which a very accurate approximation of the solution is computed : those solutions constitute the reduced basis.
In a second step (the on-line stage), the generic solutions (for other instances of the parameter)  are approximated by a linear combination of these reduced basis functions.
The interest of the approach lies in the fact that, in cases where the Kolmogorov $n$-width is small --- and there are many indications for such a smallness: (i) a priori, such as regularity of the solutions as a function of the parameters when the set of parameters lies in a finite dimensional space, and (ii) a posteriori, revealed by a large number of simulation --- these reduced model approximations require very few degrees of freedom and in some sense are close in spirit to other high order approximation methods like spectral methods, but with an improved efficiency for what concerns the computational complexity.

It is well known that high order methods generally take advantage of a global approach by using basis functions that have a large support that, combined with higher regularity of the solution to be approximated (going some times up to analyticity) allows to get, with very few degrees of freedom, a very good accuracy. For most practical design problems in engineering though,  the solution is not analytical and most of the time regularity exists but only locally which precludes the interest of global approximations. This is the reason why, for instance, by breaking the global framework to locally piecewise global approaches, the spectral element method reveals superiority with respect to the plain spectral method: a trade off between locality and globality is generally preferred as is demonstrated in e.g.  \cite{Cantwell:2011p7419} for approximation in spatial direction by spectral element approximations. The same causes imply the same effects in the parameter directions for RB approximations, this is the reason why recently, some works have been devoted to proposing ways of adding locality to the reduced methods yielding parameter subdomain domain refinement.

In this paper we shall focus on the certified RB framework for which, the construction of basis functions during the first stage of the algorithm is, as in most of the current approaches, performed  through a greedy strategy based on an a posteriori error estimator. 
For either high dimensional parameter spaces or spaces with large measure one may encounter the problem that the size of the reduced basis turns out to be larger than desired. 
Following the lines drawn above, first ideas in this context has been presented by Eftang, Patera and R\o nquist \cite{Eftang:2010p6616}
and further in \cite{NME:NME3327,Fares:2010p4573,HDO_2011} where the parameter space is decomposed into cells where different reduced basis sets are assembled.  This approach presents clear advantages in the size of the matricial system that appears in the on-line solution procedure, and corroborates the natural feeling that, in order to approximate the solution at a given parameter, primarily those solutions in the reduced basis corresponding to parameters that are close to the parameter we are interested in are to be involved in the linear approximation. 

A drawback of the current approach \cite{Eftang:2010p6616,NME:NME3327,Fares:2010p4573,HDO_2011} however is, that in two adjacent parameter-subdomains, some of the parameters that are selected may be very close or even lie both on the common interface. This is a situation that may be encountered very often when dealing with high dimensional parameter spaces, and in addition it is known from computations that the greedy algorithm has the tendency to pick sample points lying on the boundary of the parameter space.
In this case, the approximation on one of the subdomains does not benefit of the computation that was done for the other subdomain.
Parameter domain decomposition may  thus not be the ultimate approach. 
Another drawback of all classical greedy approaches is to be unable to master the size of the reduced basis and in consequence the size of the discrete system that will be solved in the on-line procedure. 

These reflections have motivated us to investigate the alternative discussed in this paper which introduces a completely new framework, which should be seen as a first layout of our ideas where it is unfortunately impossible to answer all interesting questions related to this approach without losing the big picture.

The outline of the paper is as follows.
In Section \ref{sec:StdGreedy} we present the standard greedy algorithm that is widely used in reduced basis computations in order to have a common ground for the further discussion. Then, we introduce in Section \ref{sec:LA} the concept of local approximation spaces, that is, if a reduced order approximation needs to be computed at a parameter point $\bmu$ the approximation space is spanned by the $N$ closest precomputed basis functions.
Thus, the construction of the local approximation spaces is based on a distance function which is constructed in Section \ref{sec:metric} in an empirical manner taking into account the geometry of the manifold of the solution set including the local anisotropies in the parameter space.
Section \ref{sec:hessian} introduces the concept of a varying training set that is chosen in an optimal fashion using the possibly anisotropic distance function if the computing resources are limited to handle a given (and too small) number of training points.
In Section \ref{sec:online} we explain some practical aspects of the on-line implementation, in particular how the different snapshots can be orthonormalized at the on-line stage. 
In Section \ref{sec:L2}  we present some numerical examples which uses the best approximation ($L^2$-projection) for an explicitly given family of functions whereas an example using a reduced basis framework is given in Section \ref{sec:Galerkin}.
Finally, in Section \ref{sec:conc} we draw the conclusions.

\section{The greedy algorithm for reduced basis approximations}
\label{sec:StdGreedy}
Let us first introduce a classical greedy algorithm
to have a common ground to present our ideas.
In a general setting, it is assumed that for each parameter value $\bmu$ in a parameter domain $\dom \subset \Real^p$,
a $\bmu$-dependent function $v(\bmu)\in\Wbb$ of the variable $\bm x$ in a bounded domain $\Omega\subset  \Real^d$ (e.g. $d=2$ or $d=3$)  can be computed.
The space $\Wbb$ denotes some functional space of functions defined on $\Omega$ and we denote the solution manifold by
\[
	\mathbb W_{\mathbb P} = \{v(\bmu)\in\Wbb\,|\,\bmu\in\dom\}.
\]
In order to fix the ideas this can be that $v(\bmu)$ is a solution to a parameter dependent partial differential equation, an approximation of which can be computed by the e.g. classical finite element or spectral method.
Let
\[
\Sbb_N =\{\bmu^1,\cdots,\bmu^N\}
\]
be a collection of $N$ parameters in $\dom$ and
\[
\Wbb_N = \myspan{v(\bmu^1),\cdots,v(\bmu^N)}
\]
be the approximation space associated to the set $\Sbb_N$. A projection operator $P_N:\Wbb\rightarrow \Wbb_N$ is supposed to exist that can be either a parameter dependent Ritz-projection (reduced basis methods in case of a parameter dependent PDE) \cite{Prudhomme:2002p183,Sen:2006p15}
or simply a $L^2$-projection.
Further, it is assumed that for each parameter value $\bmu \in \dom$ an error estimator $\eta(\bmu; \Wbb_N)$ of the approximation of $v(\bmu)$ by $P_N(v(\bmu))$  can be computed (this can be e.g. through a posteriori analysis of the residual, see eg \cite{Prudhomme:2002p183}).
We assume that there exists two constants $c>0$ and $C>0$ such that
\[
	c \, \eta(\bmu; \Wbb_N)
	\le 
	\|v(\bmu) - P_N(v(\bmu)) \|
	\le 
	C \, \eta(\bmu; \Wbb_N),
\]
for an appropriate norm on $\mathbb W$.
The case where $\eta(\bmu; \Wbb_N)$ describes the exact error is thus not excluded in this general framework.

The following sketch represent a typical greedy algorithm: Let $\Xitrial$ be a chosen finite training set $\Xitrial \subset\dom$ of representative points.

\mybox{
{\tt ClassicalGreedy}\\
{\tt Input: } ${\tt tol}$.\\
{\tt Output: } $N$, $\mathbb W_N$, $\mathbb S_N$.
}
{
\begin{itemize}
\item[\tt 1.] Choose (possibly randomly) $\bmu^1\in\Xitrial$, set $\Sbb_1=\{\bmu^1\}$, $\Wbb_1=\myspan{v(\bmu^1)}$ and $N=1$,
${\tt err}=\max_{\bmu \in \Xitrial}  \eta(\bmu; W_1)$.
\item[\tt 2.] {\tt While} $ {\tt err}>{\tt tol}$
\item[\tt 3.] \tab Find $\bmu^{N+1}=\mbox{argmax}_{\bmu \in \Xitrial}  \eta(\bmu; W_N)$, ${\tt err}=\max_{\bmu \in \Xitrial}  \eta(\bmu; W_N)$.
\item[\tt 4.] \tab Compute $v(\bmu^{N+1})$, set $\Sbb_{N+1} = \Sbb_N\cup\{\bmu^{N+1}\}$ and  set
\item[] \tab\tab
$\Wbb_{N+1} = \myspan{\Wbb_N ,  v(\bmu^{N+1})}$.
\item[\tt 5.] \tab Set $N := N+1$.
\item[\tt 6.] {\tt End while}.
\end{itemize}
}{
Algorithm 1: Classical greedy algorithm.
}

\vspace*{11pt}
\begin{remark}
The introduction of the finite set $\Xitrial$ is due to practical implementation, because the entire scan of $\dom$ is impossible: $\Xitrial$ has to be finite, not too large in order the previous algorithm is not too lengthy, not too small in order to represent well $\dom$. Actually the definition of this finite set can evolve during the algorithm, see \cite{HDO_2011}. We shall elaborate on this in the framework of our approach later in the paper (Section \ref{sec:hessian}).
An alternative approach is presented in \cite{BuiThanh:2008p7559}, where finding the maximum is stated as a constraint optimization problem.
\end{remark}

Let $\tN$ denote the final size of the reduced basis space $\Wbb_{\tN}$ such that the accuracy criterion 
\[
	\max_{\bmu \in \Xitrial}  \eta(\bmu; \Wbb_\tN) \le tol
\]
is satisfied. Note that this final size is a consequence of the tolerance that has been chosen, hence, for a prescribed tolerance, we have not the hand on the complexity of the RB approximation method.

\section{Locally adaptive greedy approximations}
\label{sec:LA}
The new approach pursued in this work is in the same spirit of partitioning the parameter set as in  \cite{Eftang:2010p6616,NME:NME3327,Fares:2010p4573,HDO_2011}.
However, the difference is listed below:
\begin{itemize}
\item[\mysymb{1}] 
We do not impose a clear partition of the parameter space but rather collect a global set of $K$ sample points $\Sbb_K$, they are preliminary selected in the {\it off-line} stage and the corresponding basis functions are computed and stored. This step is elaborated in Section \ref{ssec:SamplePoints}.
\item[\mysymb{2}] 
We choose a priori the size $N$ ($\le K$) of the system we want to solve {\it on-line}.
Thus, the resulting complexity of the {\it on-line} stage is user controlled.
\item[\mysymb{3}] 
When a reduced basis approximation is to be computed for a certain given parameter value $\bmu\in\dom$ in the {\it on-line} stage, we only use the $N$ precomputed basis functions whose corresponding parameter values lie in a ball around $\bmu\in\dom$.
This step is elaborated in Section \ref{ssec:LAS}.
\item[\mysymb{4}] 
The distance function used to define the ball, containing the $N$ closest basis functions, is an empirically constructed distance function that is problem dependent and can reflect anisotropies in parameter space and its local changes. 
The construction of the distance function is presented in Section \ref{sec:metric}.
\end{itemize}

\subsection{Local approximations spaces}
\label{ssec:LAS}
As already anticipated, for a certain given parameter value $\bmu\in\dom$ in the {\it on-line} stage
we only use the $N$ basis functions whose parameter values lie in a ball
\begin{equation}
	\label{eq:DefBall}
	B_{\bmu} = \{\tilde \bmu\in\dom\,|\,  d(\bmu,\tilde\bmu) \le r(\bmu)\},
\end{equation}
for a given distance function $d(\cdot,\cdot)$. 
The radius $r(\bmu)$ is tuned in such a way that there are actually $N$ basis functions in the ball. 
Denoting by $\mathbb S_K$ the set of sample points for which the corresponding functions are precomputed at the off-line stage, the local sample space is defined by
\[
	\Sbb_{\bmu} = B_{\bmu}\cap\Sbb_K = \{\tilde \bmu\in\Sbb_K \,|\, \tilde\bmu\in B_{\bmu}\}
\]
with cardinality equal to $N$\footnote{Note that $N$ has not the same meaning as in the previous section any more; hence the cardinality of $\Sbb_K$ is not $N$ but a $K>N$ generally much larger.}. Further the local reduced basis approximation space shall be defined by $\Wbb_{\bmu}=\myspan{v(\tilde\bmu)\,|\, \tilde\bmu\in\Sbb_{\bmu}}$ and its associated local projection by $P_{\bmu}:\Wbb\rightarrow \Wbb_{\bmu}$ in order to define the reduced basis approximation $ P_{\bmu}(v(\bmu))$ to $v(\bmu)$.
Also in this local framework, the Galerkin setting can be used to design a posteriori estimates to quantify the error between $v(\bmu)$ and $P_{\bmu}(v(\bmu))$, that we shall denote by $\eta(\bmu; \Wbb_\bmu)$.

A similar approach is presented in \cite{OnlineGreedy} where the snapshots are chosen also in the on-line stage but such that the approximation error is maximally reduced. Thus, in general, no local neighborhood is required.

One basic ingredient of this local reduced basis approximation is thus the distance function $d(\cdot,\cdot)$; the radius is adjusted to select the total number of modes $N$ for each value of $\bmu$ whereas the distance function can be chosen isotropic or better accounts for anisotropies in the parameter space.
The construction of an adapted local distance function is a topic itself and is addressed in Section \ref{sec:metric}.
Another ingredient is the set of sample points $\mathbb S_K$ whose construction is explained in the next section.

\subsection{Construction of global set of sample points}
\label{ssec:SamplePoints}
This consists of the off-line stage of the algorithm. The first decision is to set the number $N$ of basis functions we want to get involved in the on-line procedure (i.e. the size of the matrix system to be solved on-line).

A standard greedy algorithm, as described by Section \ref{sec:StdGreedy}, is then performed in a first stage until $N+1$ basis functions are selected. 
The corresponding parameter values define the current set $\Sbb_{N+1}$ and the solutions $v(\bmu)$, $\bmu\in\Sbb_{N+1}$, are saved.

Then, the second stage of the algorithm is performed to further enrich $\Sbb_{N+1}$. 
The construction of the further basis functions require the data of a distance function (as in the on-line stage).
As for now we assume that the distance function is available, its recursive construction is explained in Section  \ref{sec:metric} (note that, by default, we use the Euclidean isotropic metric). 
We note that we can at each iteration and for each parameter value $\bmu\in\dom$, compute an approximation with the $N$ basis functions that lie in $B_\bmu$ and estimate the error by computing $\eta(\bmu; \Wbb_\bmu)$. 
The next parameter value(s) is chosen as the one corresponding to the worse approximation properties.
It is important to notice that using local approximation spaces, several basis functions can be chosen per iteration. The following procedure is proposed. 

At the end of a given iteration where we dispose of $\bmu_1,\ldots,\bmu_K$, a new iteration may start.
The next sample point, $\bmu_{K+1}$, is chosen such that the estimated error $\eta(\bmu; \Wbb_\bmu)$ is maximized over the parameter space as in the traditional setting except that only a local approximation space $\Wbb_\bmu$ is considered to compute the error (-estimator). 
Thus ${\bmu_{K+1}}=\arg\max_{\bmu\in\Xitrial}\eta(\bmu; \Wbb_\bmu)$. 
Next, the maximum estimated error $\eta(\bmu; \Wbb_\bmu)$ is searched over all parameter values lying outside of the domain of influence of  $\bmu_{K+1}$.
By the domain of influence of $\bmu_{K+1}$ we mean 
\[
	\mathcal D_{\bmu_{K+1}}=\{\bmu\in\dom\,|\, \bmu_{K+1}\in B_{\bmu}\}\cup B_{\bmu_{K+1}}. 
\]
Outside of the domain of influence of $\bmu_{K+1}$ more basis functions can however be added. 
This condition assures that at most one sample point is added per ball. This procedure is repeated, and more and more parts of the parameter domain are excluded, until the remaining set is the empty set. 
The practical implementation of this procedure is given the following pseudo-code.
\mybox{
{\tt EnrichmentLoop}\\
{\tt Input:} $K$, $\Sbb_K$\\
{\tt Output:} updated $K$, $\Sbb_K$
}{
\begin{itemize}
\item[\tt 1.] Set $\widetilde\Xitrial = \Xitrial$, ${\tt err}=\mbox{max}_{\bmu \in \widetilde\Xitrial}\eta(\bmu; \Wbb_\bmu)$.
\item[\tt 2.] {\tt While} $\widetilde\Xitrial\ne\emptyset$ and ${\tt err}>{\tt tol}$
\item[\tt 3.] \tab Set $\bmu^{K+1} = \mbox{argmax}_{\bmu \in \widetilde\Xitrial}  \eta(\bmu; \Wbb_\bmu)$, ${\tt err}=\mbox{max}_{\bmu \in \widetilde\Xitrial}\eta(\bmu; \Wbb_\bmu)$.
\item[\tt 4.] \tab Compute $v(\bmu^{K+1})$ and set $\Sbb_{K+1} = \Sbb_K\cup\{\bmu^{K+1}\}$.
\item[\tt 5.] \tab Set $\widetilde\Xitrial := \widetilde\Xitrial\backslash
\{\bmu\in\Xitrial\,|\, \bmu\in B_{\bmu^{K+1}} \,\text{ or }\,  \bmu^{K+1}\in B_\bmu\}$.
\item[\tt 6.] \tab Set $K := K+1$.
\item[\tt 7.] {\tt End while}.
\end{itemize}
}{
Algorithm 2: Basis enrichment loop.
}

\section{Anisotropy and distance function - a general finite difference based approach}
\label{sec:metric}
In order to equip the parameter space $\dom$ with an appropriate metric to $\mathbb W_{\mathbb P}$ and distance function, one has to estimate, on the fly, the way the (unknown) function $v(\bmu)\in\mathbb W$ depends on $\bmu$. 
In this manner we aim to reduce at most the number of elements involved in the on-line reduced basis approximation.

Indeed, for instance assuming that the parameter domain $\dom$ is a subset of $\Real^2$, so that $\bmu=(\mu_1,\mu_2)$, and assuming that there exists a function $\varphi$ of one variable, such that $v(\bmu) = \varphi (\mu_1+\mu_2)$, then the optimal selection of parameters should be sought e.g. along the line $\mu_1=\mu_2$ considering only a one dimensional parameter. 
Of course the fact that $v$ has such a behavior is generally not known. 
This has to be figured out from prior computations, during the greedy algorithm. 
This is actually a quite classical quest in numerical analysis, and this valuable knowledge leads to a large computational reduction. As in other instances, e.g. adaptive approaches like the one used in finite element approximations (see e.g. \cite{Danaila:2010p7421}),  the construction of such a distance function $d(\cdot,\cdot)$ that accounts for the anisotropic behavior of changes of $v(\bmu)$ with respect to variations in $\bmu$.
Similarly as in other approaches, the distance function is derived from the knowledge of the Hessian of the function $v$ with respect to the parameter $\bmu$. 

We first propose  a general framework to obtain an approximation of the Hessian by finite differences and then give the definition of the distance function. 
In what follows, we always assume that $\Xitrial$ is such that its convex hull spans the parameter space $\dom$. 

\subsection{Definition of the Hessian} 
\label{ssec:Hessian}
The goal is to define a Hessian matrix ${H}(\bmu)$ for each point $\bmu\in \XiH$ upon which the distance function will be based on. Since the Hessian is based on the reduced basis approximation, it is updated/constructed at each iteration in the off-line greedy algorithm (Section \ref{ssec:SamplePoints}).
To do so, compute a reduced basis approximation $P_{\bmu}(v(\bmu))$ of $v(\bmu)$ based on the current local reduced basis, as explained in Section \ref{ssec:LAS}, on the stencil  $\bmu(\alpha) = \bmu +  \sum_{i=1}^p \alpha^i \delta \mu^i$, with $\bmu\in \XiH$ and where $ \delta \mu^i$, $i = 1,..,p$,  is a positive small increment in the $i^{th}$ direction in the parameter and $\alpha^i = -1, 0, 1$, such that, at most, two of them are non zero. 
Let $\{\varphi_n\}_{n=1}^N$ denotes the basis of $\Wbb_\bmu$ used to build $P_{\bmu}(v(\bmu))$ such that
\[
	P_{\bmu}(v(\bmu))= \sum_{n=1}^N v_n(\bmu) \varphi_n(\x).
\]
From the $3^p$ approximations $P_{\bmu}(v(\bmu))$ of $v(\bmu(\alpha))$ on the surrounding stencil $\bmu(\alpha)$ of $\bmu$, one can, by finite differences, define a Hessian matrix with value in $\mathbb R$ as follows
\[
	{H}_{ij} (\bmu) = \sum_{n=1}^N v_n(\bmu) \, D_{ij} v_n(\bmu),
\]
where 
\[
	D_{ij}v_n(\bmu) = \frac{v_n(\bmu+\delta\mu^i+\delta\mu^j) - v_n(\bmu-\delta\mu^i+\delta\mu^j) - v_n(\bmu+\delta\mu^i-\delta\mu^j)+ v_n(\bmu-\delta\mu^i-\delta\mu^j)}{4\delta\mu^i\delta\mu^j}.
\]
The term $D_{ij} v_n(\bmu)$ describes the anisotropy of the n-th mode while $v_n(\bmu)$ indicates its weight in the approximation. 
\begin{remark}
In the case one is only interested in a scalar output functional $s:\Wbb\to\mathbb R$ of the solution, i.e. $s(v(\bmu))$ approximated by $s(v_{\tt rb}(\bmu))$, one can also consider the Hessian do be defined as
\[
	{H}_{ij} (\bmu) = D_{ij}s(v_{\tt rb}(\bmu)).
\]
\end{remark}
\begin{remark}
The current construction of the Hessian based on a finite difference scheme involves quite a large number of points when the dimension $p$ becomes too large. Alternative exist based only on the only points in $\mathbb S_K$ that need still to be investigated and improved. 
\end{remark}

At this level we dispose for any $\bmu\in\dom$ of a method to compute ${H}_{ij} (\bmu) $.
It remains to explain how we construct the metric tensor $M$ from the knowledge of the Hessian matrix $ H$. 
Indeed, given $\bmu\in \dom$, we proceed in a standard manner and perform an eigenvalue decomposition 
\[
	 H(\bmu) = V\, \Lambda\, V^T
\]
where $V$ is an orthonormal matrix and $(\Lambda)_{ii} = \lambda_i$ a diagonal matrix consisting of the eigenvalues $\lambda_i$. Consider the diagonal matrix
\[
	|\Lambda|_{ii} = { |\lambda_i|}, \qquad i=1,\ldots,p,
\]
and the associated symmetric semi-positive definite matrix $M(\bmu)= V\, |\Lambda|\, V^T$ to define the metric tensor.

\subsection{Construction of the distance function} 
The construction of the distance between two parameter points $\bmu_1,\bmu_2\in\dom$, denoted by $d(\bmu_1,\bmu_2)$, in the metric space given by the metric tensor $M$ is ideally defined by the length of the geodesic curve linking $\bmu_1$ and $\bmu_2$. Indeed, if $\bm\gamma:[0,1]\to \dom$ is a parametrization of the geodesic from $\bmu_1$ to $\bmu_2$ we would like to compute
\[
	\int_0^1 \sqrt{(\bm\gamma'(t))^T\,M(\bm\gamma(t))\,\bm\gamma'(t)}\, dt.
\]
The above integral is however not computable due to the lack of knowledge of the geodesic curve. 
The geodesic curve itself could be approximated, but this would slow down the on-line computation drastically. 
We therefore propose to approximate the above integral roughly by the following trapezoidal rule to define the distance
\[
	d(\bmu_1,\bmu_2) = \frac12 \sqrt{(\bmu_2-\bmu_1)^T\, M(\bmu_1)\,(\bmu_2-\bmu_1)} 
		+ \frac12  \sqrt{(\bmu_2-\bmu_1)^T\, M(\bmu_2)\,(\bmu_2-\bmu_1)}
\] 
where we also took the approximated $\bm\gamma'(t) \approx (\bmu_2-\bmu_1)$ for $t=0,1$ (which is correct when dealing with the Euclidian metric).

We note that the construction of the Hessian $H(\bmu)$ and thus the distance function depends on the approximation $P_{\bmu}(v(\bmu))$, and thus on $\mathbb S_K$. That is, for each reduced basis model we can define a distance function and  we recursively update the distance function at each step of the stage 2 of the greedy algorithm. 
Therefore, the distance function is, simultaneously with $P_{\bmu}(v(\bmu))$, converging to a final distance function that is used in the on-line procedure.

Let us make the following comments:
\begin{itemize}
\item[\mysymb{1}] 
It is easy to see that our distance function is symmetric by construction.
\item[\mysymb{2}] 
The matrix $M(\bmu)$ is semi-definite positive. Allowing for zero eigenvalues is indeed a desired property since it allows to shrink certain directions in parameter space where there is no variation in the solution with respect to the parameter.
\item[\mysymb{3}]
Defining a mesh consisting of $p$-dimensional simplices with vertices being the set $\XiH$ allows to interpolate/represent quantities only defined on $\XiH$ by means of piecewise linear functions on the mesh for any parameter value $\bmu\in\dom$, since the convex hull of $\XiH$ coincides with the parameter space $\dom$.
We can therefore easily access the interpolated metric tensor $M (\bmu)$, for any $\bmu\in\dom$, by means of its interpolated component functions $M_{ij} (\bmu)$.
\end{itemize}

\begin{remark}
As an alternative to interpolating the coefficients in a $p$-dimensional space one can construct the Hessian on-line as explain above for any $\bmu\in\dom$.
\end{remark}

\section{Adaptive training sets based on Hessian}
\label{sec:hessian}
The goal of this section is to present an approach that further allows to reduce the off-line computational cost. Here, we focus on the training set $\Xitrial$. It is aimed to keep its cardinality as small as possible, but large enough to capture the local geometry of the parametrized system. 
We aim to give an answer to the following question: Given the computational resources of handling $Q_M$ points in $\Xitrial$, how do we place them optimally in the parameter space?

We propose to place the $Q_M$ points uniformly in $\dom$ with respect to the empirically constructed metric tensor $M(\bmu)/r(\bmu)^2$.
The factor $1/r(\bmu)^2$ is motivated by the fact that a ball defined in \eqref{eq:DefBall} corresponds (approximatively) to a unit ball in the metric defined by $M(\bmu)/r(\bmu)^2$.
Since the metric gets updated at each iteration of the greedy algorithm, we also need to adapt the training points $\Xitrial$ at each iteration. 
Further, we propose that the cardinality of $\Xitrial$ is an increasing function of the inverse of the actual error (-estimation) $\tt err$. 
Let $Q_m$ be the minimal cardinality of $\Xitrial$ at the beginning of the greedy algorithm and $Q_M$ the maximal target cardinality at the end of the algorithm once the tolerance criteria ${\tt err} < {\tt tol}$ is reached, then we define
\[
	Q({\tt err}) =  \left\lceil\frac{Q_M-Q_m}{\log({\tt tol})} \log({\tt err})  +  Q_m\right\rceil
\]
to be the number of points at the next iteration in $\Xitrial$.

To summarize, the framework that we propose is that the training set $\Xitrial$ of cardinality $Q({\tt err})$ is constructed such that the points are uniformly distributed in the empirically constructed metric tensor $M(\bmu)/r(\bmu)^2$.

Our implementation of this concept is based on the Delaunay mesh defined by the points of $\Xitrial$. We use the software FreeFem++ \cite{Hecht:2008p27} which provides a uniform mesh with respect to $M(\bmu)$.
Other implementations are of course also possible.

\begin{remark}
Enriching the training set during the sampling process is not a
new idea. An adaptively enriching greedy sampling strategy was proposed in [7] however it is not based, as it is the case in the present paper, on a learning process of the geometry of the manifold. In addition our approach is more precise since it not only detects
where the manifold of the solutions (as a function of the parameters) is complex but
also the directions that are important to sample and those that are not. Second in our
approach the size of the discrete system to be solved on-line is determined by the user
leading to a controlled simulation cost.

\end{remark}

\section{On-line orthonormalization of basis functions}
\label{sec:online}
In the standard approach, it is well known that the stable implementation of the reduced basis method requires the orthonormalization of the snapshots $v(\bmu)$, $\bmu\in \Sbb_N$ \cite{Patera:2006p6}.  This orthonormalization is prepared off-line and the vectorial space $\Wbb_N$ spanned by the series $\{v(\bmu)\,|\, \bmu\in\Sbb_N\}$ is not affected by this change of basis. 

In our alternative method,  there are a large number of approximation spaces, almost one for each $\bmu$ : $\Wbb_{\bmu}=\myspan{v(\tilde\bmu)\,|\, \tilde\bmu\in\Sbb_{\bmu}}$. We cannot compute all various orthonormalized basis sets for all possible  $\Wbb_{\bmu}$ off-line. This would be too costly and using way too much memory.

Nevertheless, we can still {\it prepare} the orthonormalization process so that, the on-line orthonormalization is feasible. For any $\hat\bmu_i$ and $\hat\bmu_j$ in $\Sbb_K$, let us compute off-line the scalar products 
\begin{equation}
M_{ij} = \langle \hat v(\bmu_j), v(\hat \bmu_i)\rangle,\qquad i,j=1,\ldots,K
\end{equation}
where $\langle\cdot,\cdot\rangle$ represents an appropriate scalar product : say either $L^2$ or $H^1$ type. Note that we can even compute only those scalar products corresponding to pairs of parameters $(\hat \bmu_i,\hat\bmu_j)$ that are close, indeed if $\hat\bmu_i$ and $\hat\bmu_j$
are distant in  $\dom$, they will never be in the same $\Sbb_{\bmu}$ for any $\bmu\in\dom$.
 
In the on-line stage, once the approximation in $\Wbb_{\bmu}$ will be required, we can orthonormalize the basis set $\{v(\tilde\bmu)\,|\, \tilde\bmu\in\Sbb_{\bmu}\}$ first by declaring an order in the parameters
\begin{equation}
\left\{ \tilde\bmu\in\Sbb_{\bmu}\right\} =\left\{ \bmu_1,  \bmu_2, . . . ,  \bmu_N\right\} 
\end{equation}
then perform a classical Gram-Schmidt orthonormalization process
\begin{eqnarray}
\label{GS}
\zeta_1 &=&  \beta_{11}v(\bmu_1),\cr
\zeta_2 &=&\beta_{22} v(\bmu_2) + \beta_{21} \zeta_1 ,\cr
&\vdots& \cr
\zeta_N &=&\beta_{NN} v(\bmu_N) + \sum_{i=1}^{N-1} \beta_{Ni} \zeta_i,
\end{eqnarray}
where the coefficients $\beta_{ni}$ are chosen so that $\zeta_n$ is orthogonal to $$\myspan{v(\bmu_i)\,|\, i=1,..,n-1} = \myspan{\zeta_i\,|\, i=1,..,n-1}$$ and the coefficients $\beta_{nn}$ are chosen so that $\zeta_n$ is norm 1. 
It is well known that these coefficients exist, are unique and that, in order to compute them, we have to know each scalar product $\langle  v(\bmu_n), \zeta_i\rangle$, $i=1,..,n-1$ that can be easily computed mimicking the Gram-Schmidt procedure \eqref{GS} but using the precomputed coefficients $M_{ij}$. For the construction of the coefficients $\beta_{nj}$, we proceed as follows (it helps to have in mind that $\tilde \beta_{kn}$ stands for $- \langle v(\bmu_k),\zeta_n\rangle)$ : for all $n=1,\ldots,N$ apply recursively

\begin{eqnarray}
	\label{eq:beta}
 \alpha &=& \left( M_{nn} - \sum_{j=1}^{n-1}|\tilde\beta_{nj} |^2 \right)^\frac12, \cr
 \beta_{nn} &=& \frac{1}{\alpha},   \cr
 \beta_{nj} &=& \frac{\tilde\beta_{nj}}{\alpha}, \qquad k=1,\ldots,n-1,   \cr
 \tilde\beta_{kn} &=& - {\beta_{nn}}\,M_{nk} + \sum_{j=1}^{n-1} {\beta_{nj}}\,\tilde\beta_{kj}, \qquad k=1,\ldots,N.
\end{eqnarray}
Note that the construction of $\alpha$ is not subject to instabilities since the sum involves only positive values.

Next, we notice that the basis functions $\{\zeta_n\}_{n=1}^N$ can alternatively be expressed by the following linear transformation
\begin{eqnarray}
\label{zetas}
  \zeta_1 &=&  \gamma_{11} v(\bmu_1),   \cr
  \zeta_2 &=& \gamma_{22} v(\bmu_2) + \gamma_{21} v(\bmu_1) ,  \cr
&\vdots& \cr
  \zeta_N &=& \sum_{i=1}^{N} \gamma_{Ni} v(\bmu_i),
\end{eqnarray}
from the set of basis functions $\{v(\bmu_n)\}_{n=1}^N$. 
It is our final intention to construct the coefficients $\gamma_{ni}$ representing the change of basis functions from $\{v(\bmu_n)\}_{n=1}^N$ to $\{\zeta_n\}_{n=1}^N$ as defined by \eqref{zetas} in a stable way and without any operation depending on the length of the vectors representing $v(\bmu_n)$ to guarantee stability and good on-line performance.  

\begin{remark}
Note that the coefficients $\gamma_{n,i}$ could be obtained based on \eqref{zetas} directly. Indeed, each $\zeta_n$ is orthogonal to $\myspan{v(\bmu_i)\,|\, i=1,..,n-1}$ allowing to compute the coefficients $\gamma_{ni}$, $i=1,\ldots,n-1$. 
This however involves solving a linear systems based on the mass matrix with coefficients $M_{ij}= \langle v(\bmu_j), v(\bmu_i)\rangle$, $i,j=1,..,n-1$, which can be ill conditioned due to the fact that the set of functions $\{v(\bmu_i)\}_i$ is heavily linearly dependent.
\end{remark}

The set of coefficients $\{\gamma_{ni}\}_{n,i=1}^N$ can be obtained by ${{\mathcal O}}(N^3)$ operations from the coefficients $\beta_{nj}$, through a triangular process.
Indeed, the coefficients $\gamma_{ni}$ can be obtained by the following recursive formula
\begin{equation}
	\label{eq:gamma}
	\gamma_{ni} = \left\{\begin{array}{ll}
	\sum_{k=1}^{n-1} \beta_{nk}\gamma_{ki} & \mbox{if}\quad i<n,\\
	\beta_{nn} & \mbox{if}\quad i=n, \\
	0 &  \mbox{if}\quad i>n,
	\end{array}
	\right.
	\qquad \mbox{for all } n=1,\ldots,N.
\end{equation}

As mentioned above, the set of basis functions $\{\zeta_n\}_{n=1}^N$ leads to well-conditioned matrices, the precomputed quantities are nevertheless expressed in the basis $\{v(\bmu_n)\}_{n=1}^N$, but the change of basis functions is described by the coefficients $\{\gamma_{ni}\}_{n,i=1}^N$ as stated in equation \eqref{zetas}.

As an illustration assume that the reduced basis problem to be solved is about the approximation of the solution to a parametrized PDE, written in a variational form as : for a given parameter value $\bmu$, find $u(\bmu)\in\Wbb_\bmu$ such that
\[
	a(u(\bmu),v;\bmu) = f(v;\bmu),\qquad\forall v\in\Wbb_\bmu,
\]
where it is assumed, for the sake of simplicity that $a$ and $f$ are affine decomposable
\begin{equation}
	\label{eq:affdec}
	a(w,v;\bmu) = \sum_{q=1}^{Q_a} g_q(\bmu) \, a_q(w,v),\quad  f(v;\bmu) \sum_{q=1}^{Q_f} h_q(\bmu) \, f_q(v).
\end{equation}
The ultimate goal is to construct the well-conditioned stiffness matrix $a(\zeta_j,\zeta_i;\bmu)$ derived from the off-line pre-computation of the series  $A^q_{ij} = a_q	(v(\bmu_j),v(\bmu_i))$, for $1\le i,j\le N$. 
To do so, follow the next steps:
\vspace{5pt}
\begin{itemize}
\item[\mysymb{1}] 
	Compute the coefficients $\{\beta_{ni}\}_{n,i=1}^N$ following \eqref{eq:beta}.
\item[\mysymb{2}] 
	Compute the coefficients $\{\gamma_{ni}\}_{n,i=1}^N$ following \eqref{eq:gamma}.
\item[\mysymb{3}] 
	For each $q=1,\ldots,Q_a$, the matrix $a_q(\zeta_j,\zeta_i)$ can be computed by  $a_q(\zeta_j,\zeta_i) = \sum_{k,l=1}^N \gamma_{jk}\gamma_{il} A^q_{lk}$.
\item[\mysymb{4}] 
	Build the solution matrix by summing $a(\zeta_j,\zeta_i;\bmu)=\sum_{q=1}^{Q_a} g_q(\bmu) \, a_q(\zeta_j,\zeta_i)$.
\end{itemize}
\vspace{5pt}

Similarly, the evaluation of each component of the vector $f(\zeta_i;\bmu)$, $1\le i\le N$ can be performed in ${\mathcal O}(Q_fN)$ operations, and the inversion of the associated matricial problem is done in ${\mathcal O}(N^3)$ further operations.

\section{Overview of the framework}
\label{sec:BigPicture}

Before we start with the numerical tests, an abstract description of both the off-line and on-line parts of the locally adaptive greedy algorithm is given in the following box:

\mybox{
{\tt LocallyAdaptiveGreedy Off-line} \\
{\tt Input:} $N$\\
{\tt Output:} $K$, $\Wbb_K$, $\Sbb_K$
}{
\begin{itemize}
\item[] \hspace{-0.5cm}========= {\small {\tt Stage 1}} ========= 
\item[\tt 1.] Perform a classical greedy algorithm to select $N+1$ basis functions.
\item[] \hspace{-0.5cm}========= {\small {\tt Stage 2}} ========= 
\item[\tt 2.] For each $\bmu\in\Xitrial$, do (parameter space search)
\begin{itemize}
\item[{\tt a)}] Compute the error estimate $\eta(\bmu,\Wbb_\bmu)$.
\item[{\tt b)}] Update the Hessian matrix ${H}_{ij} (\bmu)$.
\end{itemize}
\item[\tt 4.] Enrich the set of basis functions ({\tt EnrichmentLoop}, Section \ref{ssec:SamplePoints}).
\item[\tt 5.] Create a new training set $\Xitrial$ (Section \ref{sec:hessian}).
\item[\tt 6.] Go to {\tt 2.} until tolerance $\tt tol$ is achieved.
\end{itemize}
}{
Algorithm 3: Locally adaptive greedy algorithm (Off-line).
}
\vspace*{11pt}
\mybox{
{\tt LocallyAdaptiveGreedy On-line} \\
{\tt Input:} $K$, $\Wbb_K$, $\Sbb_K$, $\bmu$\\
{\tt Output:} $P_\bmu(v(\bmu))$
}{
\begin{itemize}
\item[\tt 1.] Find the approximation space $\Wbb_\bmu$:
\begin{itemize}
\item[{\tt a)}] Compute the distance $d(\bmu,\tilde\bmu)$ for each $\tilde\bmu\in\Sbb_K$ (Section \ref{ssec:Hessian}).
\item[{\tt b)}] Choose the $N$ closest $\bar\bmu\in\Sbb_K$ to form $\Sbb_{\bmu}$ and $\Wbb_\bmu =\myspan{v(\tilde\bmu)\,|\, \tilde\bmu\in\Sbb_{\bmu}}$ (Section \ref{ssec:LAS}).
\end{itemize}
\item[\tt 2.] Orthonormalize the basis functions $\{v(\tilde\bmu)\,|\, \tilde\bmu\in\Sbb_{\bmu}\}$ to build a new basis $\{\xi_1,\ldots,\xi_N\}$ (Section \ref{sec:online}).
\item[\tt 3.] Compute the projection $P_\bmu(v(\bmu))$ based on the basis  $\{\xi_1,\ldots,\xi_N\}$.
\end{itemize}
}{
Algorithm 4: Locally adaptive greedy algorithm (On-line).
}


\section{Numerical results using the $L^2$-projection}
\label{sec:L2}
Before we pass to a full Reduced Basis model based on a parametrized Galerkin-projection, we use a parametrized family of explicitly given functions combined with the $L^2$-projection to define the projection operator $P_N$ in order to test our framework. 
Reminding Cea's lemma which is a consequence of the Galerkin approach, studying the best-approximation is sufficient to test the method of locally adaptive approximation spaces.
The error is computed in an exact manner and no {\em a posteriori} estimation is used.

Further, in a first step, Section \ref{ssec:FTS}, we only test the feature of local anisotropic approximation balls without adapting the training set using the distance function. That is we use the same fixed $\Xitrial$ throughout all steps of the algorithm.
In a second step, Section \ref{ssec:ATS}, we test the entire algorithm.

For all numerical tests, $\Omega$ and $\dom$ are discretized by a regular lattice of $75\times 75$ points and the target tolerance is set to $10^{-4}$ if not otherwise mentioned.

\subsection{Numerical results with fixed training set}
\label{ssec:FTS}
\subsubsection*{Test 1}
\begin{figure}[!ht]
  \captionsetup[subfigure]{labelformat=empty}
   \captionsetup[subfigure]{labelformat=empty,width=\two}
 \centering
  \subfloat[Comparison of different approaches.]{
  \includegraphics[width=\two]{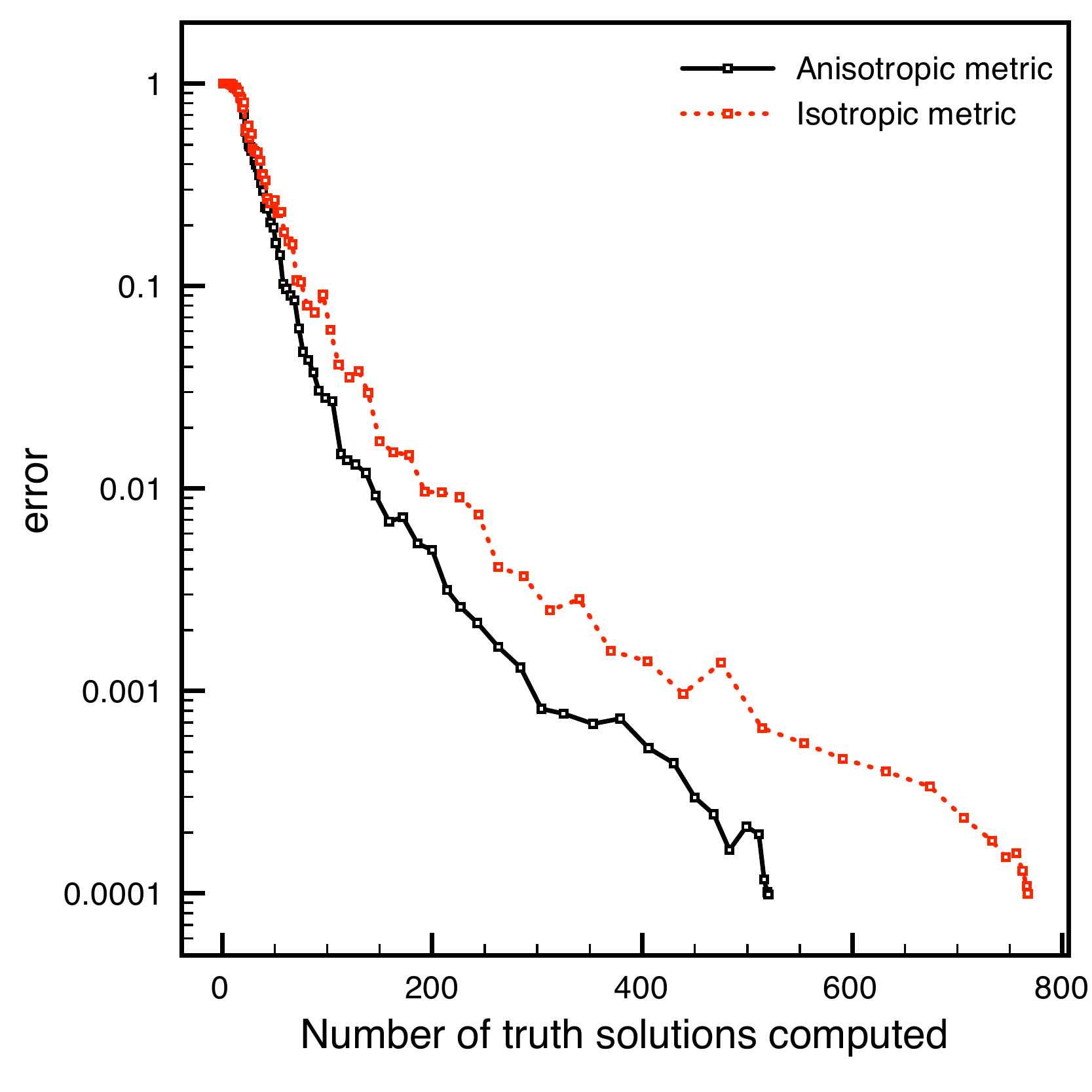}}
  \subfloat[Comparison under variation of $N$.]{
  \includegraphics[width=\two]{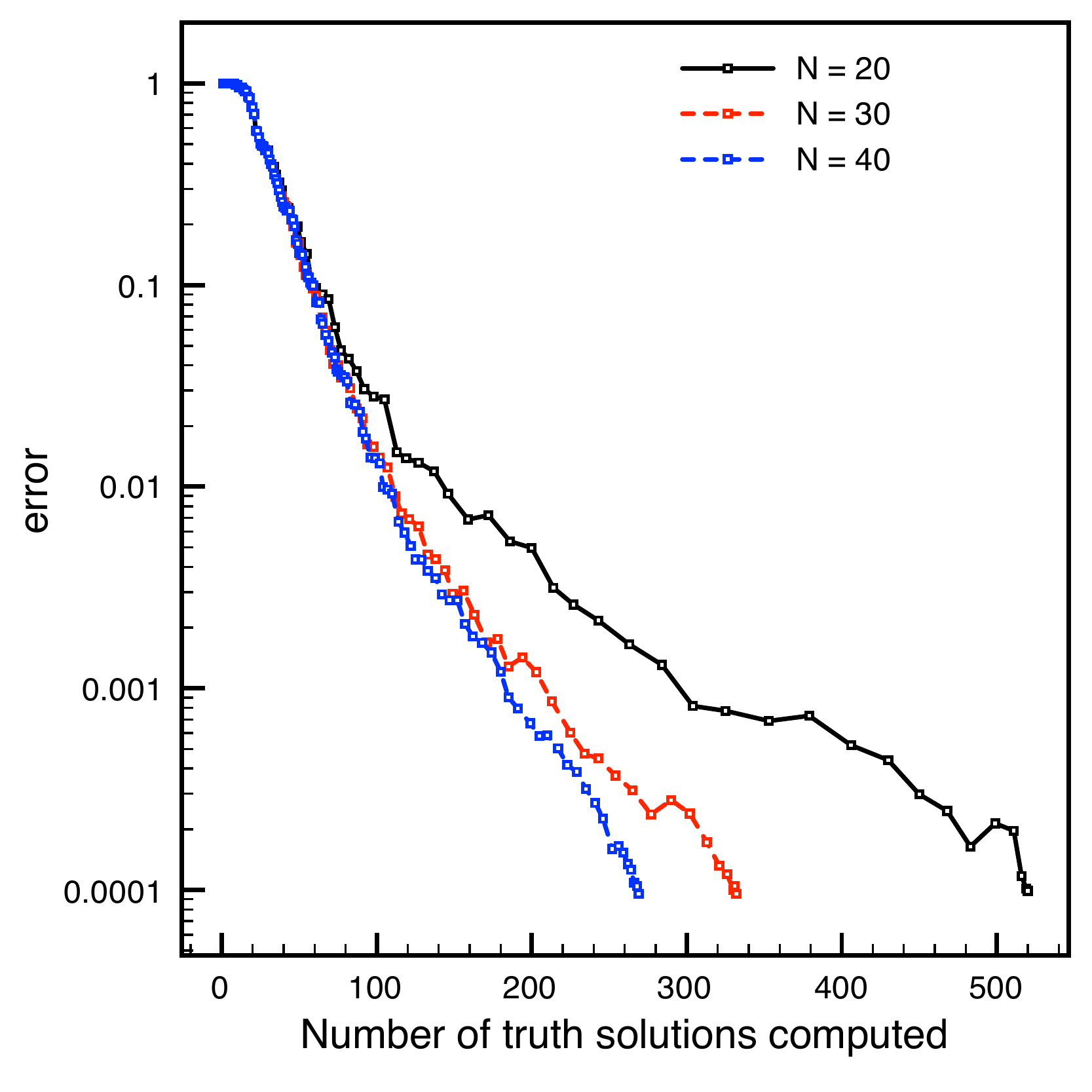}
  }
  \caption{Test 1:  Accuracy with respect to the number of truth solutions to be computed in comparison with the isotropic approach for $N=20$ (left) and the for different values of $N=20,30,40$ (right).}
  \label{fig:Ex1_conv}
\end{figure}

We start with presenting a numerical example to illustrate the benefit of the local anisotropic approximation spaces. Consider the function
\begin{eqnarray*}
	f_1(\x;\bmu) = \exp\left[
	\frac{- (x_1-0.1(\mu_1-\mu_2))^2}{0.01} 
	 - 	\frac{(x_2-(\mu_1+\mu_2))^2}{0.01}\right],& \\
	\x\in \Omega=(-1,1)^2,&\,\bmu \in \dom=[-0.5,0.5]^2
\end{eqnarray*}
that exhibits a constant anisotropy of parameters over the whole parameter space. The dependency in the $(\mu_1+\mu_2)$-direction is ten times stronger than in the $(\mu_1-\mu_2)$-direction. 
Figure \ref{fig:Ex1_conv} (left) illustrates the evolution of the number of basis functions that are necessary to be computed versus the achieved accuracy in the $L^\infty$-norm for $N=20$ for the anisotropic approach compared to the same algorithm but using isotropic local approximation spaces as comparison. 
One can observe that using the anisotropic approach is beneficial in terms of the number of truth solutions to be computed at the off-line stage.
Figure \ref{fig:Ex1_conv} (right) presents the same quantity but for varying $N=20,30,40$. Not surprisingly the number of needed basis functions decreases while increasing $N$ and more and more the exponential convergence establishes as is manifested for the standard greedy algorithm in this example. We can also observe that for a lower $N$ more basis functions can be included per iteration. 

\begin{figure}[!ht]
  \captionsetup[subfigure]{labelformat=empty,width=\three}
  \centering
\subfloat[Local approximation spaces.]{
\begin{minipage}[t]{\three}
\includegraphics[Trim=\trimval,clip,width=\textwidth]{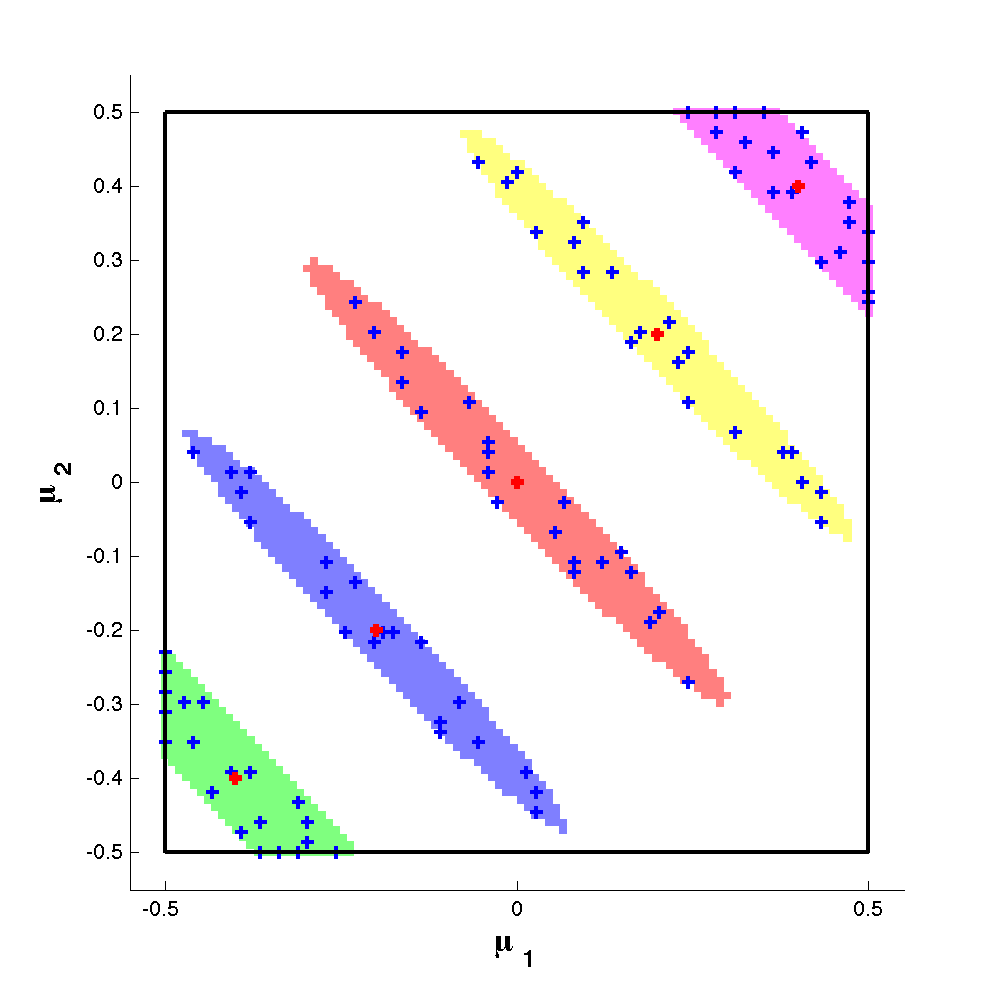}
\\
\includegraphics[Trim=\trimval,clip,width=\textwidth]{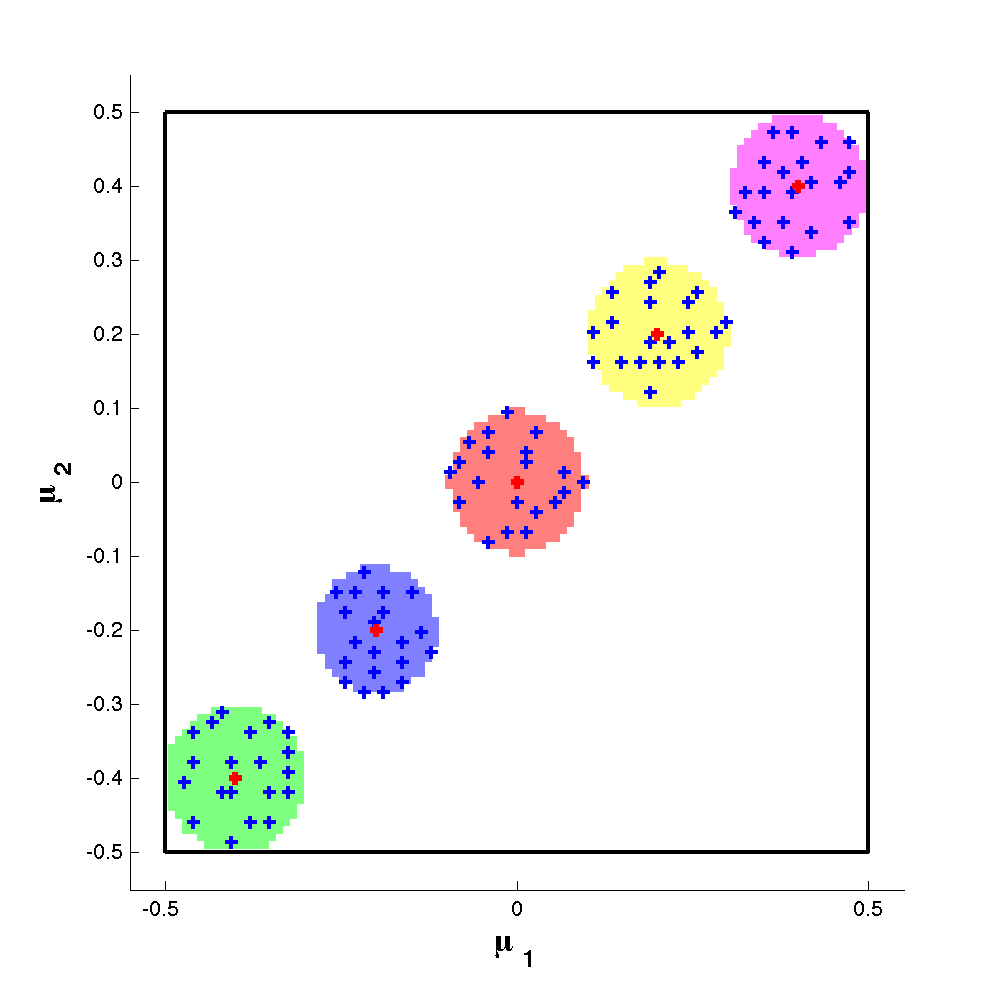}
\end{minipage}
}
\subfloat[Radius.]{
\begin{minipage}[t]{\three}
\includegraphics[Trim=\trimvalr,clip,width=\textwidth]{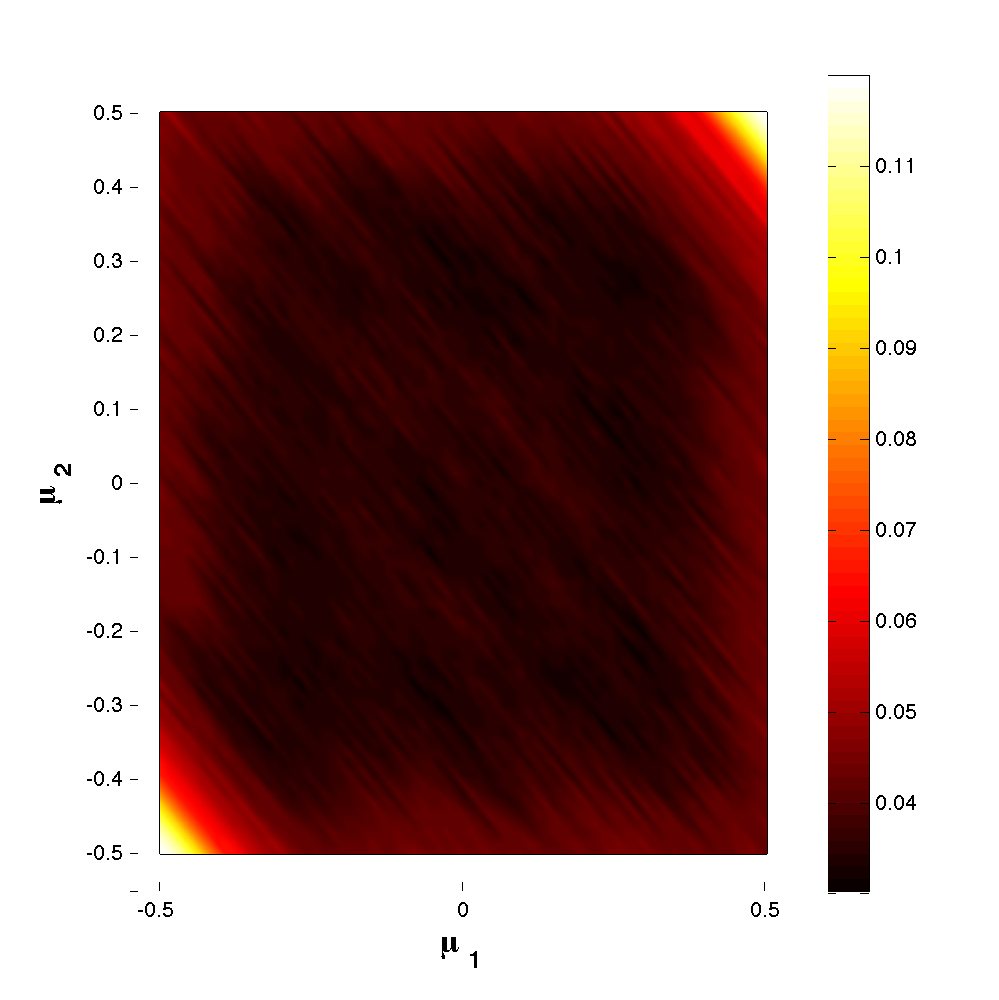}
\\
\includegraphics[Trim=\trimvalr,clip,width=\textwidth]{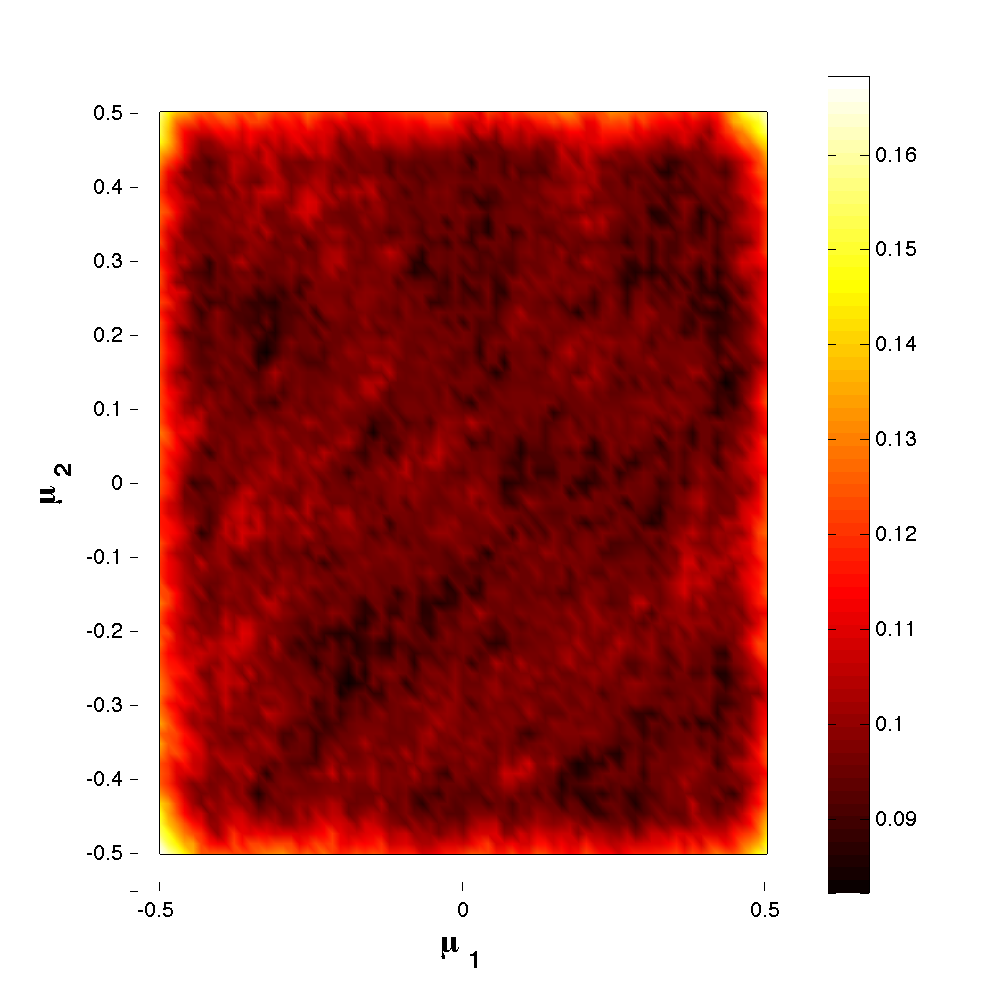}
\end{minipage}
}
\subfloat[Sample points.]{
\begin{minipage}[t]{\three}
\includegraphics[Trim=\trimval,clip,width=\textwidth]{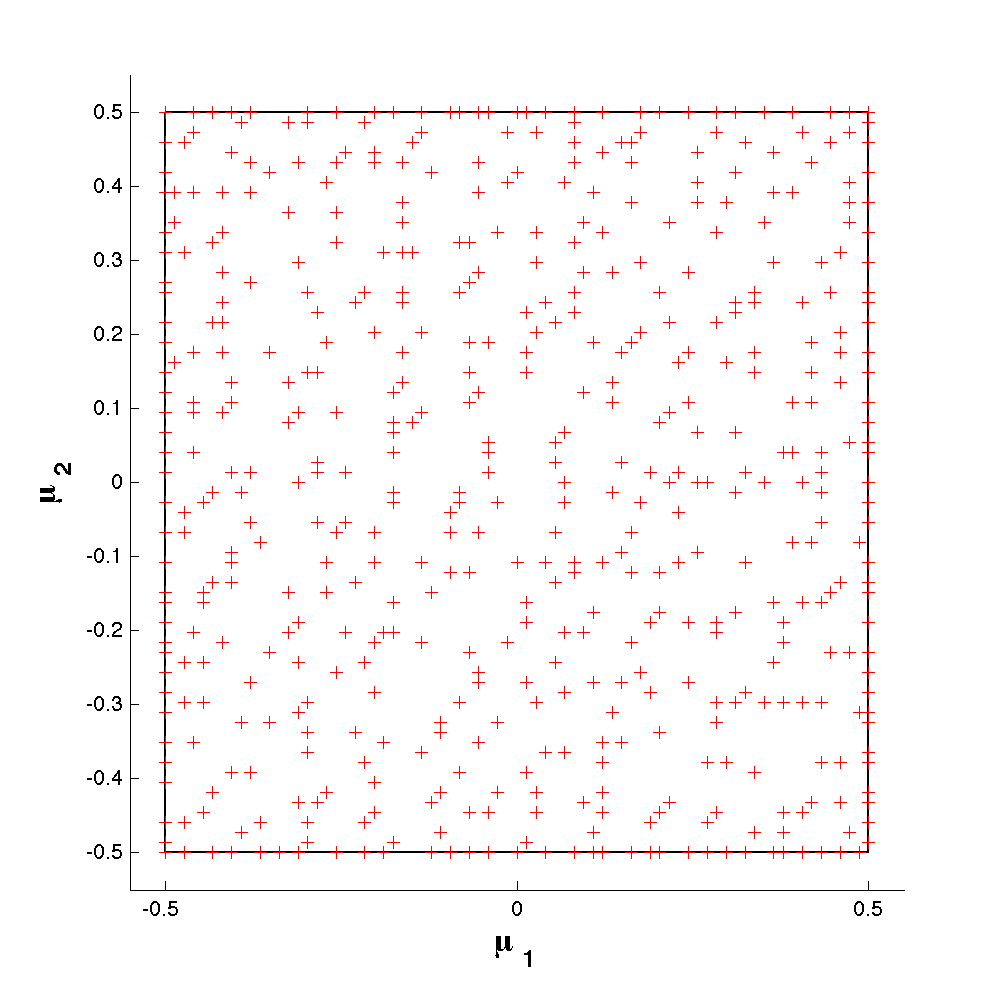}
\\
\includegraphics[Trim=\trimval,clip,width=\textwidth]{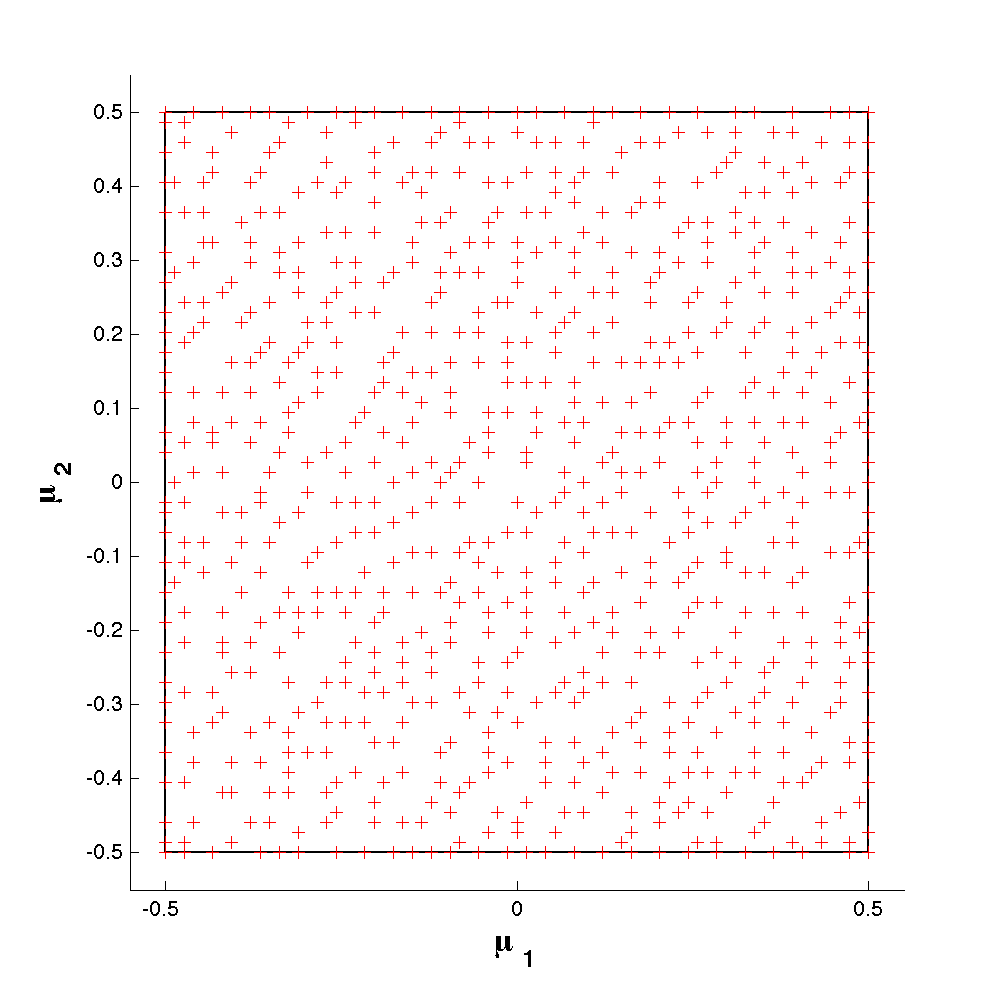}
\end{minipage}
}
\caption{
Test 1: Local approximation spaces for selected parameter values (left), radius as a function of the parameters (middle) and sample points (right) for anisotropic approach (top) and the isotropic version (bottom, as comparison) with $N=20$.
}
\label{fig:Test1}
\end{figure}

Figure \ref{fig:Test1} illustrates the local approximation spaces expressed in $\mathbb P$, the radius as a function of the parameter and the sample set $\Sbb_K$ for the anisotropic and the isotropic approach.

\subsubsection*{Test 2}
The previous example is in some sense idealized since the anisotropy is not varying throughout the parameter domain. 
In this regard, the next example
\begin{eqnarray*}
	f_2(\x;\bmu) = \exp
	\left[
	-
	\frac{
	(x_1-(\mu_1^2+\mu_2^2))^2}{0.01
	}
	-
	\frac{
	(x_2-(\mu_1^2+\mu_2^2))^2}{0.01
	}
	\right],& \\
	\x\in \Omega=(-1,1)^2,&\,\bmu \in \dom=[-0.5,0.5]^2
\end{eqnarray*}
is more interesting. It presents a family of parametrized functions where the functions (as functions of $\x$) are constant along concentric circles around the origin in parameter space. As an example, all four corners of the parameter space define the same function. 

\begin{figure}[!ht]
  \captionsetup[subfigure]{labelformat=empty,width=\two}
  \centering
\subfloat[Anisotropy vs. Isotropy]{
\includegraphics[width=\two]{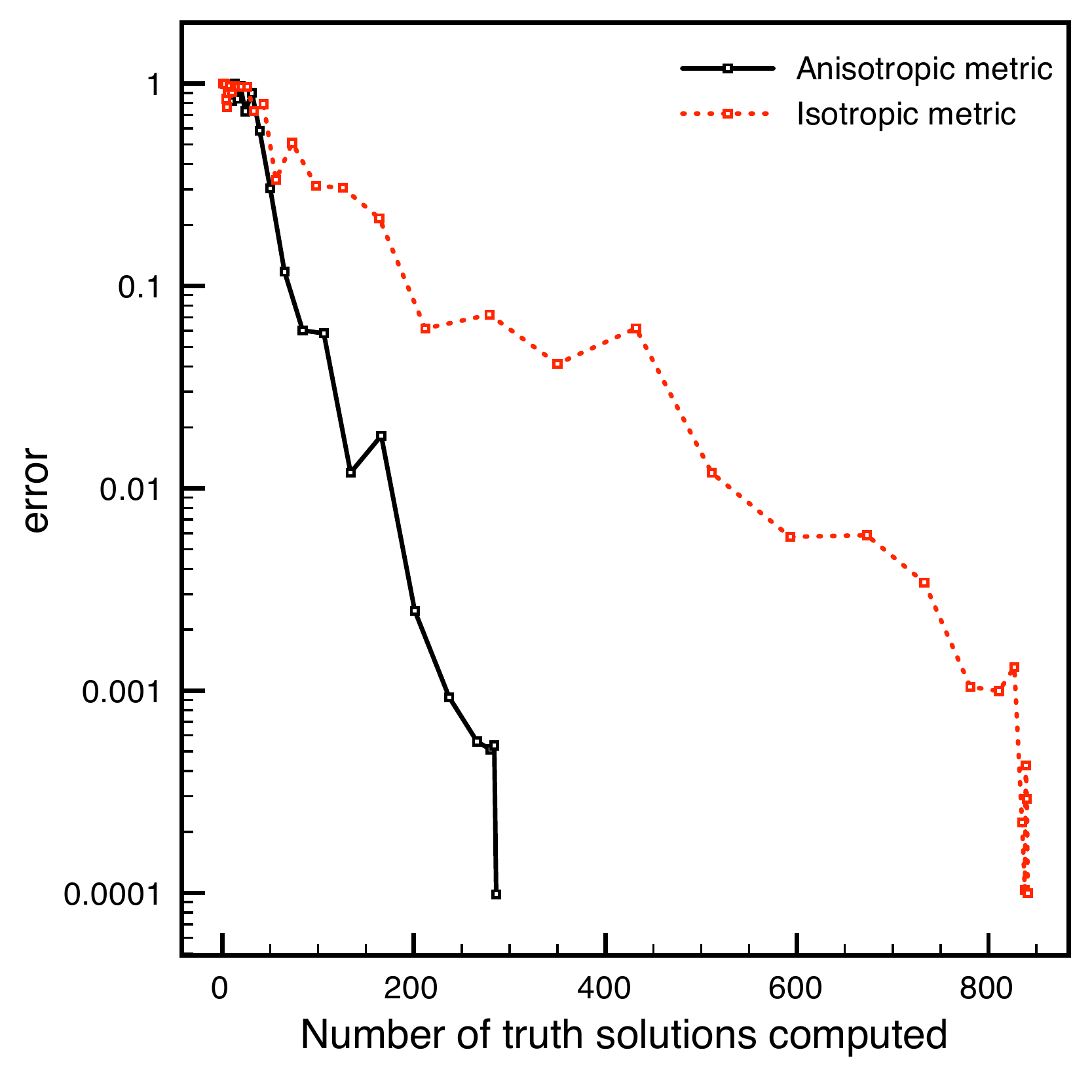}
}
\subfloat[Different tolerances]{
\includegraphics[width=\two]{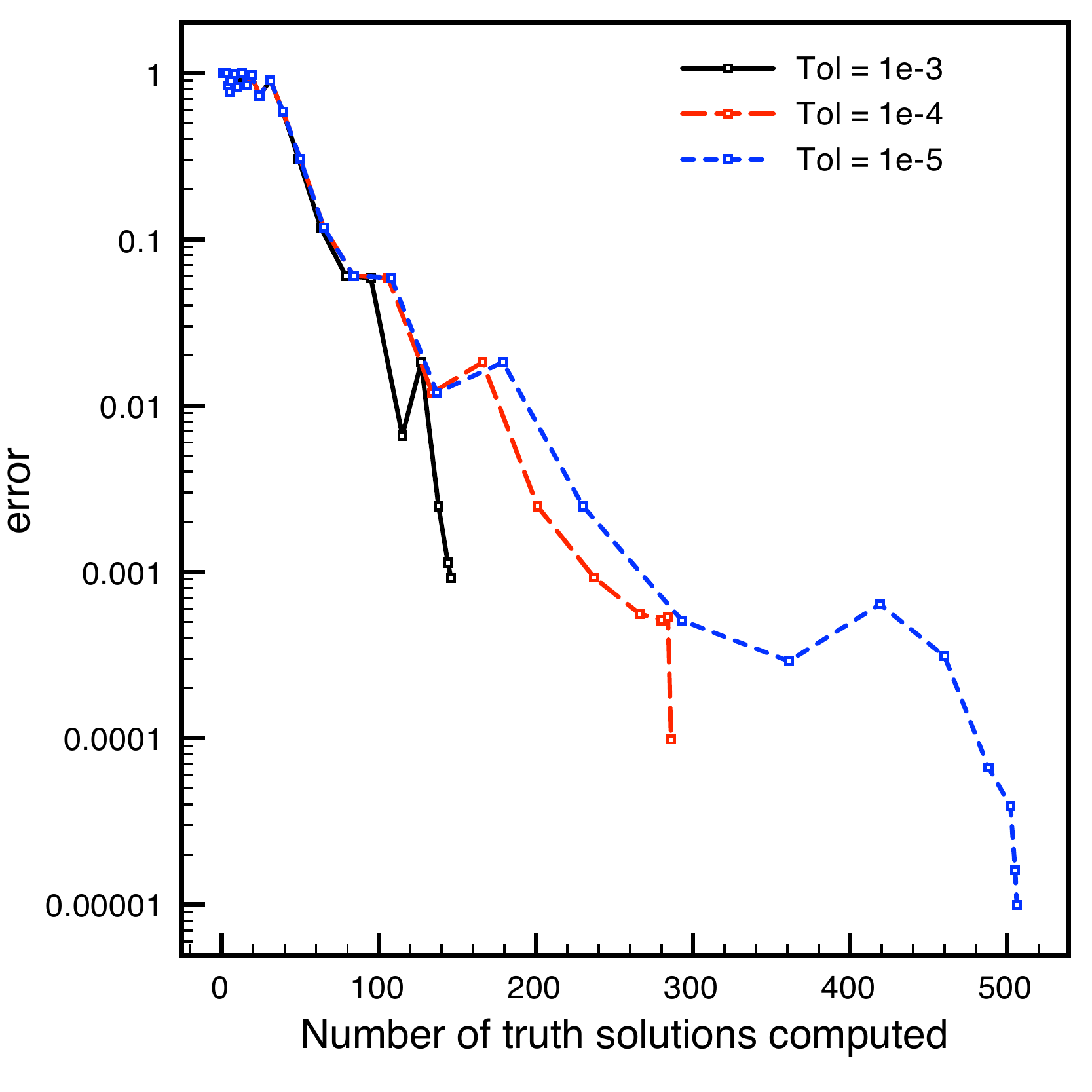}
}
\caption{Test 2: Accuracy with respect to the number of truth solutions to be computed in comparison with the isotropic approach (left) and the for different end tolerances (right). In both cases there holds $N=5$.}
\label{fig:Ex2_conv}
\end{figure}

\begin{figure}[!ht]
  \captionsetup[subfigure]{labelformat=empty,width=\three}
  \centering
\subfloat[Local approximation spaces.]{
\begin{minipage}[t]{\three}
\includegraphics[Trim=\trimval,clip,width=\textwidth]{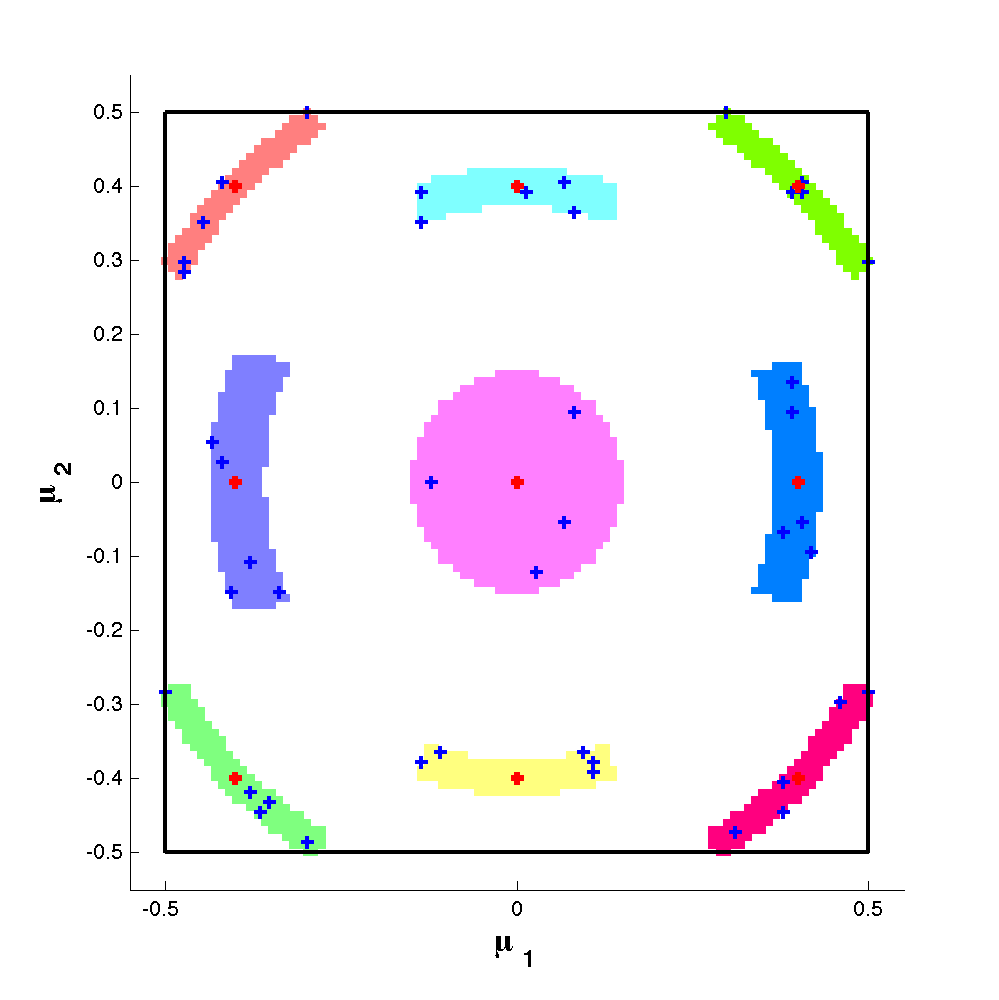}
\\
\includegraphics[Trim=\trimval,clip,width=\textwidth]{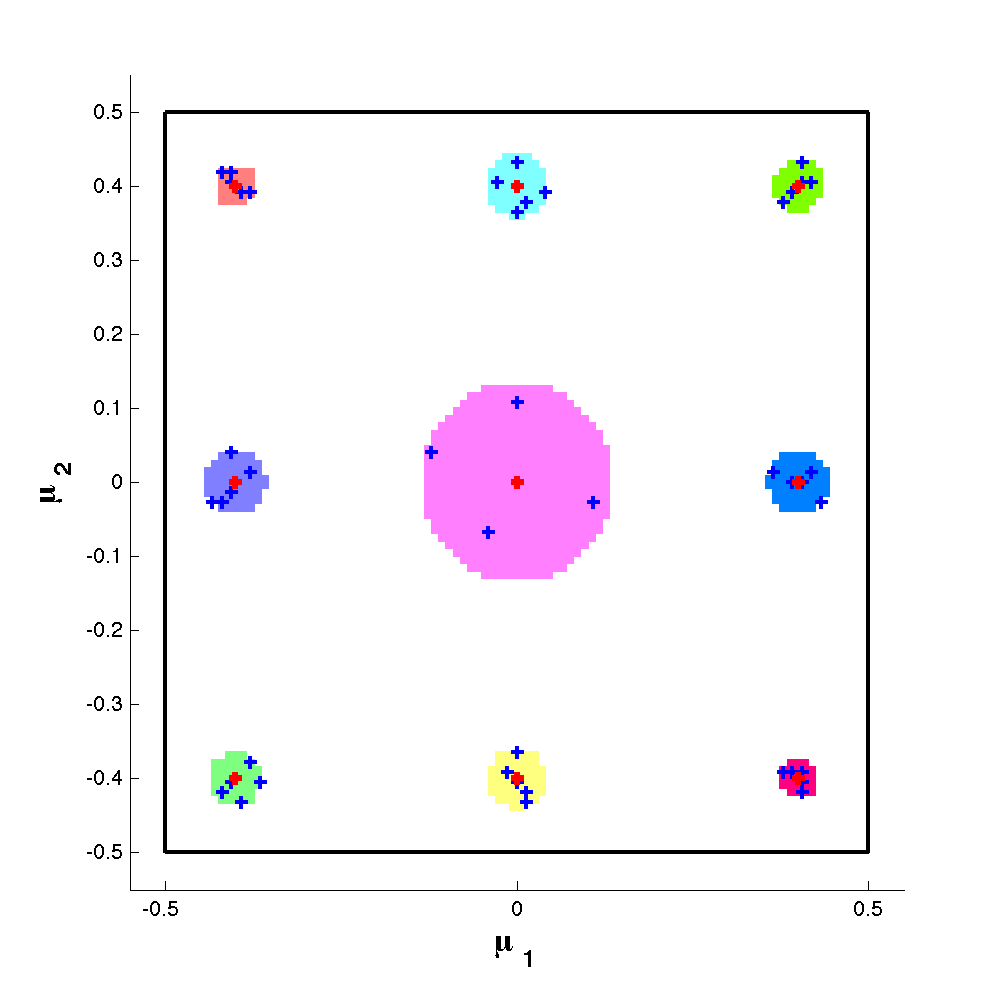}
\end{minipage}
}
\subfloat[Radius.]{
\begin{minipage}[t]{\three}
\includegraphics[Trim=\trimvalr,clip,width=\textwidth]{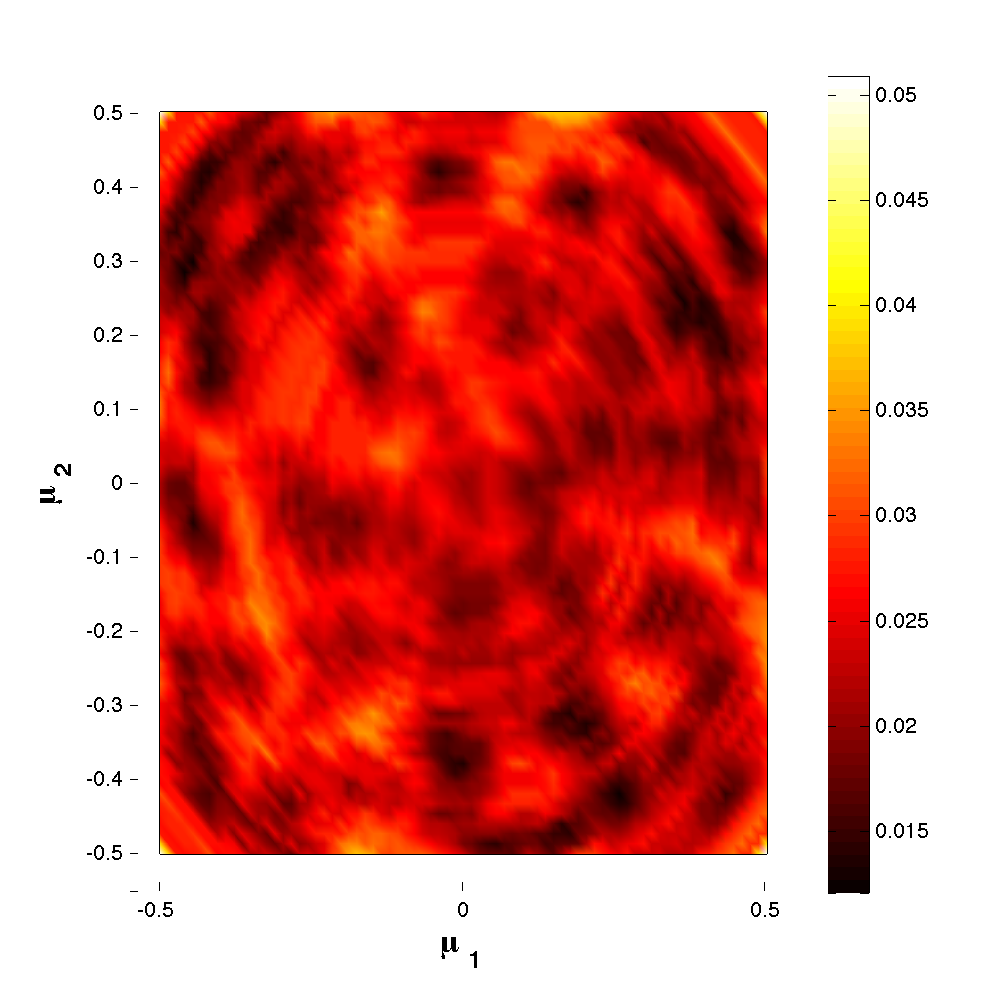}
\\
\includegraphics[Trim=\trimvalr,clip,width=\textwidth]{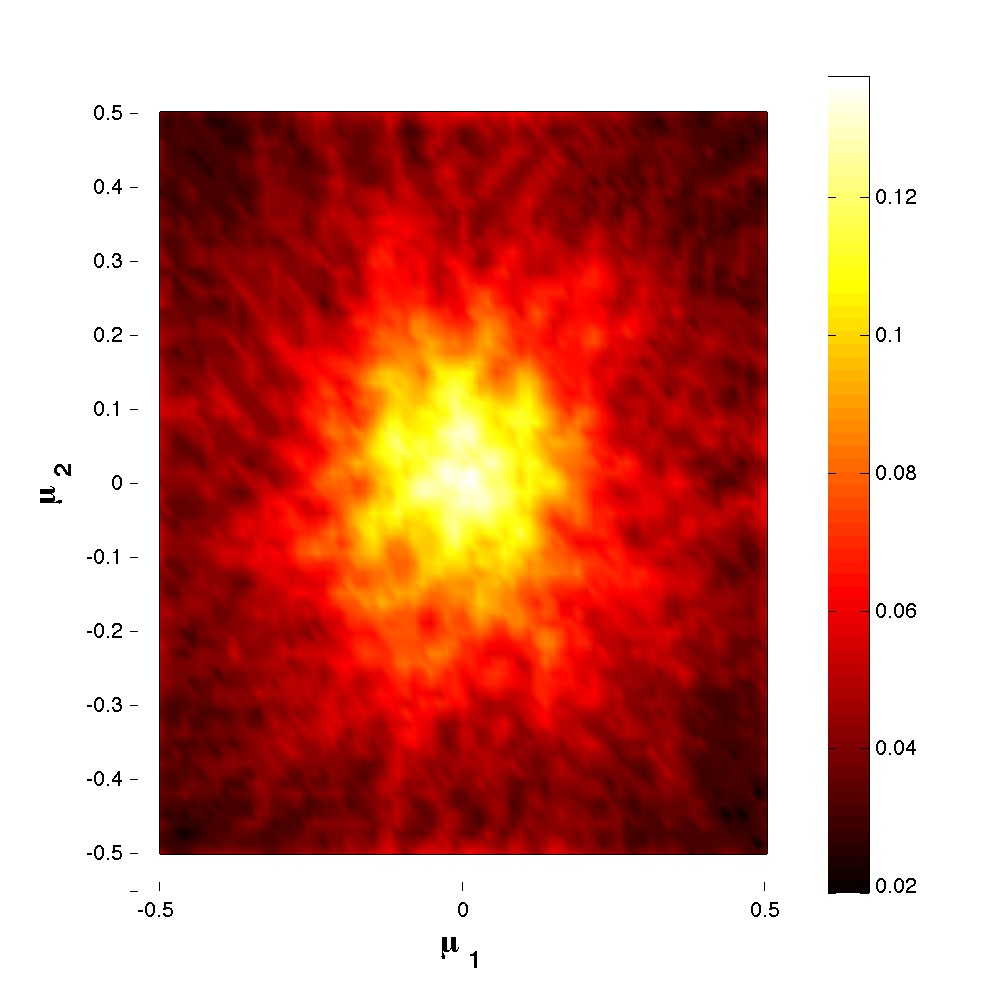}
\end{minipage}
}
\subfloat[Sample points.]{
\begin{minipage}[t]{\three}
\includegraphics[Trim=\trimval,clip,width=\textwidth]{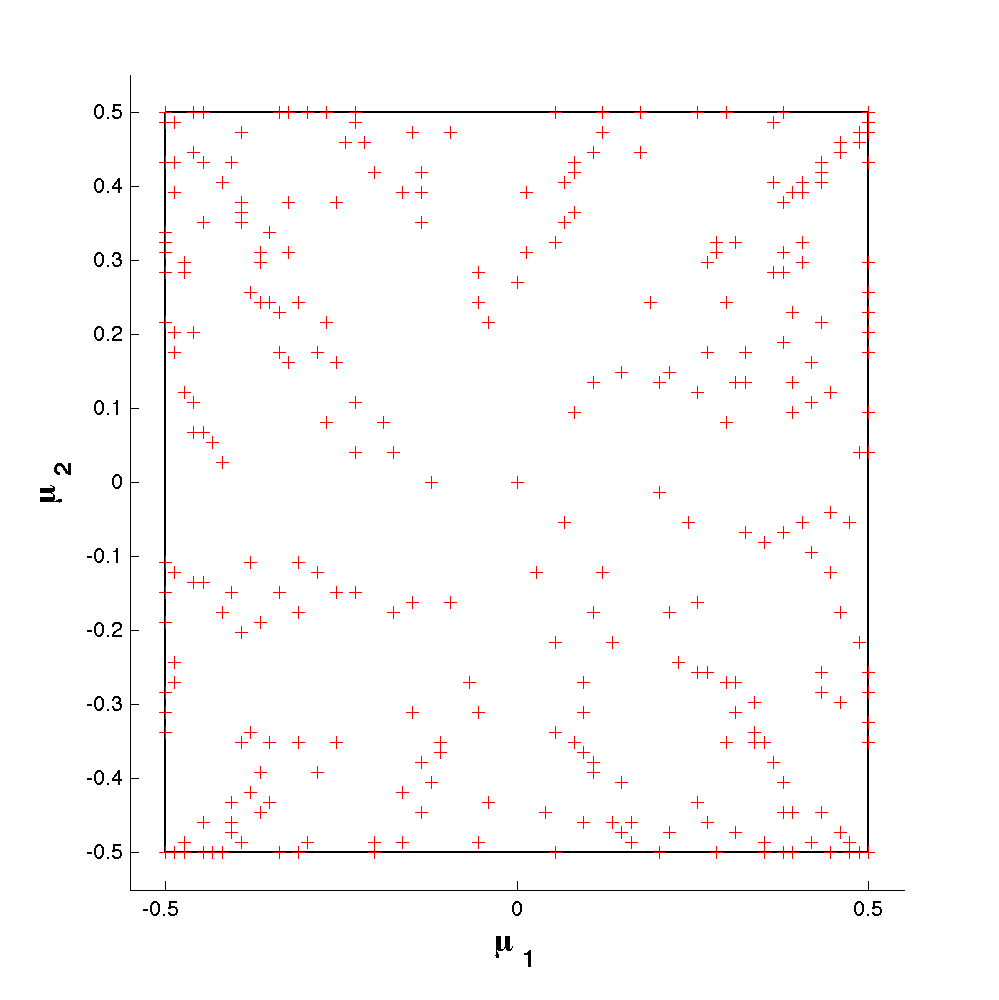}
\\
\includegraphics[Trim=\trimval,clip,width=\textwidth]{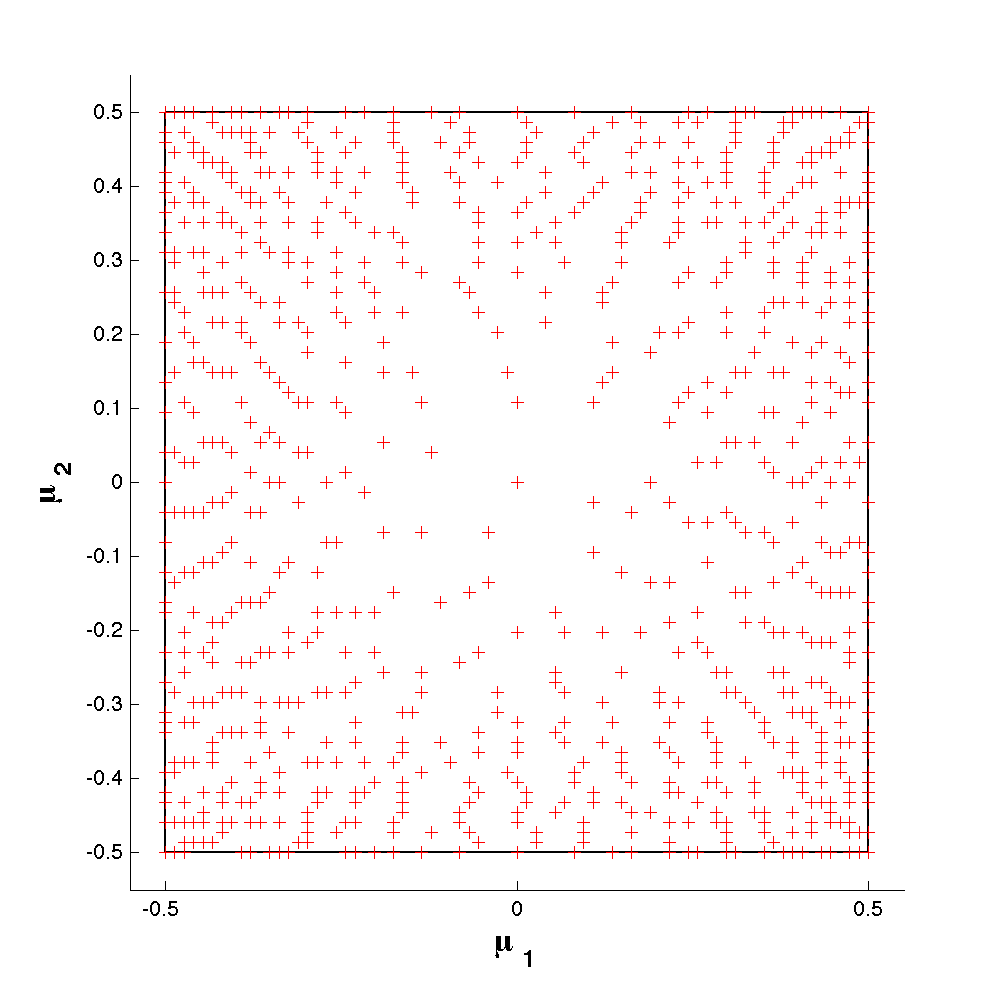}
\end{minipage}
}
\caption{
Test 2: Local approximation spaces for selected parameter values (left), radius as a function of the parameters (middle) and sample points (right) for the presented approach (top) in comparison wit the isotropic version (bottom) for $N=5$.
}
\label{fig:Test2}
\end{figure}

We compare again the number of required truth solutions for the anisotropic approach compared with the presented algorithm but using isotropic approximation spaces, which is presented in Figure \ref{fig:Ex2_conv} (left) for $N=5$. 
We observe that 60\% of the number of truth solutions to be computed can be saved by constructing the problem-dependent approximation spaces.
The corresponding local approximation spaces, radii and sample points $\Sbb$ are illustrated in Figure \ref{fig:Test2}.
Figure \ref{fig:Ex2_conv} (right) also presents the number of required truth solutions for $N=5$ but for different tolerance levels. 
We observe that starting from a critical value $N_c$ of computed truth solutions, that depends on the tolerance, the curves show different behavior. 
During the iterations before $N_c$ the complete parameter region has to be scanned because the tolerance is achieved nowhere in the parameter space.
During the iterations after $N_c$ however the tolerance is achieved in an increasing region of the parameter space and this part needs no longer to be not scanned.
The enrichment of new truth solutions is then limited in space and the final phase of the greedy algorithm converges faster.

\begin{figure}[!ht]
  \captionsetup[subfigure]{labelformat=empty,width=\three}
  \centering
\subfloat[Tol=$10^{-3}$]{
\begin{minipage}[t]{\three}
\includegraphics[Trim=\trimval,clip,width=\textwidth]{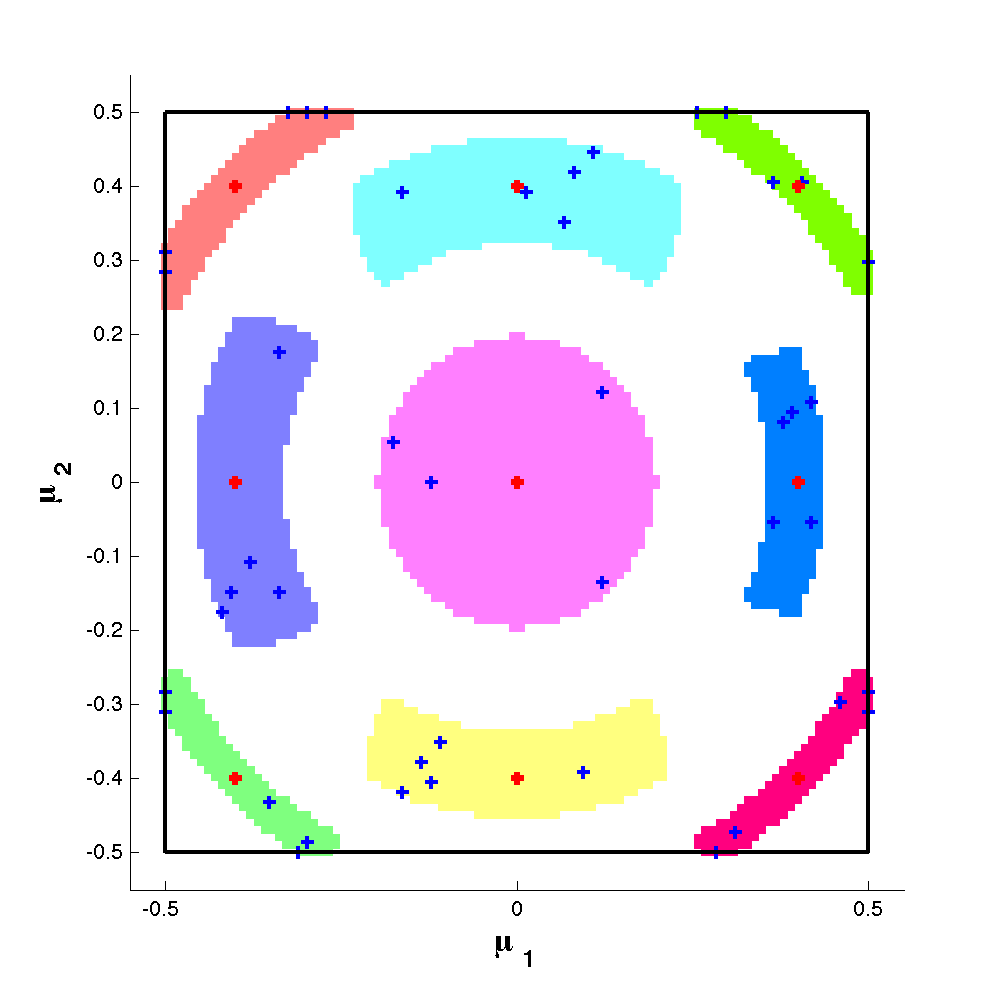}
\end{minipage}
}
\subfloat[Tol=$10^{-4}$]{
\begin{minipage}[t]{\three}
\includegraphics[Trim=\trimval,clip,width=\textwidth]{pics/LAS_TEST2_N5_TOL4.png}
\end{minipage}
}
\subfloat[Tol=$10^{-5}$]{
\begin{minipage}[t]{\three}
\includegraphics[Trim=\trimval,clip,width=\textwidth]{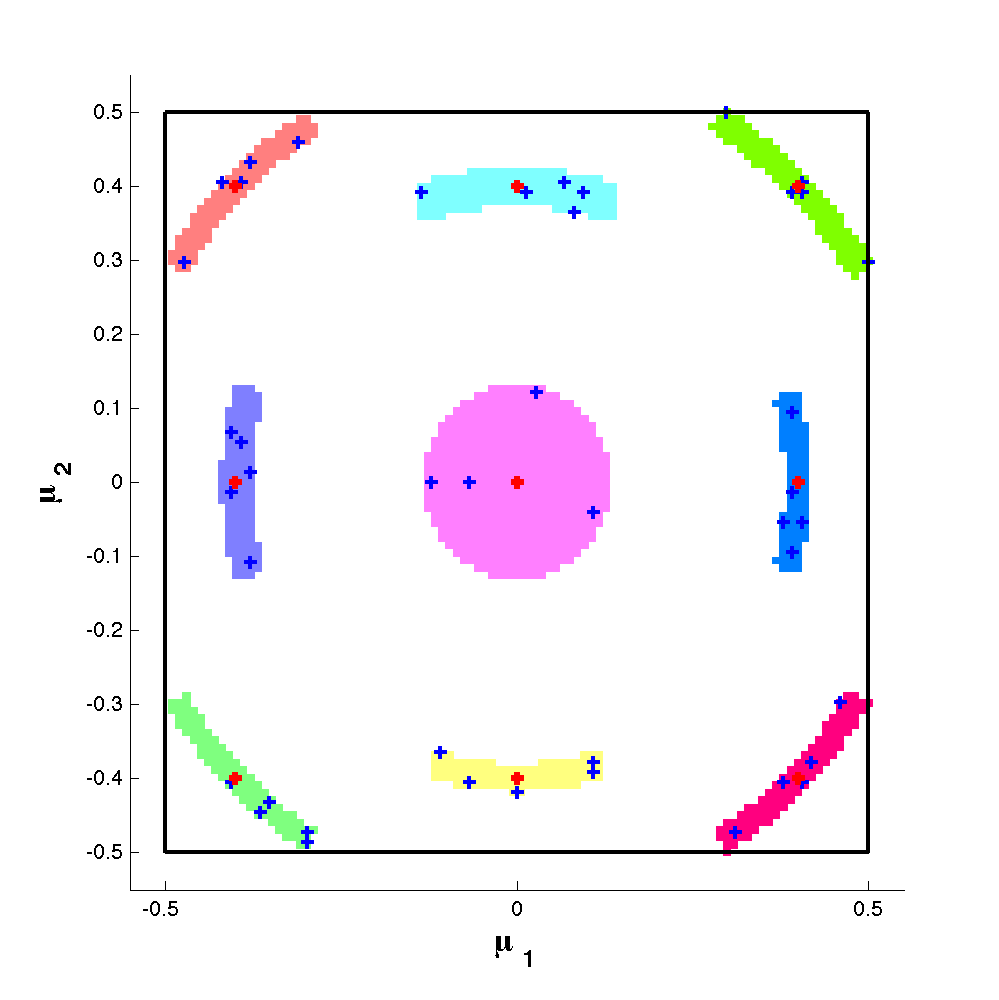}
\end{minipage}
}
\caption{
Test 2: Local approximation spaces for different tolerance levels.}
\label{fig:Ex2_LAS_tol}
\end{figure}

Additionally we present in Figure \ref{fig:Ex2_LAS_tol} the local approximation spaces for the different tolerance levels. 
One can clearly see how the approximation spaces take into account the local geometry that also allow non-convex approximation spaces.


\subsubsection*{Test 3}
\begin{figure}[!ht]
  \captionsetup[subfigure]{labelformat=empty,width=\two}
  \centering
\subfloat[]{
\includegraphics[trim =0mm 0mm 0mm 0mm, clip, width=\two]{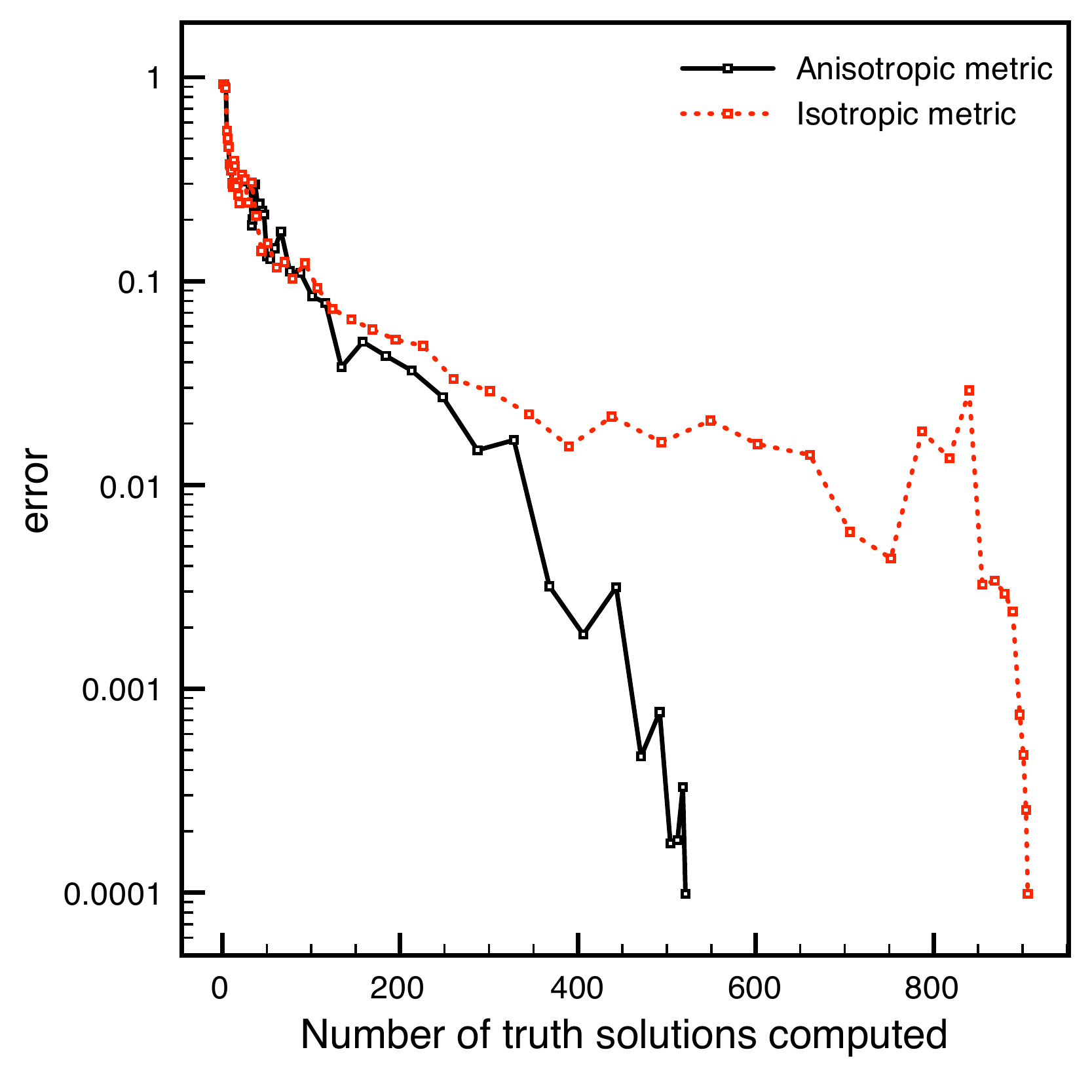}
}
\vspace{-0.8cm}
\caption{Test 3: Accuracy with respect to the number of truth solutions to be computed in comparison with the isotropic approach for $N=10$.}
\label{fig:Ex3_conv}
\end{figure}

The next example
\begin{eqnarray*}
	f_3(\x;\bmu) = \exp
	\left[
	-
	\frac{
	(x_1-(\mu_1+3\mu_2))^2}{0.1+5|\mu_1+3\mu_2|
	}
	-
	\frac{
	(x_2-(3\mu_1-\mu_2))^2}{0.1+5|3\mu_1-\mu_2|
	}
	\right],& \\
	\x\in \Omega=(-1,1)^2,&\,\bmu \in \dom=[-0.5,0.5]^2
\end{eqnarray*}
is interesting in the sense that it presents an almost singularity in parameter space at the origin. The performance of the anisotropic approach, compared to the isotropic version, is presented in Figure \ref{fig:Ex3_conv}. Again, the anisotropic approaches outperform the isotropic one and about a 56\% of computations of truth solutions can be saved. 

\begin{figure}[!ht]
  \captionsetup[subfigure]{labelformat=empty,width=\three}
  \centering
\subfloat[Local approximation spaces]{
\begin{minipage}[t]{\three}
\includegraphics[Trim=\trimval,clip,width=\textwidth]{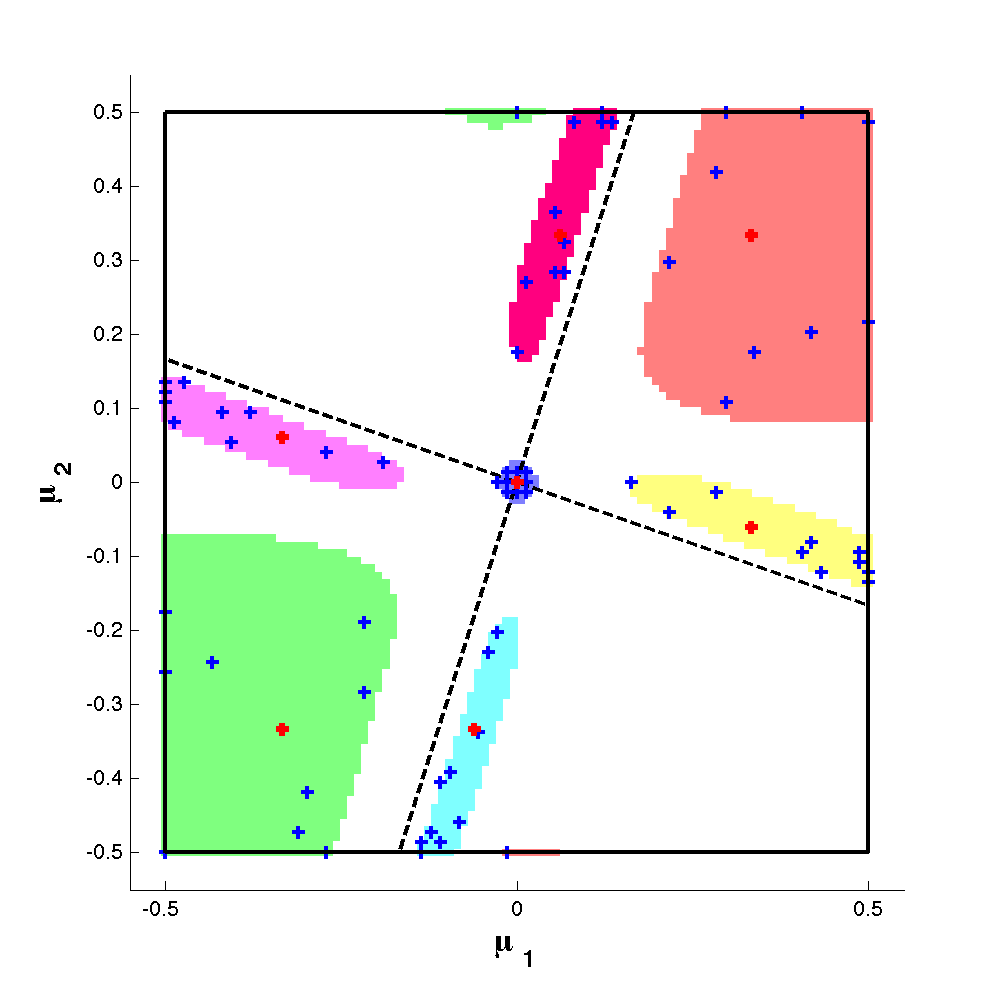}
\\
\includegraphics[Trim=\trimval,clip,width=\textwidth]{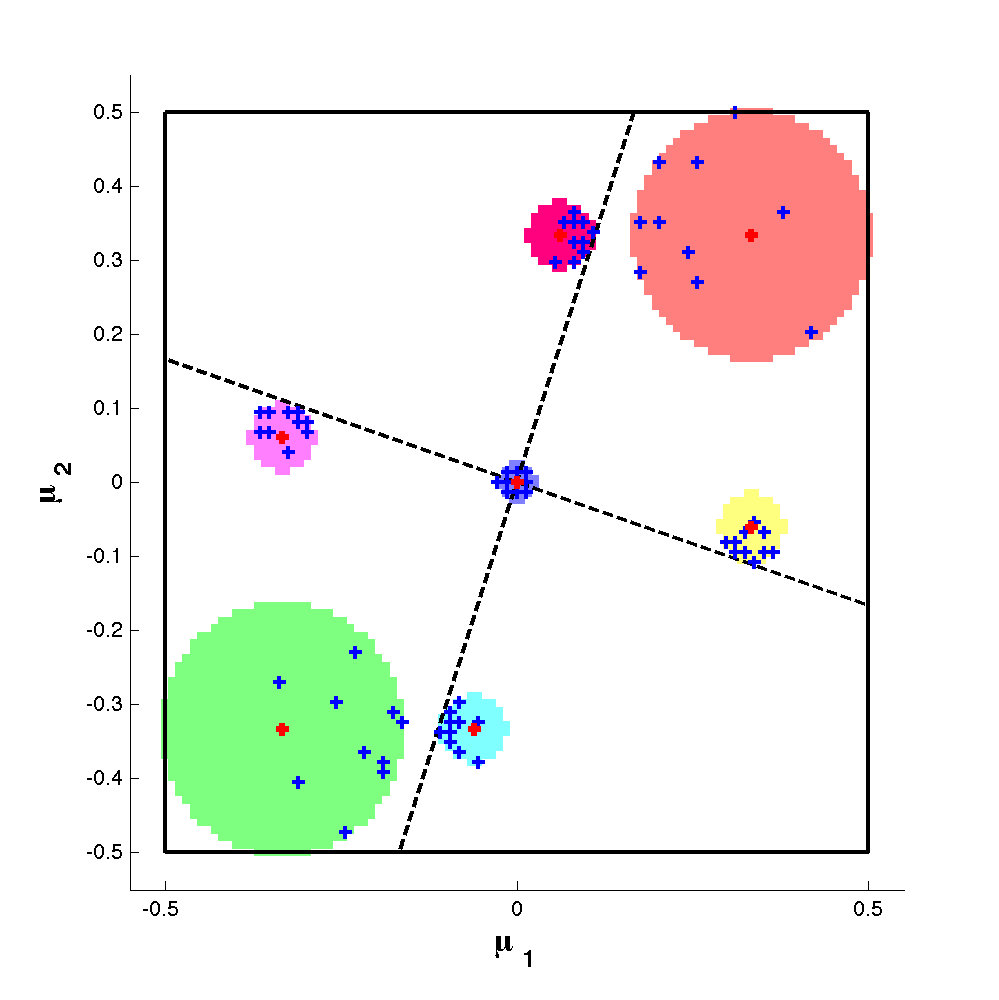}
\end{minipage}
}
\subfloat[Radius]{
\begin{minipage}[t]{\three}
\includegraphics[Trim=\trimvalr,clip,width=\textwidth]{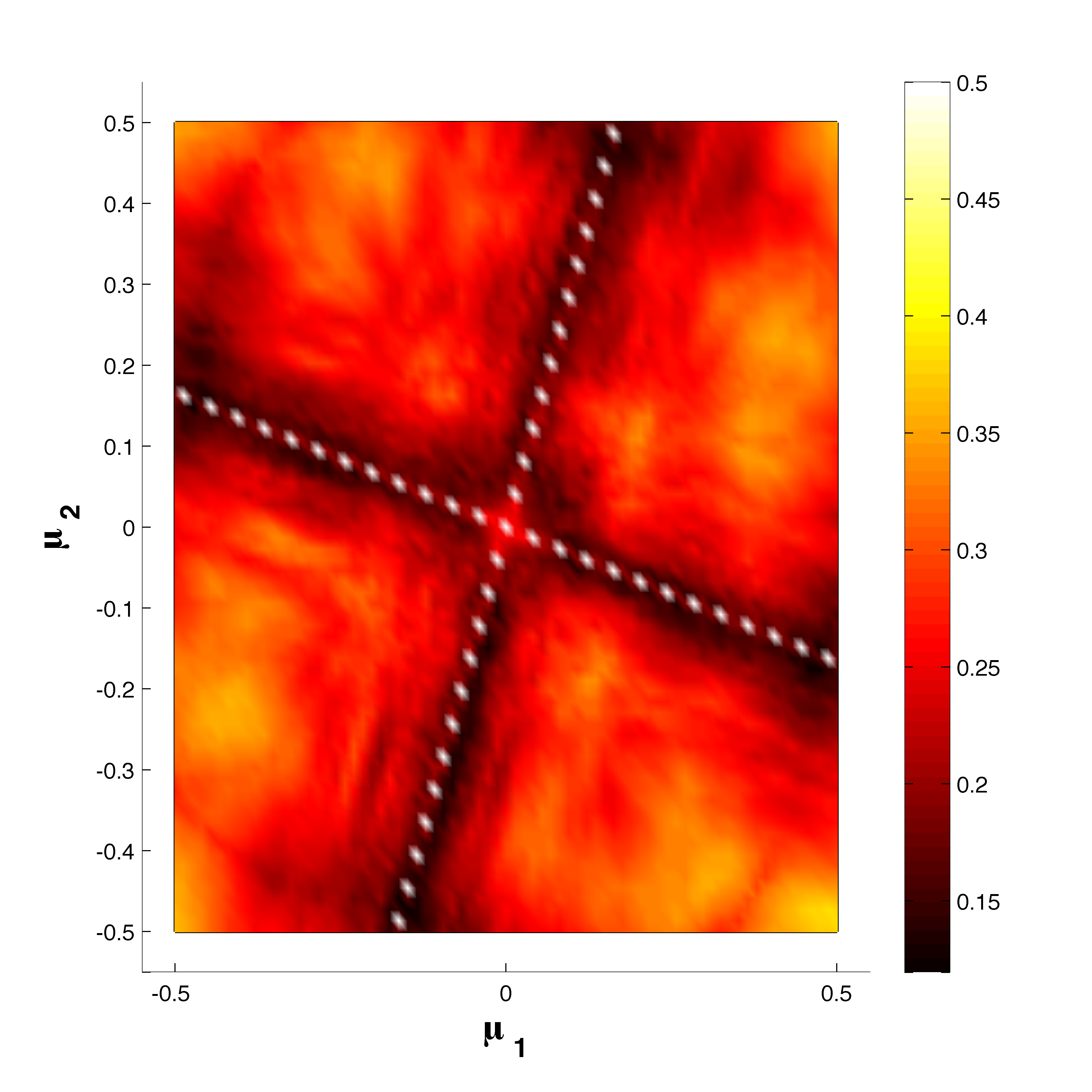}
\\
\includegraphics[Trim=\trimvalr,clip,width=\textwidth]{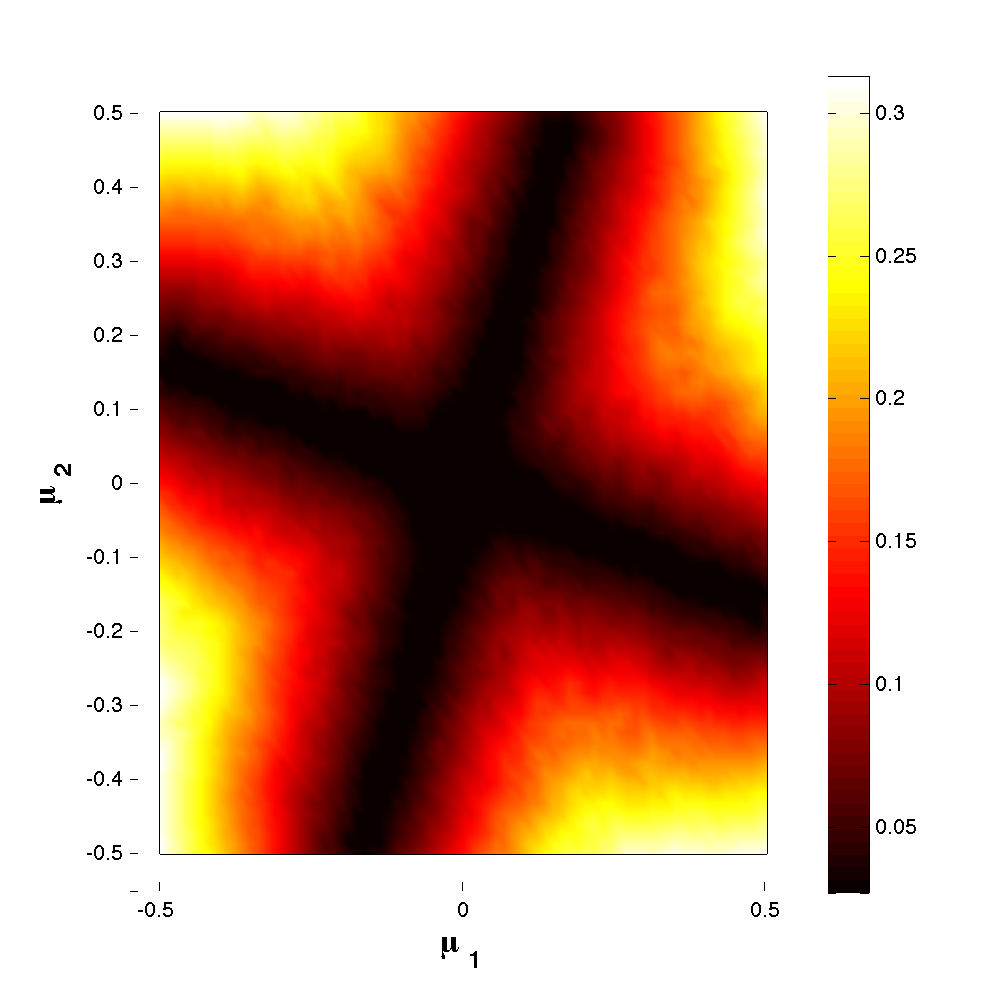}
\end{minipage}
}
\subfloat[Sample points]{
\begin{minipage}[t]{\three}
\includegraphics[Trim=\trimval,clip,width=\textwidth]{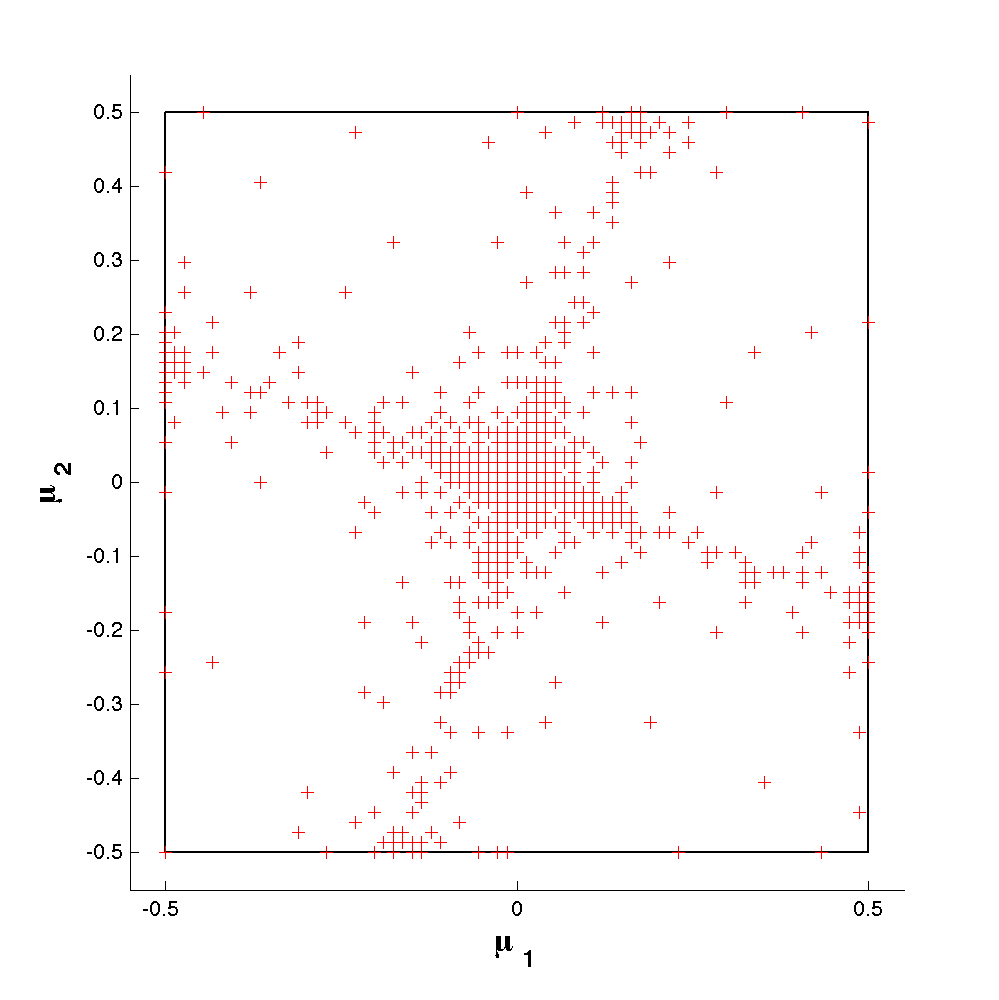}
\\
\includegraphics[Trim=\trimval,clip,width=\textwidth]{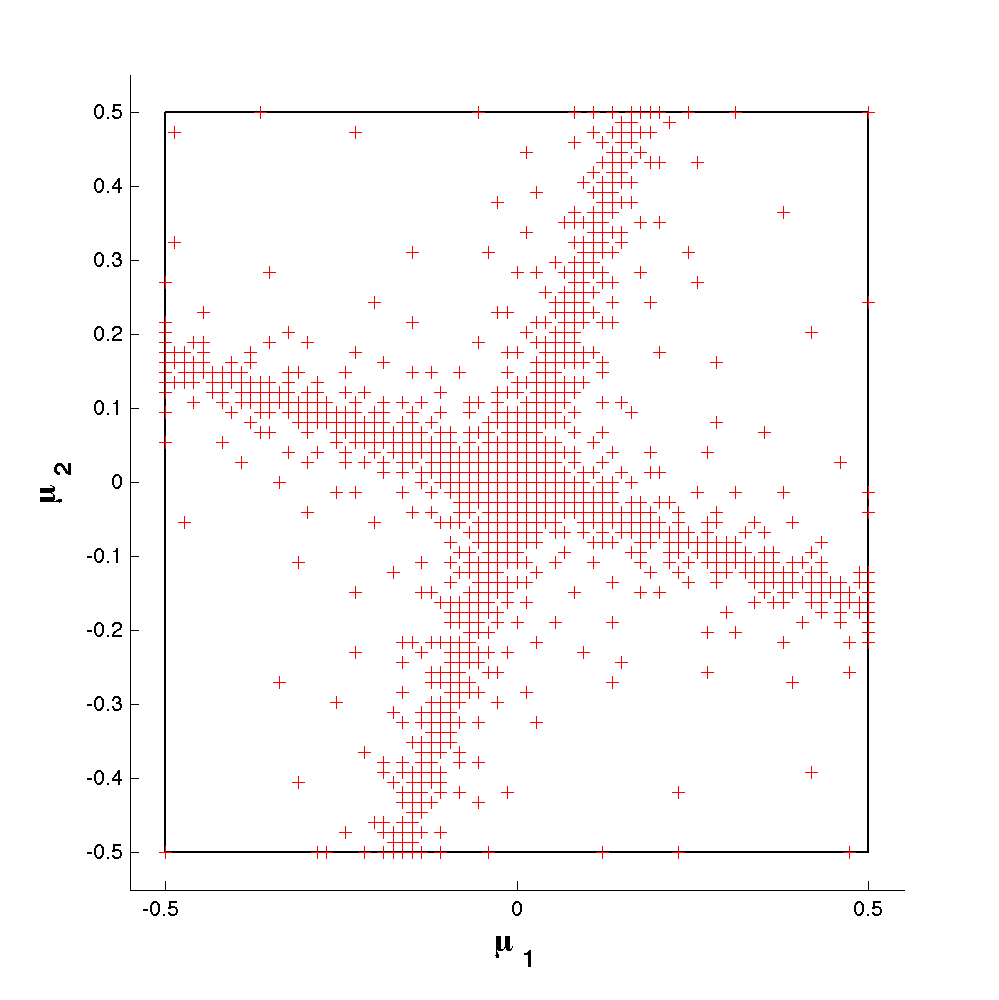}
\end{minipage}
}
\caption{
Test 3: Local approximation spaces for selected parameter values (left), radius as a function of the parameters (middle) and sample points (right) for the presented approach (top) in comparison wit the isotropic version (bottom) for $N=10$.
}
\label{fig:Test3}
\end{figure}

The local approximation spaces, the radius and the sample points with $N=10$ are presented in Figure \ref{fig:Test3}.
We observe {\em a posteriori} that the training space $\Xitrial$ was not sampled fine enough. Indeed in the region around the origin (and the cross for the isotropic version), every training point is included in the set of sample points. 
This also explains the sudden drop of the convergence in Figure \ref{fig:Ex3_conv}, in particular for the isotropic version. 
This is a known difficulty in greedy methods and our approach with adaptive training sets presents a solution to this problem. The corresponding numerical results are presented in the next section.

\begin{figure}[!ht]
  \captionsetup[subfigure]{labelformat=empty,width=\three}
  \centering
\subfloat[Tol=$10^{-3}$]{
\begin{minipage}[t]{\three}
\includegraphics[Trim=\trimval,clip,width=\textwidth]{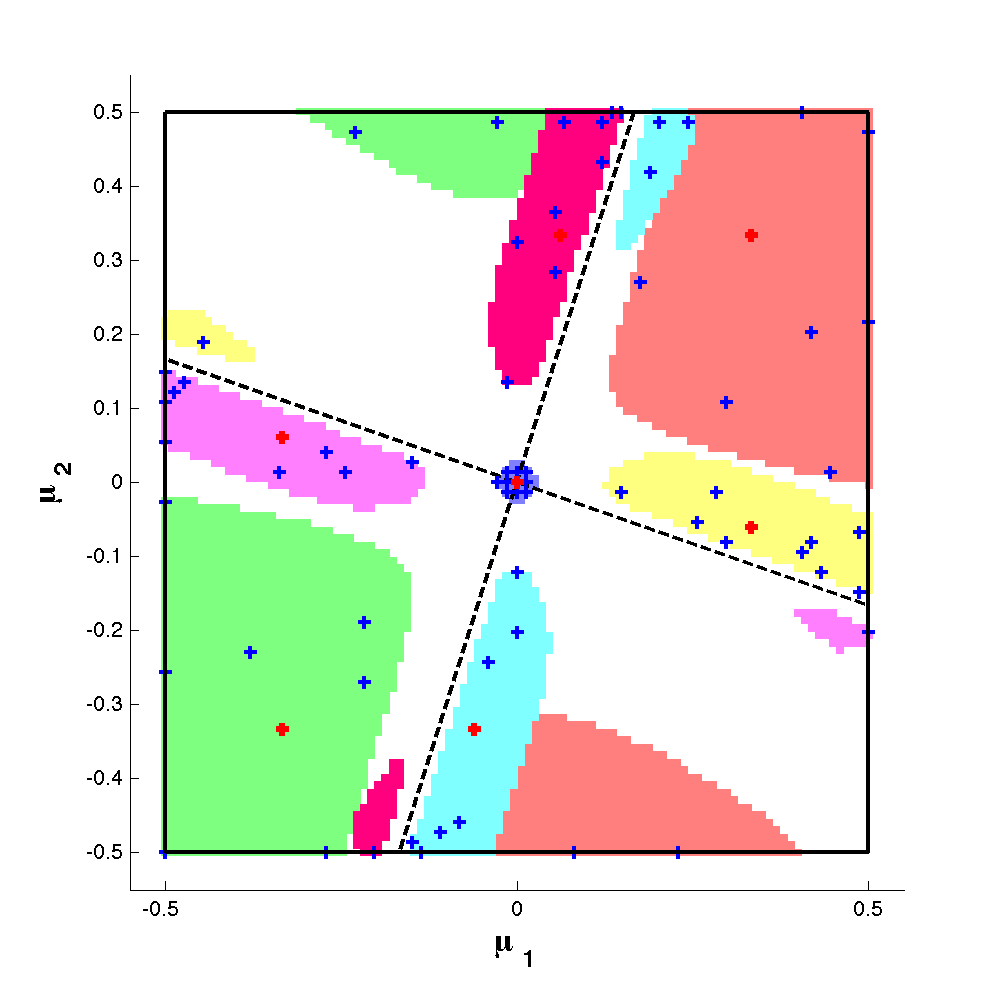}
\end{minipage}
}
\subfloat[Tol=$10^{-4}$]{
\begin{minipage}[t]{\three}
\includegraphics[Trim=\trimval,clip,width=\textwidth]{pics/LAS_TEST3_N10_TOL4.png}
\end{minipage}
}
\subfloat[Tol=$10^{-5}$]{
\begin{minipage}[t]{\three}
\includegraphics[Trim=\trimval,clip,width=\textwidth]{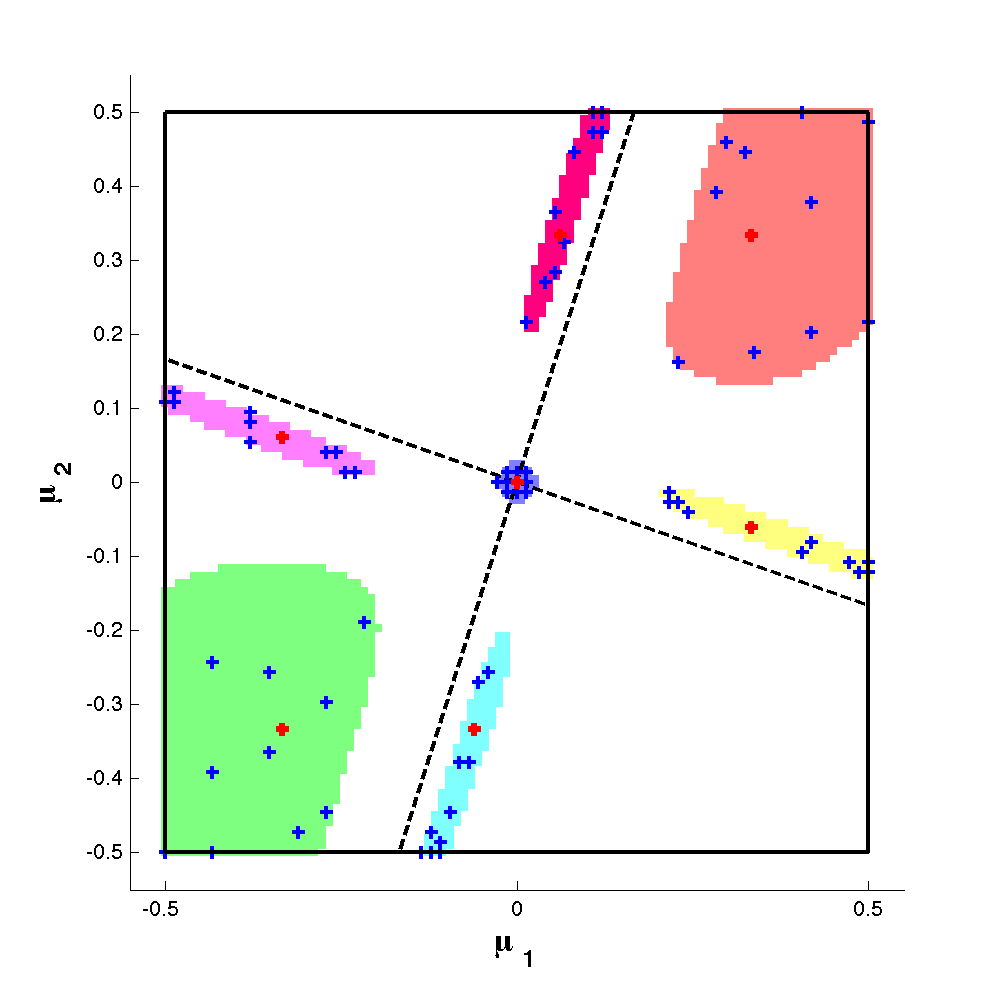}
\end{minipage}
}
\caption{
Test 3: Local approximation spaces for different tolerance levels.}
\label{fig:Ex3_LAS_tol}
\end{figure}

Further, we plot the local approximation spaces for different values of the tolerance levels in Figure \ref{fig:Ex3_LAS_tol}.
We observe that the originally non-connected approximation spaces are becoming connected with increased tolerance.

\begin{figure}[!ht]
  \captionsetup[subfigure]{labelformat=empty,width=\three}
  \centering
\subfloat[$\xi_1=\mu_1+3\mu_2$, $\xi_2=\frac12\mu_1-\mu_2$]{
\begin{minipage}[t]{\three}
\includegraphics[Trim=\trimval,clip,width=\textwidth]{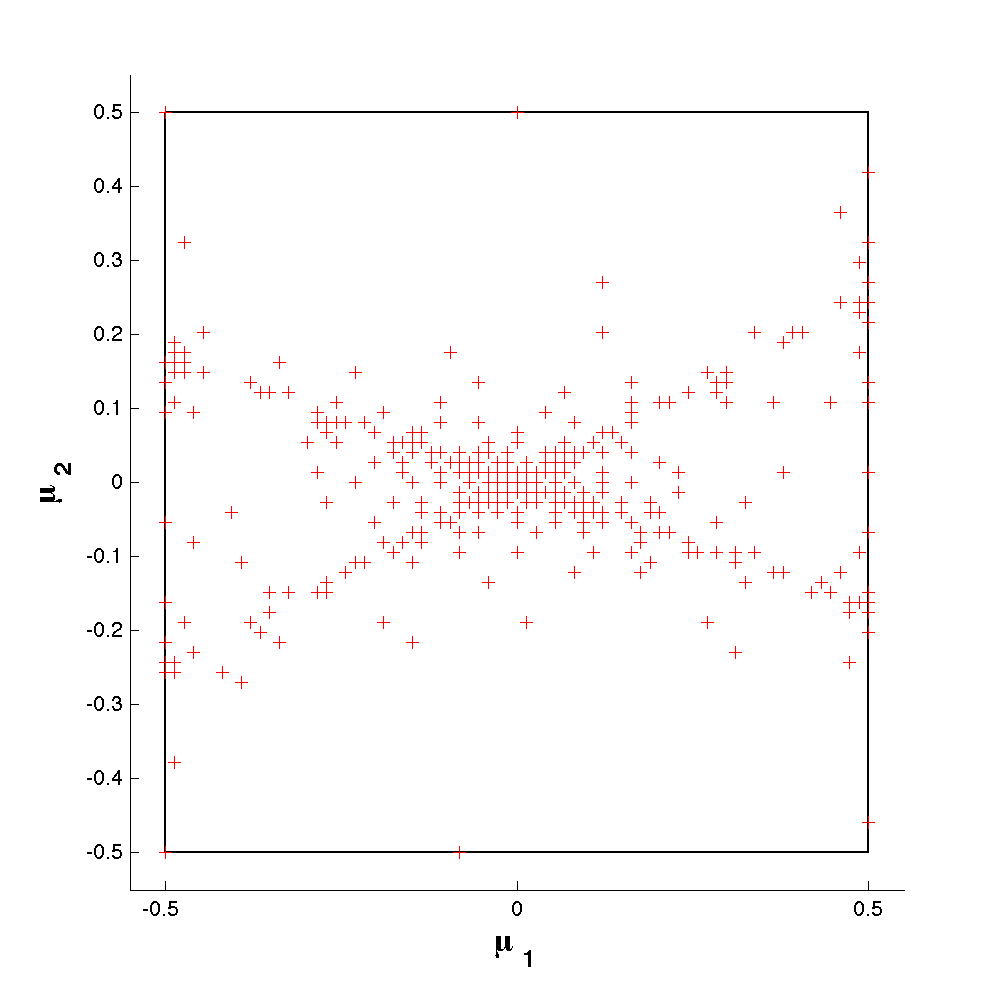}
\end{minipage}
}
\hspace{1cm}
\subfloat[$\xi_1=\mu_1+\mu_2$, $\xi_2=\mu_1-\mu_2$]{
\begin{minipage}[t]{\three}
\includegraphics[Trim=\trimval,clip,width=\textwidth]{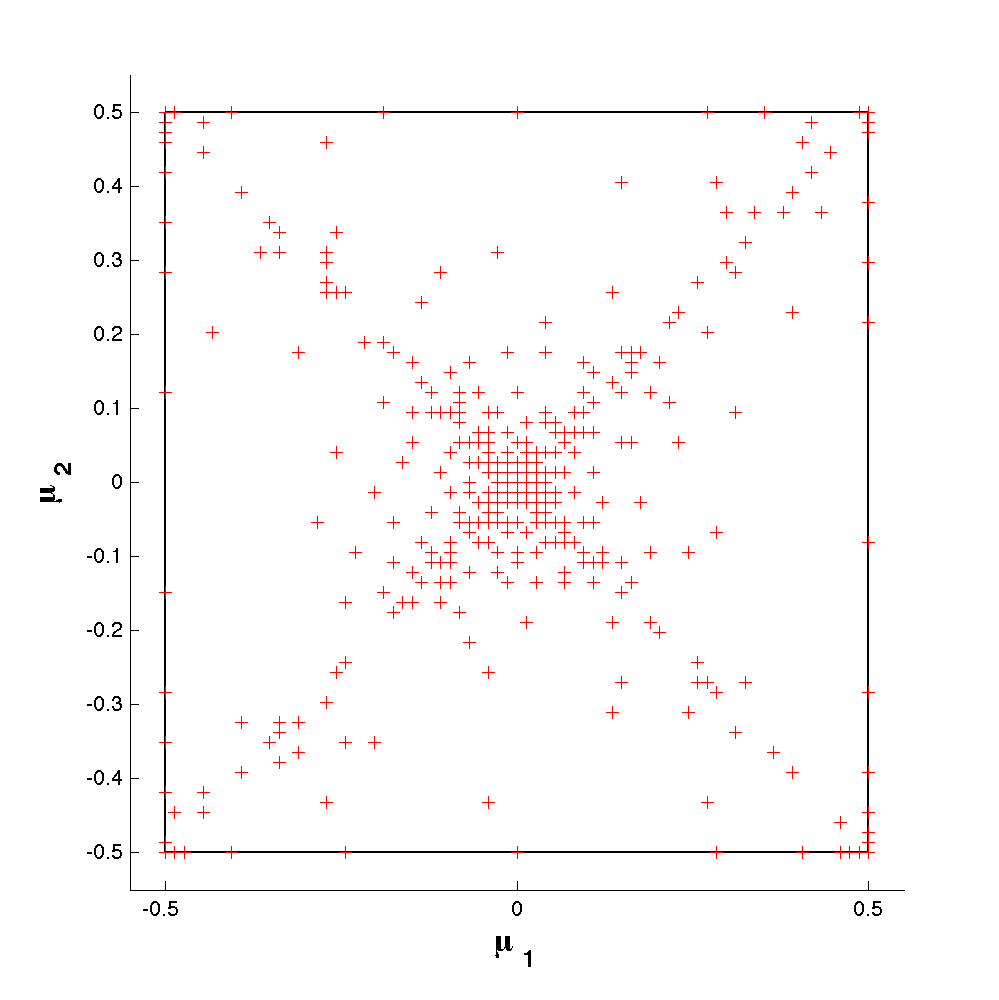}
\end{minipage}
}
\caption{
Sample set patterns for different type of parametrized functions with singular behavior on two crossing line with different slopes.
}
\label{fig:Test3_xi}
\end{figure}

We also recognize that the scheme detects the cross where the behavior of the function is most singular and present in Figure \ref{fig:Test3_xi} the sample points that were chosen for different functions of the type
\begin{eqnarray*}
	f_{3,\xi_1,\xi_2}(\bm x;\bmu) =
	\exp
	\left[
	-
	\frac{
	(x_1-\xi_1(\bmu))^2}{0.1+5|\xi_1|
	}
	-
	\frac{
	(x_2-\xi_2(\bmu))^2}{0.1+5|\xi_2|
	}
	\right]& \\
	\x\in \Omega=(-1,1)^2,&\,\bmu \in \dom=[-0.5,0.5]^2
\end{eqnarray*}
where $\xi_1$ and $\xi_2$ are functions of $\bmu=(\mu_1,\mu_2)$ using an error tolerance of $10^{-3}$ in this case.

\subsection{Numerical results with adapted training sets}
\label{ssec:ATS}
In order to compare the accuracy of this version  with adapted training sets we introduce a fixed test set of sample points $\Xitest$. In the following examples it consists of a lattice of $75\times 75$ points that corresponds to the fixed training set of the tests of the previous section. We use the command {\tt adaptmesh} in FreeFem \cite{Hecht:2008p27} in order to adapt the sample points by constructing a desired uniform mesh in the new metric defined by the Hessian.

The target advantage of our new approach with an adapting training set is not to use fewer computed solutions than the approach using a fixed training set, at least in low parameter's dimension but rather provide an improved reliability. Indeed
 a non-adapted training set may leave large errors in parameter's regions that are not sufficiently sampled by the  
 fixed training set. On the contrary the adapted set tracks these unexplored parameter's regions where the RB approximation error remains large. The points of the training set are chosen wisely using the information available from the local geometry of the system. 
%
Note however that for problem where the manifold $\mathbb W_{\mathbb P}$ happens to be uniformly regular, 
the benefit of our approach will be on the computational efficiency as a result of the use of a small training set in the beginning that is gradually increasing.

Figure \ref{fig:ATS_conv} illustrates the total number of basis functions that are necessary to be computed versus the achieved accuracy of the proposed algorithm with adaptive training sets for the test cases 1, 2 and 3. 
In each plot, we present the error of the version using an adaptive training set on the adaptive training set and the error evaluated on the fixed testing set $\Xitest$. The accuracy of the version with a fixed training set (which equals $\Xitest$) is given as comparison.

\subsubsection*{Test 1}
\begin{figure}[!ht]
  \captionsetup[subfigure]{labelformat=empty,width=\threep}
  \centering
\subfloat[Test 1: $N=20$]{
\includegraphics[Trim=\trimvalc,clip,width=\threep]{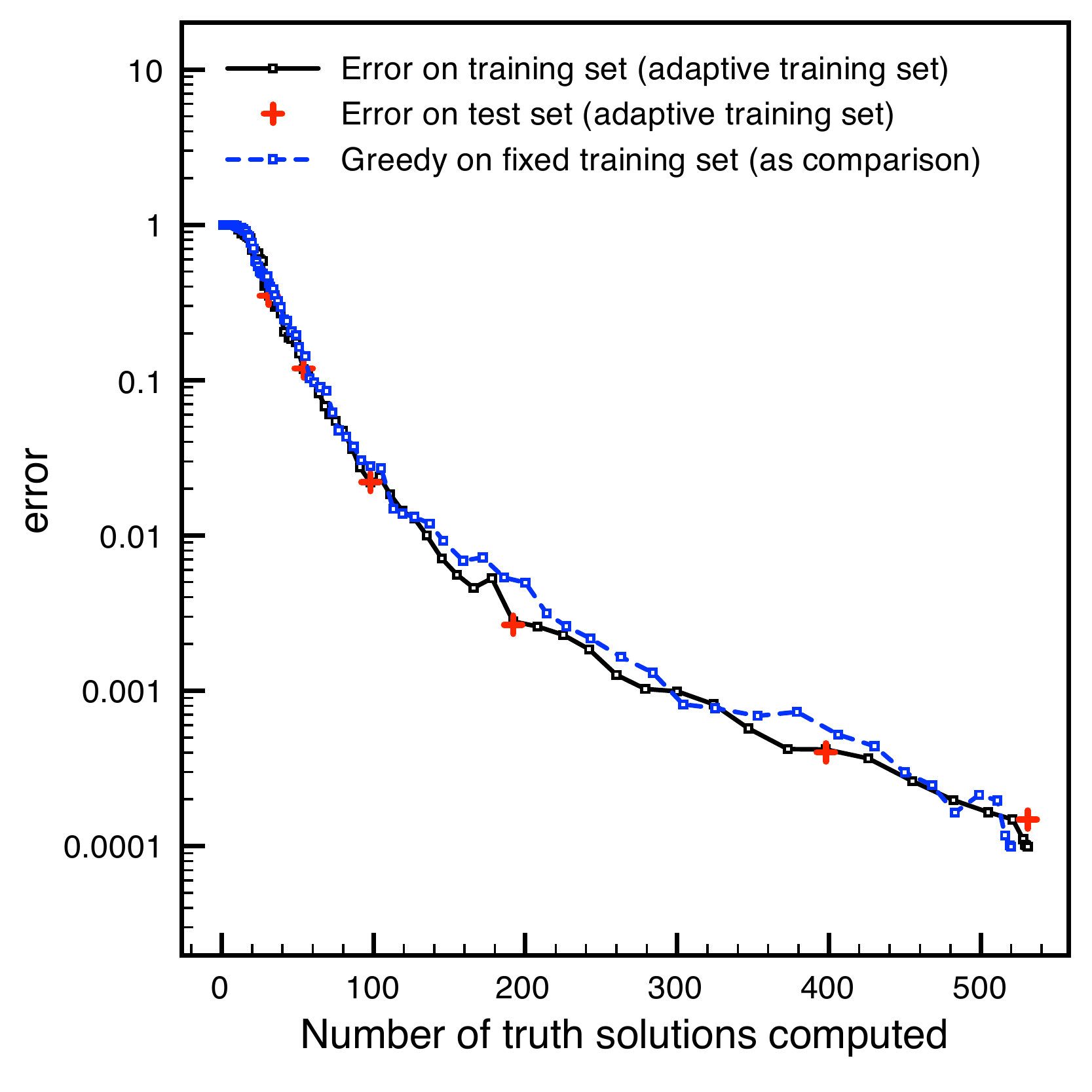}
}
\subfloat[Test 2: $N=5$]{
\includegraphics[Trim=\trimvalc,clip,width=\threep]{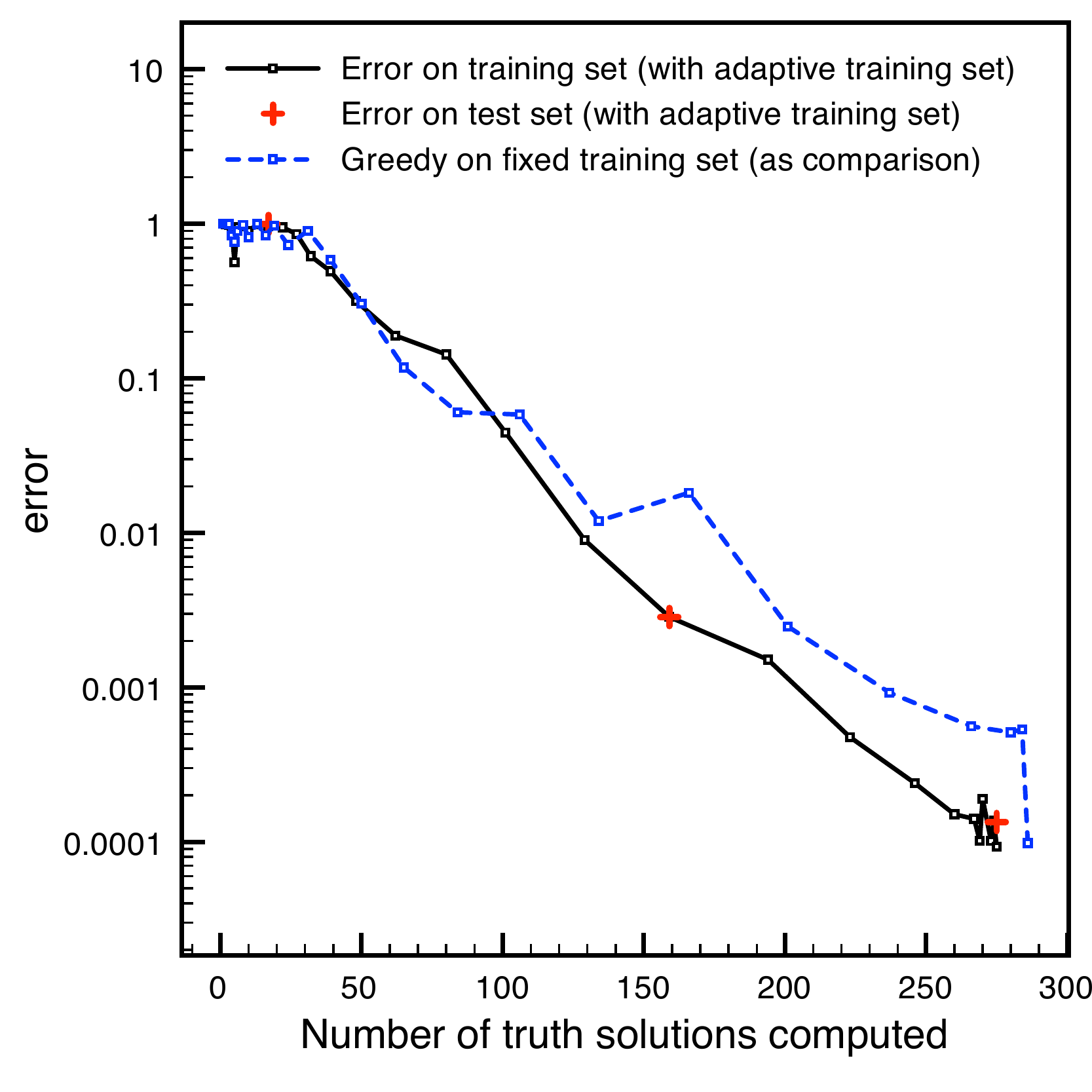}
}
\subfloat[Test 3: $N=10$]{
\includegraphics[Trim=\trimvalc,clip,width=\threep]{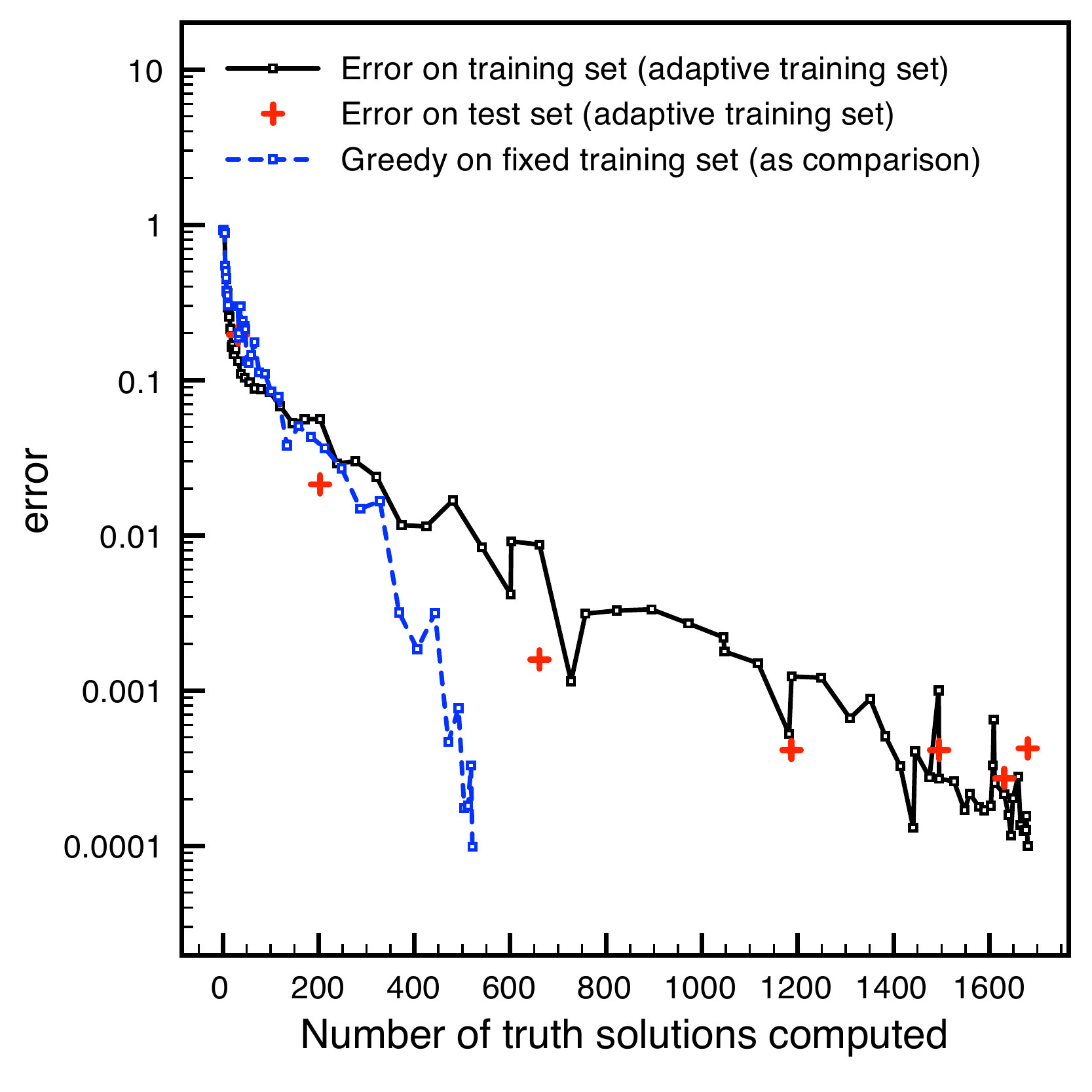}
}
\caption{Accuracy with respect to the number of truth solutions that need to be computed for all three test case.}
\label{fig:ATS_conv}
\end{figure}

From Figure \ref{fig:ATS_conv} (left), one can observe that the error on the training set and on the test set of the approach with an adaptive training set are similar.
This illustrates that the accuracy is not only guaranteed on the adaptive training set but also satisfied on the test grid.
Both errors are also similar to the error of the algorithm using a fixed training set as presented in Section \ref{ssec:FTS}.

The adaptive version is as good as the version on a fixed training set, but uses less evaluation of the error estimate as is shown in the following table:
\vspace{3pt}
\begin{center}
\begin{tabular}{lc}
\hline
\vspace{-8pt}\\
Number of training point evaluations using adaptive training sets: & 150'923\\
Number of training point evaluations using fixed training sets: & 455'625
\vspace{3pt}
\\
\hline
\end{tabular}
\end{center}
\vspace{5pt}
We observe a gain of 65\% fewer error evaluations. 

\begin{figure}[!ht]
  \captionsetup[subfigure]{labelformat=empty,width=\four}
  \centering
\subfloat[Training points.]{
\begin{minipage}[t]{\four}
\includegraphics[Trim=\trimval,clip,width=\textwidth]{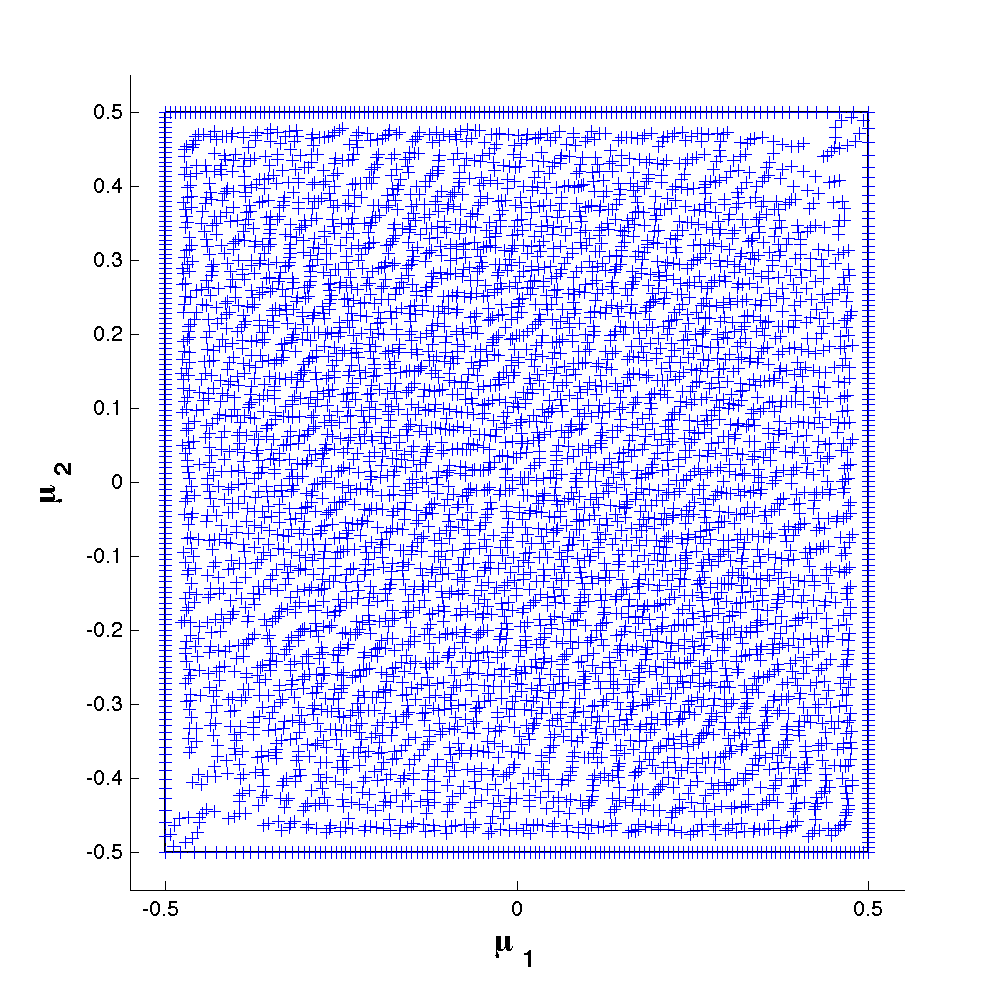}
\end{minipage}
}
\subfloat[Local approximation spaces.]{
\begin{minipage}[t]{\four}
\includegraphics[Trim=\trimval,clip,width=\textwidth]{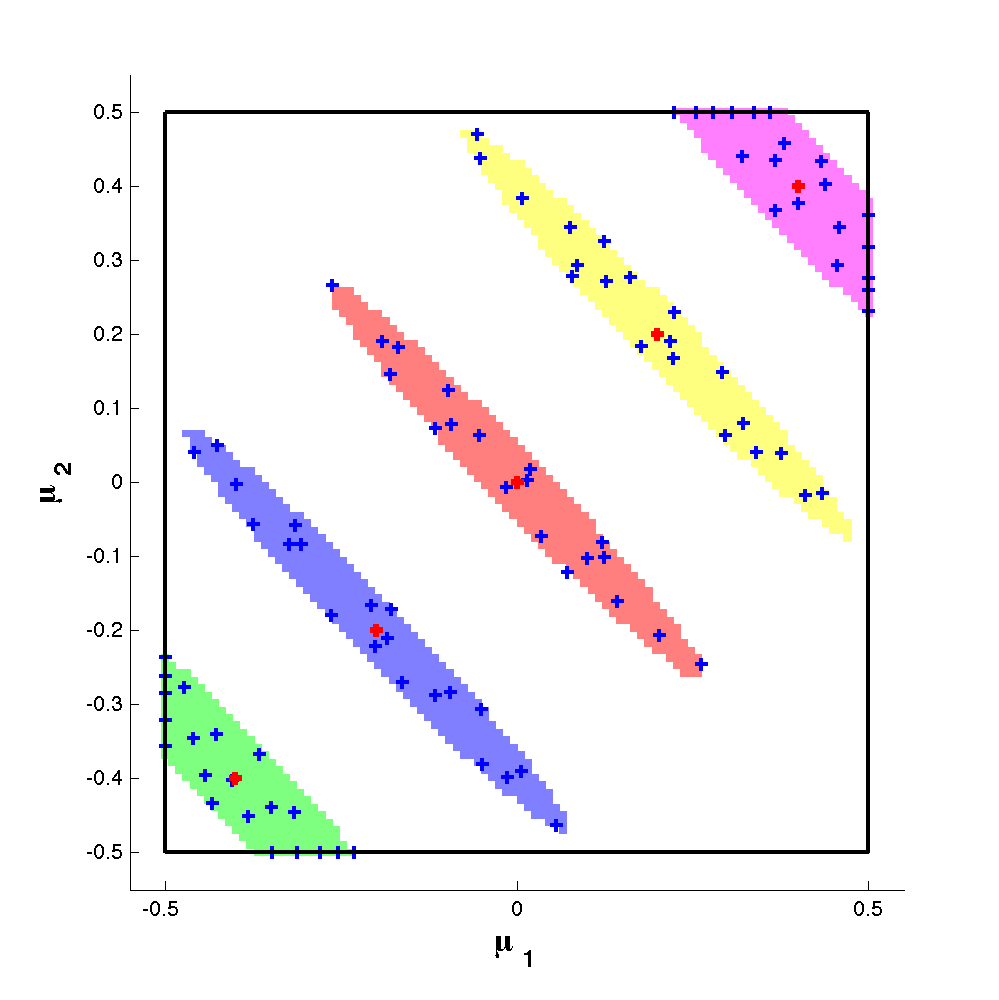}
\end{minipage}
}
\subfloat[Radius.]{
\begin{minipage}[t]{\four}
\includegraphics[Trim=\trimvalr,clip,width=\textwidth]{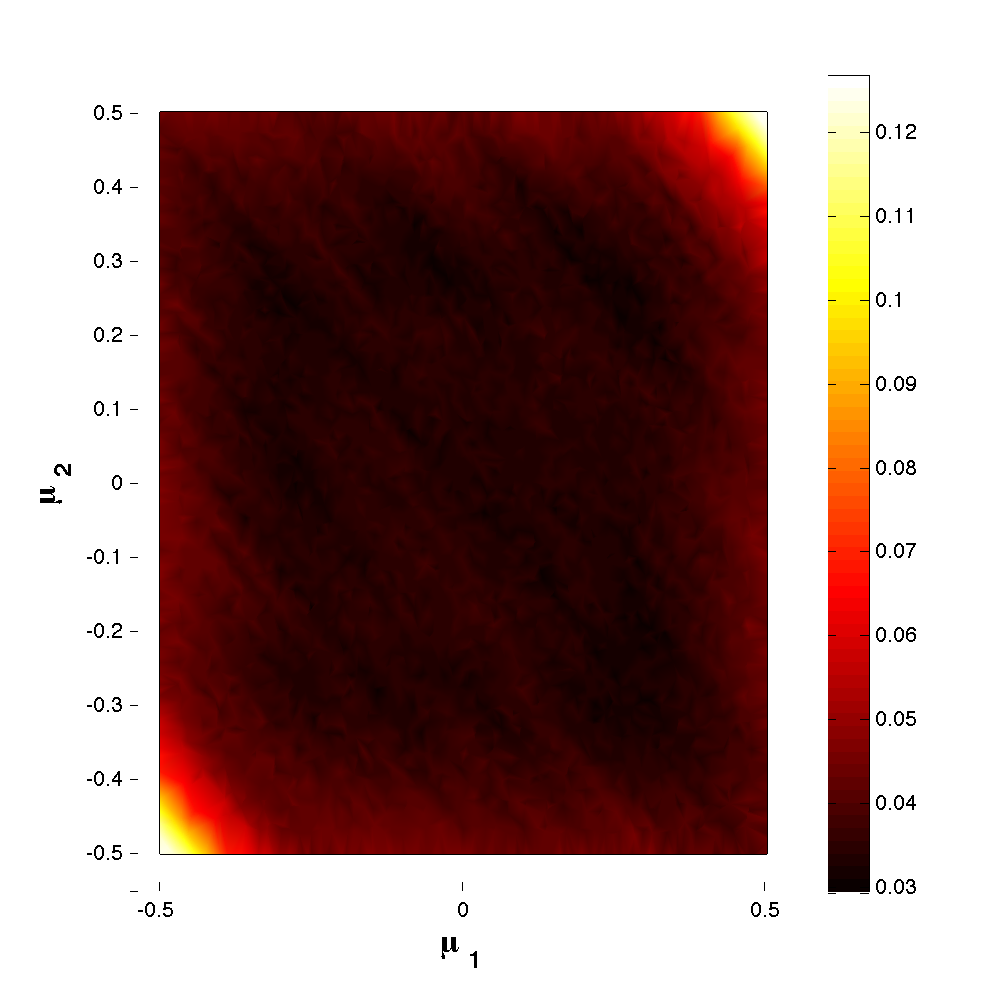}
\end{minipage}
}
\subfloat[Sample points.]{
\begin{minipage}[t]{\four}
\includegraphics[Trim=\trimval,clip,width=\textwidth]{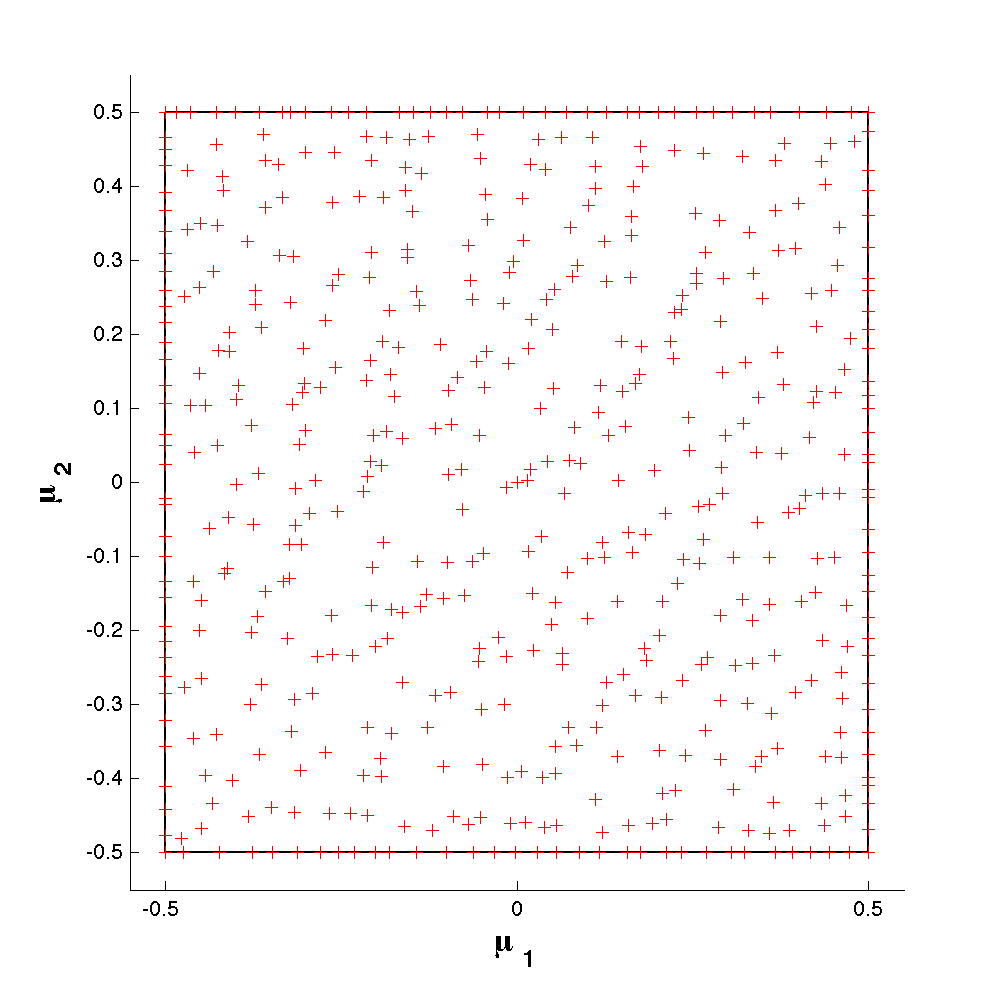}
\end{minipage}
}
\caption{
Test 1: Training points, local approximation spaces for some specific parameter values, radii and sample points for greedy algorithm with adaptive training set and $N=20$. 
}
\label{fig:Test1_ATS}
\end{figure}

Figure \ref{fig:Test1_ATS} presents the training set, local approximation spaces, radii and sample points for this calculation. 

\subsubsection*{Test 2}
Next, we consider the second example of the previous section. Figure \ref{fig:ATS_conv} (middle) compares this version (with adaptive training sets) with the version of the previous section using a fixed training set. One can observe again that the performance is similar than the version using a fixed training set, and that the accuracy is also satisfied on the fixed test set $\Xitest$.

Comparing again with the version using a fixed training set, the version with an adaptive training set requires less evaluation of the error estimate as is shown in the following table:
\vspace{3pt}
\begin{center}
\begin{tabular}{lc}
\hline
\vspace{-8pt}\\
Number of training point evaluations using adaptive training sets: & 60'202\\
Number of training point evaluations using fixed training sets: & 151'875
\vspace{3pt}
\\
\hline
\end{tabular}
\end{center}
\vspace{5pt}
We observe a gain of 70\% fewer error evaluations. 

\begin{figure}[!ht]
  \captionsetup[subfigure]{labelformat=empty,width=\four}
  \centering
\subfloat[Training points.]{
\begin{minipage}[t]{\four}
\includegraphics[Trim=\trimval,clip,width=\textwidth]{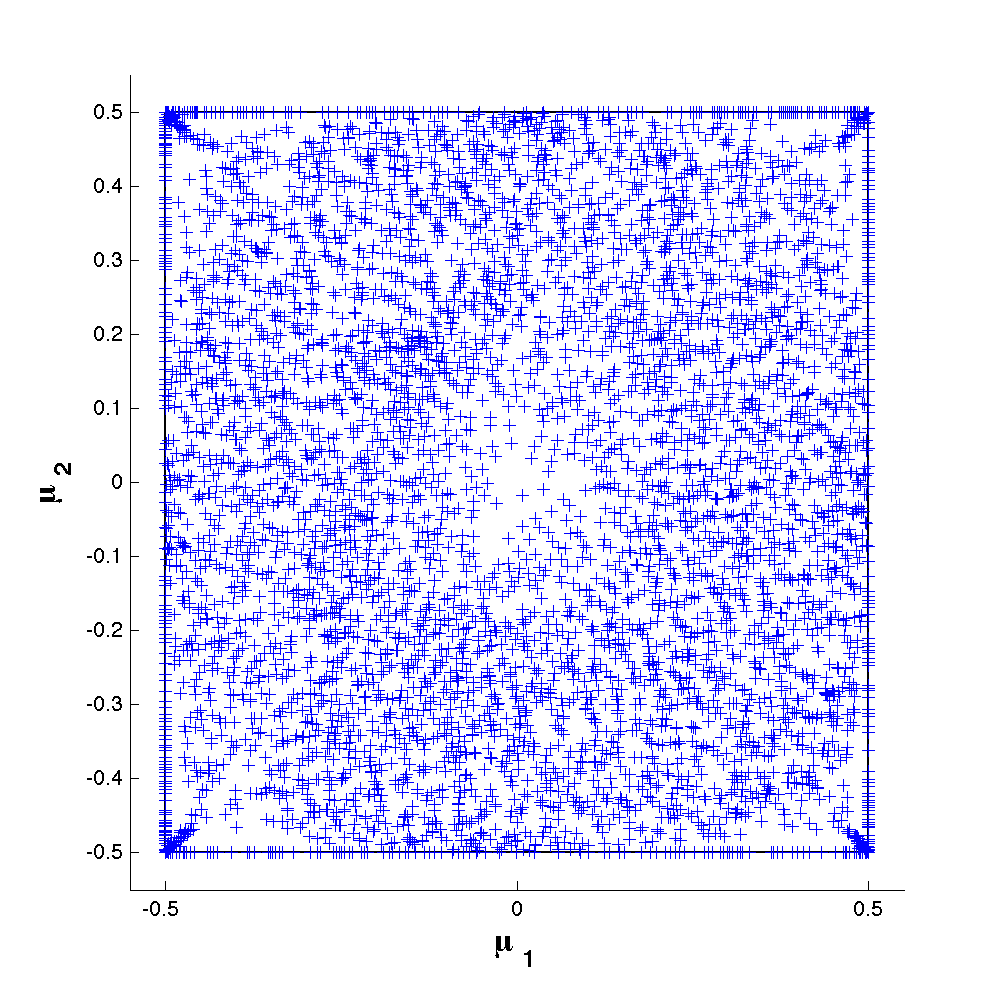}
\end{minipage}
}
\subfloat[Local approximation spaces.]{
\begin{minipage}[t]{\four}
\includegraphics[Trim=\trimval,clip,width=\textwidth]{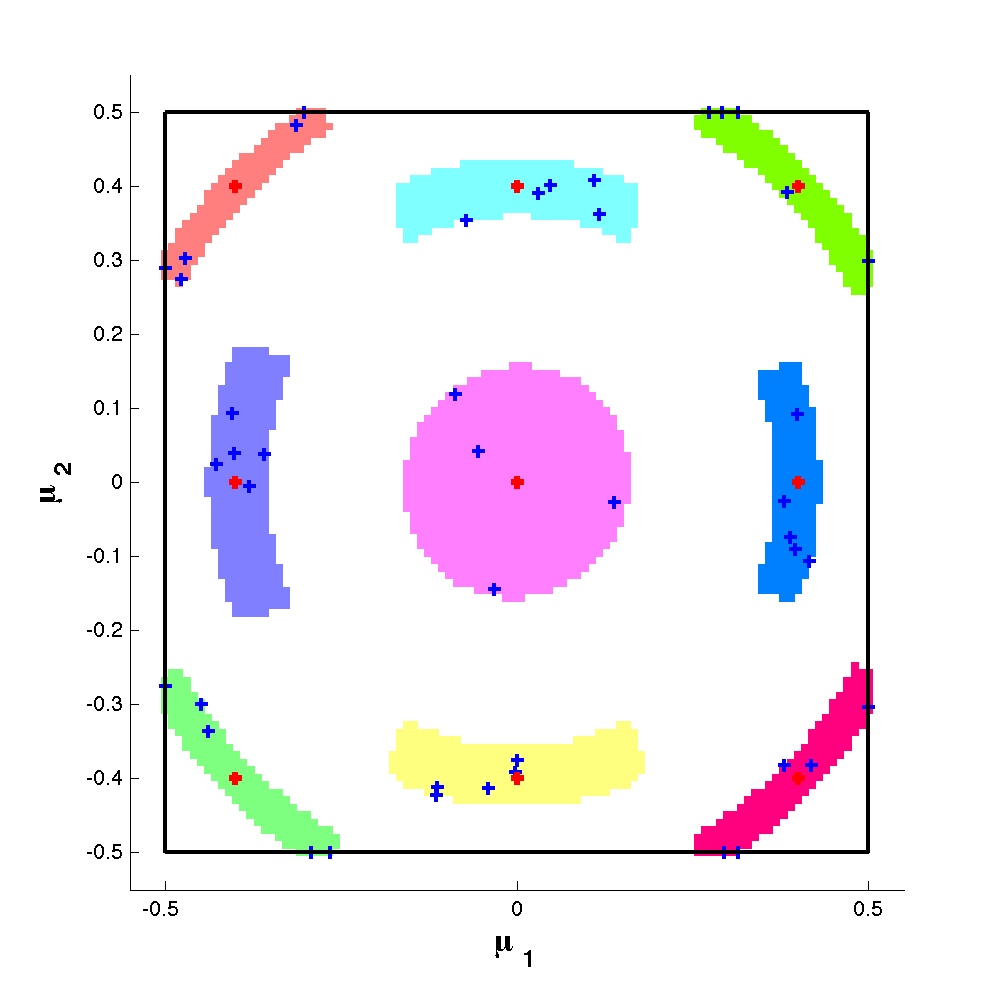}
\end{minipage}
}
\subfloat[Radius.]{
\begin{minipage}[t]{\four}
\includegraphics[Trim=\trimvalr,clip,width=\textwidth]{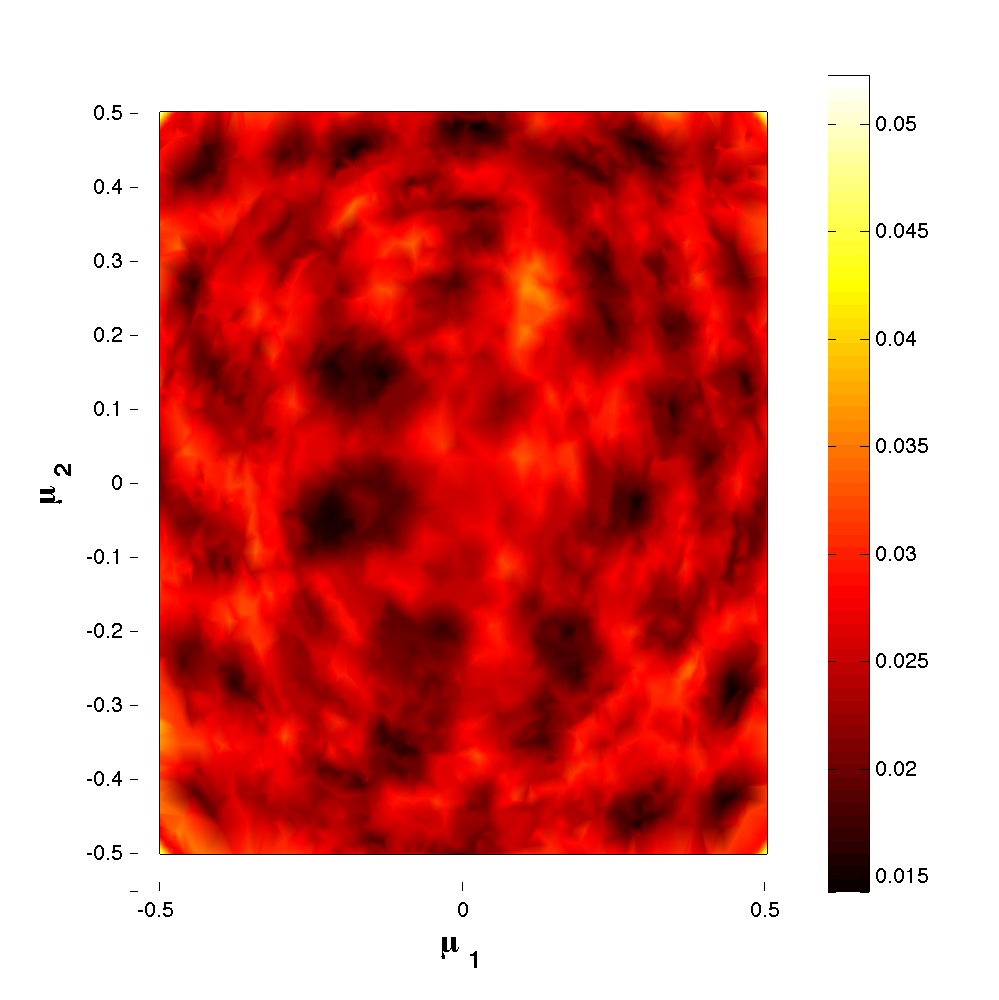}
\end{minipage}
}
\subfloat[Sample points.]{
\begin{minipage}[t]{\four}
\includegraphics[Trim=\trimval,clip,width=\textwidth]{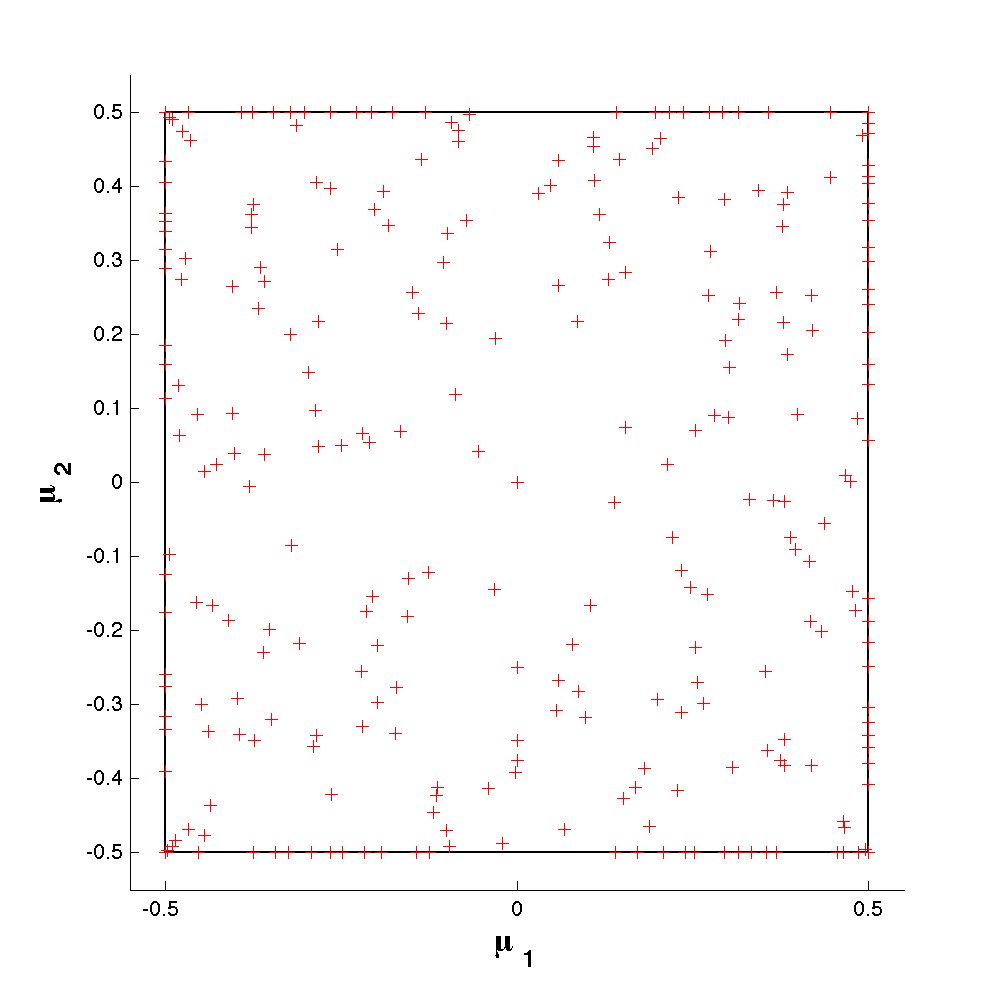}
\end{minipage}
}
\caption{
Test 2: Training points, local approximation spaces for some specific parameter values, radii and sample points for greedy algorithm with adaptive training set and $N=5$. 
}
\label{fig:Test2_ATS}
\end{figure}

Figure \ref{fig:Test2_ATS} presents the training set, local approximation spaces, radii and sample points of this calculation.

\subsubsection*{Test 3}
\begin{figure}[!ht]
  \captionsetup[subfigure]{labelformat=empty,width=\four}
  \centering
\subfloat[Training points.]{
\begin{minipage}[t]{\four}
\includegraphics[Trim=\trimval,clip,width=\textwidth]{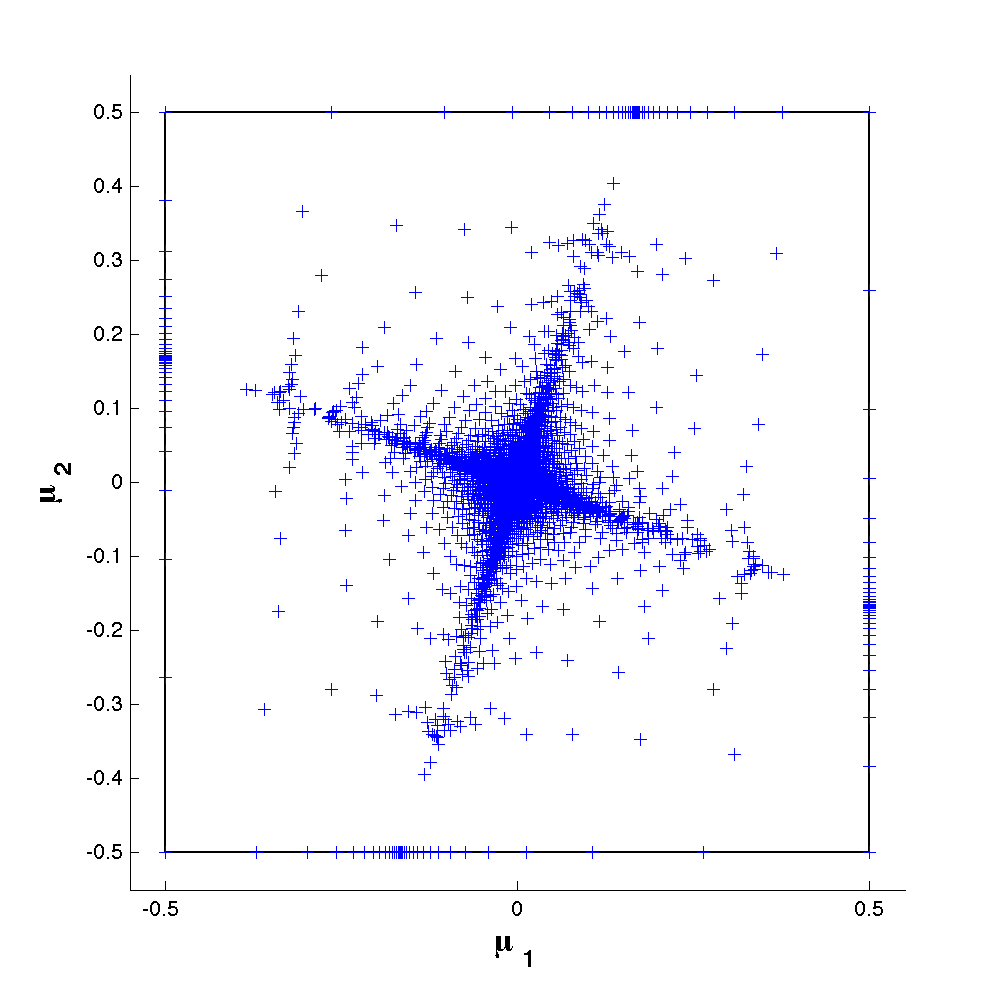}
\end{minipage}
}
\subfloat[Local approximation spaces.]{
\begin{minipage}[t]{\four}
\includegraphics[Trim=\trimval,clip,width=\textwidth]{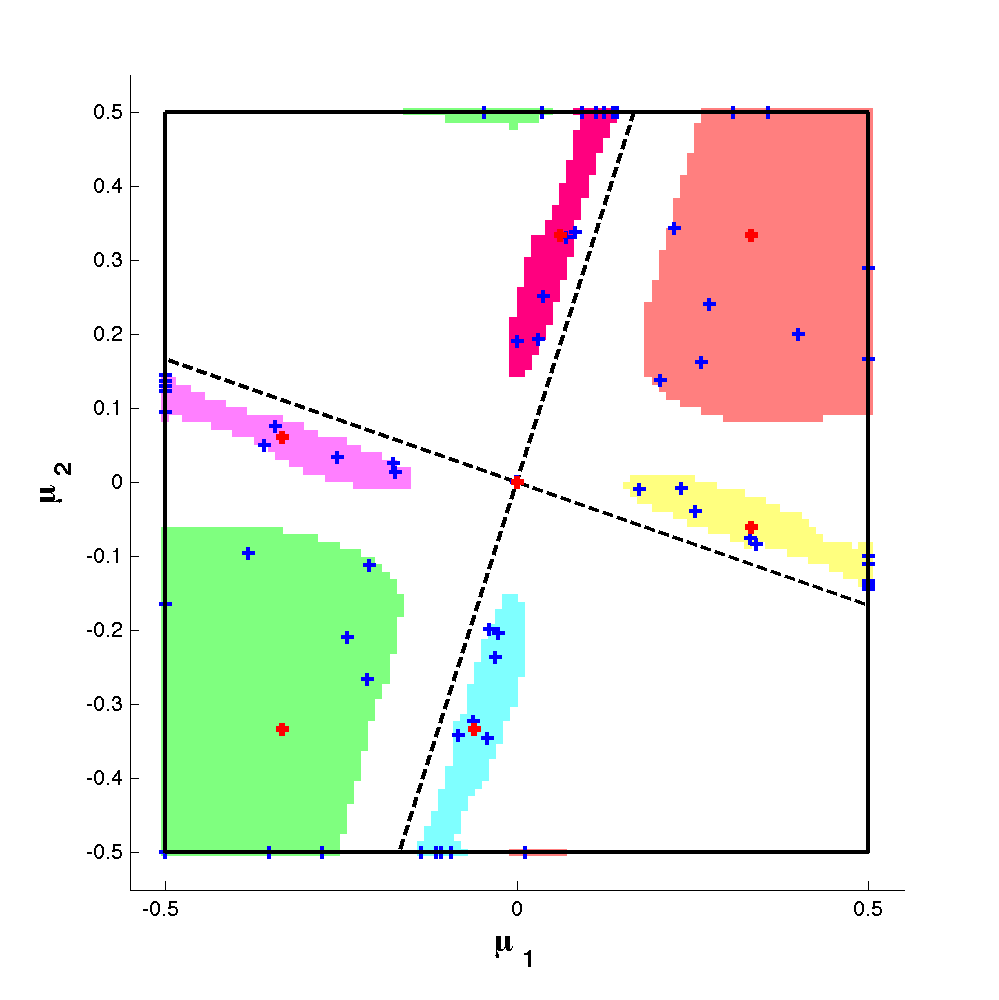}
\end{minipage}
}
\subfloat[Radius.]{
\begin{minipage}[t]{\four}
\includegraphics[Trim=\trimvalr,clip,width=\textwidth]{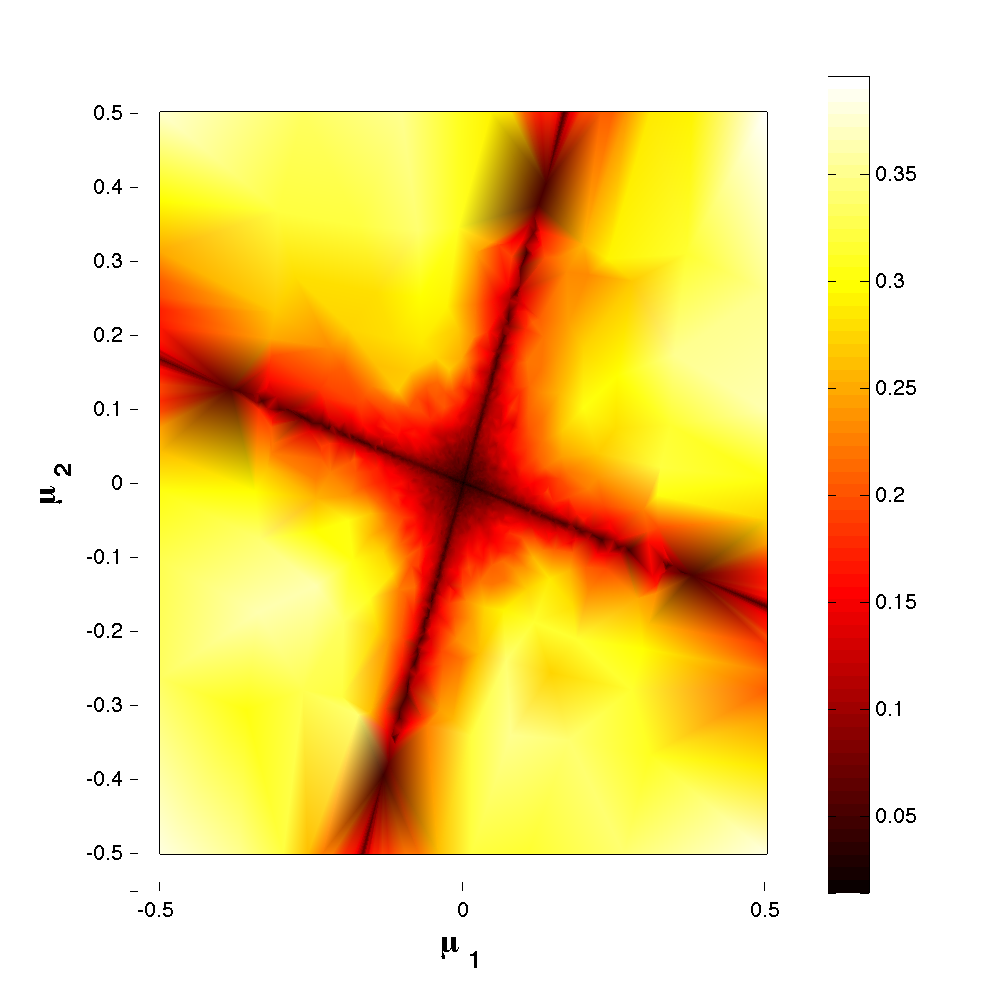}
\end{minipage}
}
\subfloat[Sample points.]{
\begin{minipage}[t]{\four}
\includegraphics[Trim=\trimval,clip,width=\textwidth]{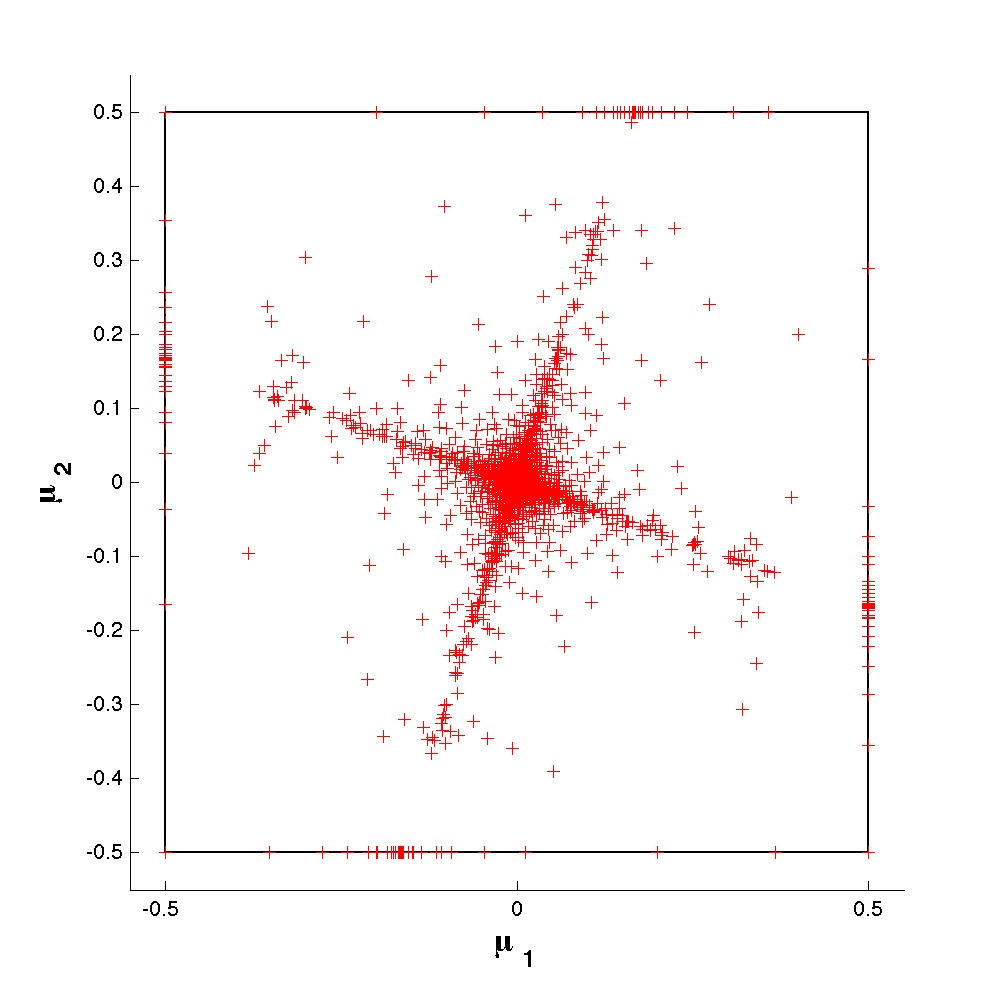}
\end{minipage}
}
\caption{
Test 3: Training points, local approximation spaces for some specific parameter values, radii and sample points for greedy algorithm with adaptive training set and $N=10$. 
}
\label{fig:Test3_ATS}
\end{figure}

Finally, we apply the new algorithm with an adaptive training set to the third example. As already mentioned, the solution of this example shows very strong changes with respect to the parameter and the version with a fixed training set does not resolve this since the training set was not sampled fine enough.  
We expect the version with the adaptive training set to perform accurately on the entire parameter space, and thus in contrast using more sample points than the version with a fixed training set. 
Figure \ref{fig:ATS_conv} (right) presents the evolution of the number of sample points versus the accuracy of the algorithm. One can observe that this approach (with adaptive training set) uses indeed more sample points than the version on a fixed training set as explained above. 
Finally, Figure \ref{fig:Test3_ATS} presents the training set, local approximation spaces, radii and sample points.
We observe that the training points are in accordance with the sample points, and that a more dense sampling is indeed required around the origin.  We observe that the local approximation spaces are now connected also for the tolerance of $10^{-4}$.

\begin{figure}[!ht]
  \captionsetup[subfigure]{labelformat=empty,width=\two}
  \centering
\subfloat[Fixed training set.]{
\begin{minipage}[t]{\two}
\includegraphics[Trim=\trimvalc,clip,width=\textwidth]{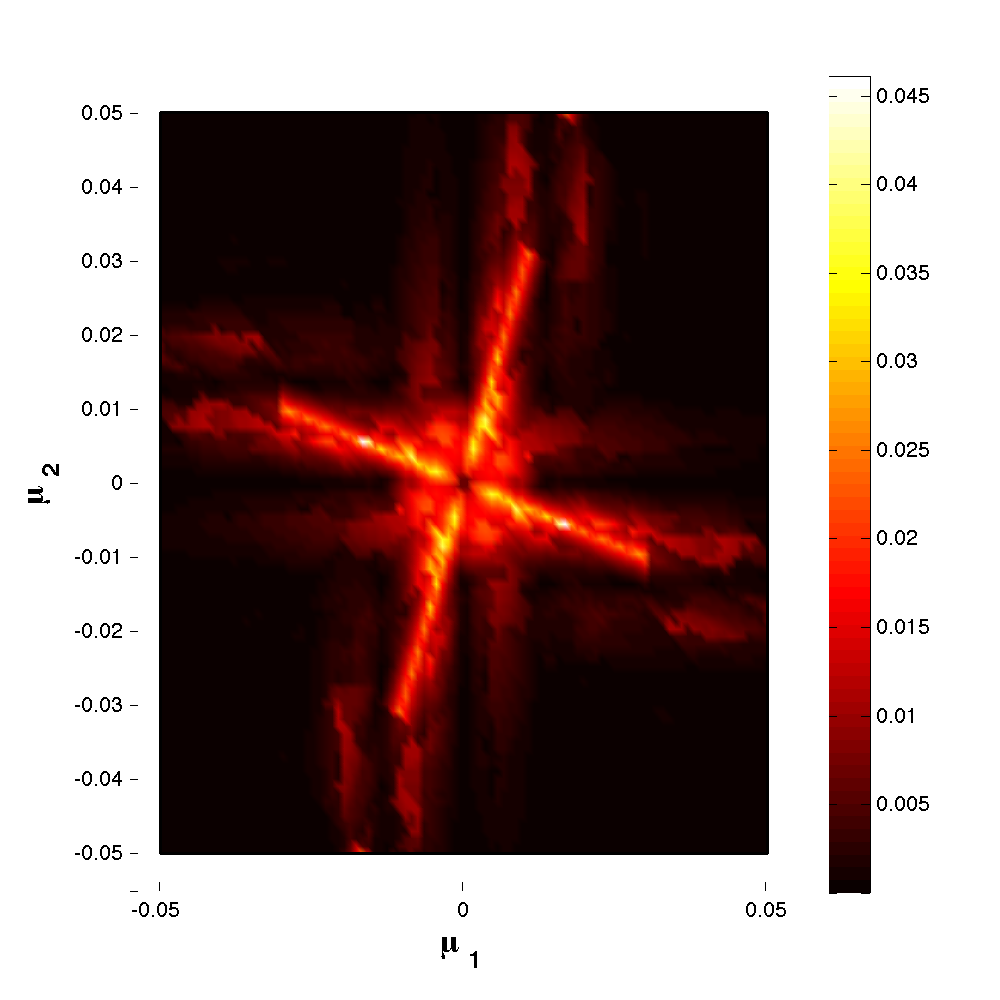}
\end{minipage}
}
\subfloat[Adaptive training set.]{
\begin{minipage}[t]{\two}
\includegraphics[Trim=\trimvalc,clip,width=\textwidth]{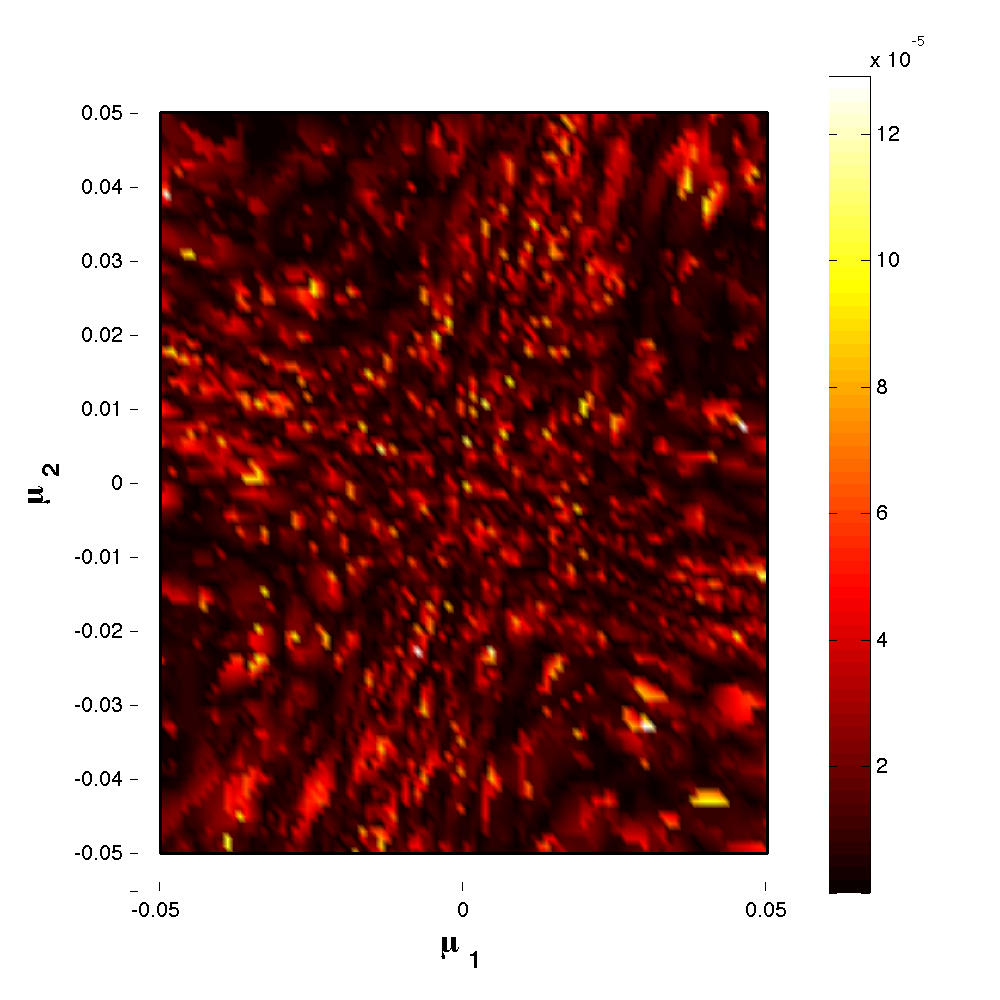}
\end{minipage}
}
\caption{
Test 3: Error distribution on a region around the origin of the greedy version using a fixed (left) and an adaptive (right) training set. 
}
\label{fig:Test3_ATS_patch}
\end{figure}

In order to illustrate the benefit of using adaptive training sets also in this case we consider a test sample of $100\times 100$ uniformly distributed points in the region $[-0.05,0.05]^2$ around the origin. Figure \ref{fig:Test3_ATS_patch} illustrates the error distribution using the online procedure generated using a fixed (left) and an adaptive (right) training set. The maximum error is 0.046 resp. $1.29\cdot 10^{-4}$. While the error tolerance is almost satisfied in the latter case, it is clearly not the case for the former approach. 
Further, the number of error evaluations is given in the following table:
\vspace{3pt}
\begin{center}
\begin{tabular}{lc}
\hline
\vspace{-8pt}\\
Number of training point evaluations using adaptive training sets: & 181'672\\
Number of training point evaluations using fixed training sets: & 247'500
\vspace{3pt}
\\
\hline
\end{tabular}
\end{center}
\vspace{5pt}
Thus, the adaptive version uses still less error estimator evaluations and is more accurate in the region around the origin. Also, the error is more equally distributed.

\section{Numerical results using a Galerkin-projection}
\label{sec:Galerkin}

\subsubsection*{Test 4}
\begin{figure}[!ht]
  \captionsetup[subfigure]{labelformat=empty,width=\three}
  \centering
\subfloat[$\bmu=(-2.5,-\pi/4)$.]{
\begin{minipage}[t]{\three}
\includegraphics[Trim=\trimval,clip,width=0.9\textwidth]{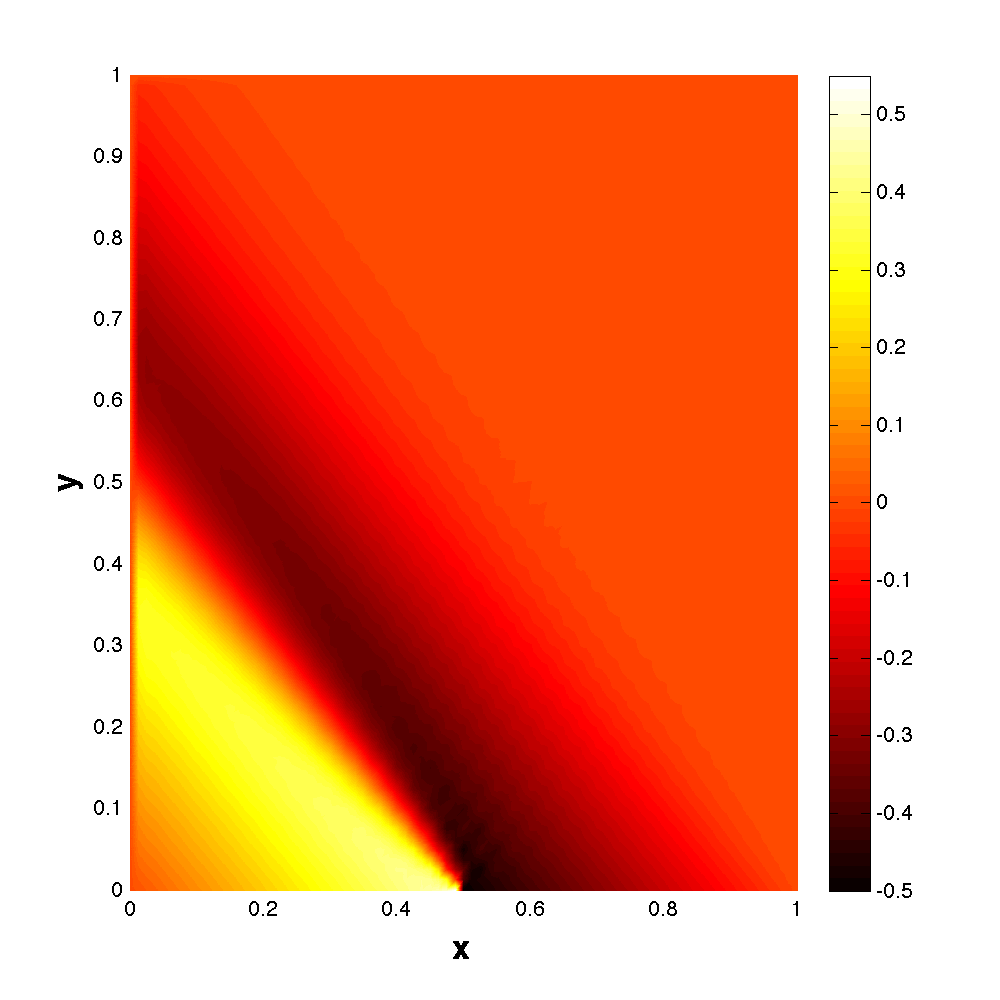}
\end{minipage}
}
\subfloat[$\bmu=(-2,0)$.]{
\begin{minipage}[t]{\three}
\includegraphics[Trim=\trimval,clip,width=0.9\textwidth]{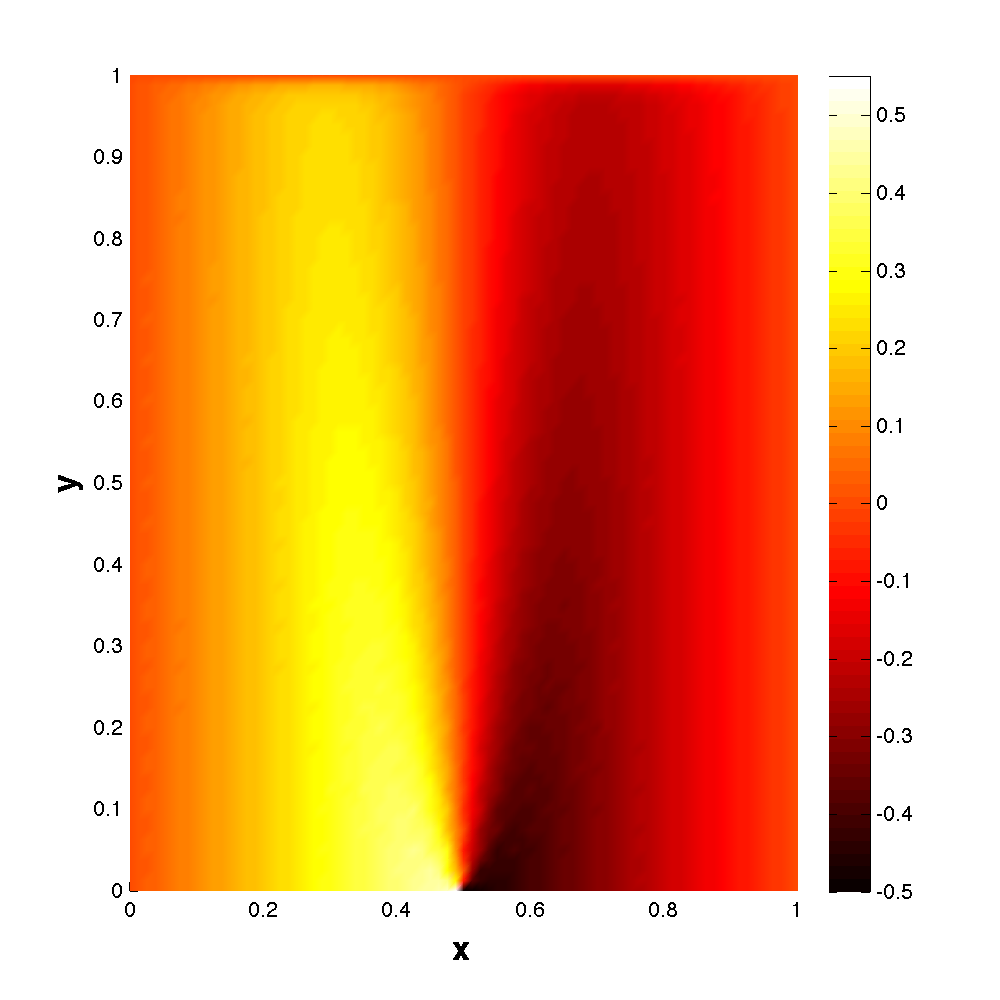}
\end{minipage}
}
\subfloat[$\bmu=(-1,\pi/4)$.]{
\begin{minipage}[t]{\three}
\includegraphics[Trim=\trimval,clip,width=0.9\textwidth]{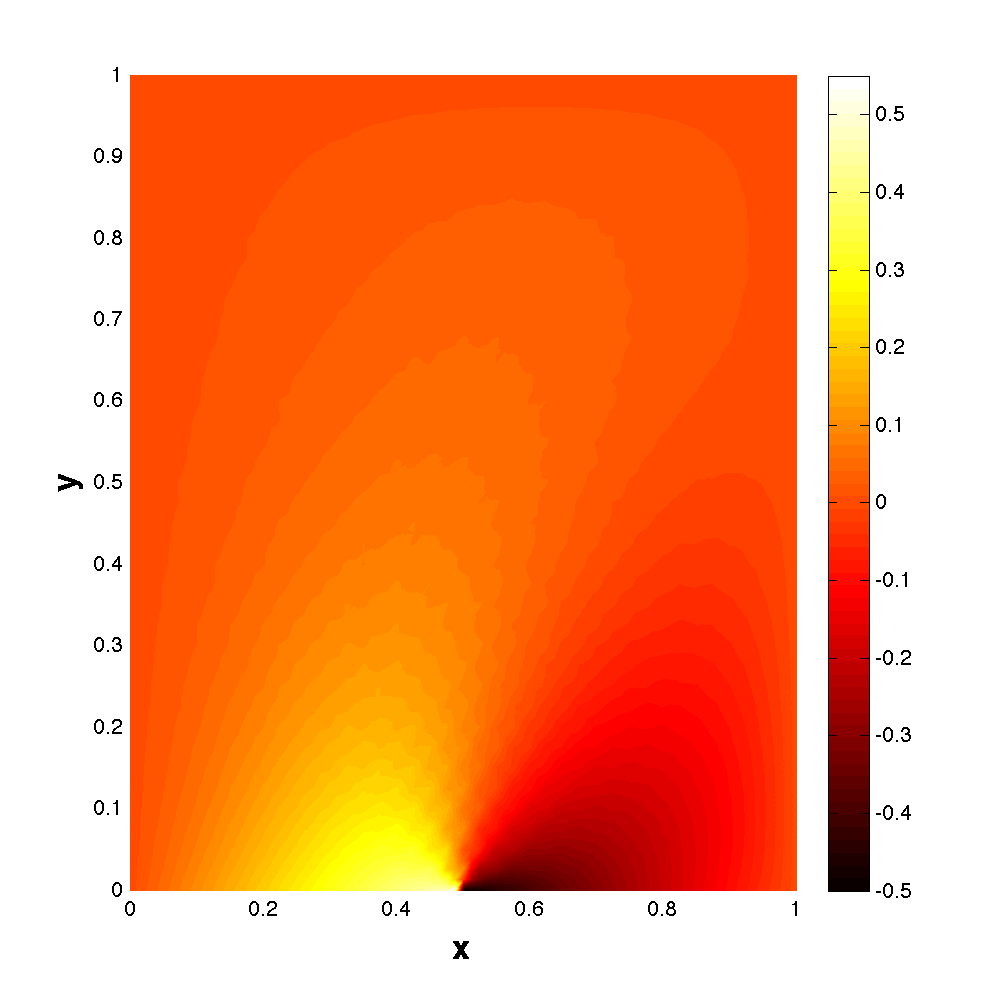}
\end{minipage}
}
\caption{
Test 4: Truth solutions for different parameter values which consist of the diffusion coefficient and the direction of the convection field. 
}
\label{fig:RB_truth}
\end{figure}

Finally, we present an example involving the full reduced basis framework. 
We consider a two-dimensional steady convection-diffusion equation on $\Omega=(0,1)^2$. Homogeneous Dirichlet boundary conditions are imposed on the left, upper and right boundary, i.e. on $\Gamma_0 =\{(x,y)\in\partial \Omega\,|\, x=0, \mbox{ or } x=1, \mbox{ or } y=1\}$.
On the lower boundary, i.e. on $\Gamma_g = \{(x,0)\in\partial \Omega\,|\, x\in[0,1]\}$, we impose the following Dirichlet data
\[
	u(x,0) = g(x) := x \, \chi_{[0,0.5)} + (-40\,x+21) \, \chi_{[0.5,0.55]} + (x-1) \, \chi_{[0.525,1]},
\]
(essentially a slightly regularized saw tooth function) where $\chi_X$ denotes the characteristic function for the mono-dimensional region $X$.
Thus, the boundary condition presents a strong gradient between $(0.5,0)$ and $(0.525,0)$.
Next, we introduce two parameters, the diffusion coefficient $\varepsilon=10^{\mu_1}$ and the direction of the convection parametrized by an angle $\mu_2$:
\[
	\bm\beta(\mu_2) = (\sin(\mu_2),\cos(\mu_2))^{\rm T}.
\]
Denoting the parameter vector by $\bmu=(\mu_1,\mu_2)\in \dom = [-2.5,0]\times[-\pi/4,\pi/4]$, the parametrized problem writes as: for any $\bmu\in\dom$, find $u(\bmu)\in H^1(\Omega)$ such that
\begin{align}
	\label{eq:StrongConvDiff}
	-10^{\mu_1}\Delta u(\bmu) + \nabla\cdot(\bm\beta(\mu_2) u(\bmu)) &=0 ,\qquad\mbox{in } \Omega,\\
	u(\bmu) &= g, \qquad\mbox{on } \Gamma_g,\\
	u(\bmu) &= 0, \qquad\mbox{on } \Gamma_0.
\end{align}
The solutions presents a convected near discontinuity, that is also subject to a diffusion process, along a parametrized velocity field and parametrized diffusion coefficient. 
To solve the truth problem, we use a finite element solver using a uniform quadrilateral mesh and tensorized polynomials of degree 3 with mesh size $h=0.025$.
Different solutions corresponding to different parameter values are illustrated in Figure \eqref{fig:RB_truth} and we highlight that the discretization is fine enough in order to avoid strong oscillations due to under-resolution related to the dominant convection for all parameter values.
For the error computation/estimation that is used in the greedy algorithm, we use for sake of simplicity the exact error between the truth solution and the reduced basis approximation.

This problem can be cast as a variational problem with bilinear- and linear forms that are affine decomposable as 
described in \eqref{eq:affdec} with
\begin{align*}
	a_1(w,v) &= \int_{\Omega} \nabla w\,\nabla v\,d\bm x, & g_1(\bmu) &= 10^{\mu_1},\\
	a_2(w,v) &= \int_{\Omega} \partial_x w\, v\,d\bm x, & g_2(\bmu) &= \sin(\mu_2),\\
	a_3(w,v) &= \int_{\Omega} \partial_y w\,v\,d\bm x, & g_3(\bmu) &= \cos(\mu_2).
\end{align*}
The force field is equal to zero in this case, however in the practical implementation we impose the boundary condition nodally and project the equation onto the nullspace of the (discrete) Dirichlet trace operator. In this manner, the boundary condition appears as a right hand side. 
Indeed, the weak form of \eqref{eq:StrongConvDiff} can be stated as:
for any $\bmu\in\dom$, find $\mathring{u}(\bmu)\in H^1_0(\Omega)$ such that
\[
	a(\mathring{u}(\bmu),v;\bmu) = -a(u_g,v;\bmu),\qquad \forall v\in H^1_0(\Omega),
\]
where $u(\bmu) = \mathring{u}(\bmu)+u_g$ and $u_g|_{\Gamma_g}=g$.
Thus we have a modified linear form $\tilde f(v;\bmu)$ with an affine decomposition inherited of the one from the bilinear form, i.e.
\[
	\tilde f(v;\bmu) = -a(u_g,v;\bmu) = -\sum_{q=1}^{3} g_q(\bmu) \, a_q(u_g,v).
\]

\begin{figure}[!ht]
  \captionsetup[subfigure]{labelformat=empty,width=\three}
  \centering
\subfloat[Local approximation spaces.]{
\begin{minipage}[t]{\three}
\includegraphics[Trim=\trimval,clip,width=\textwidth]{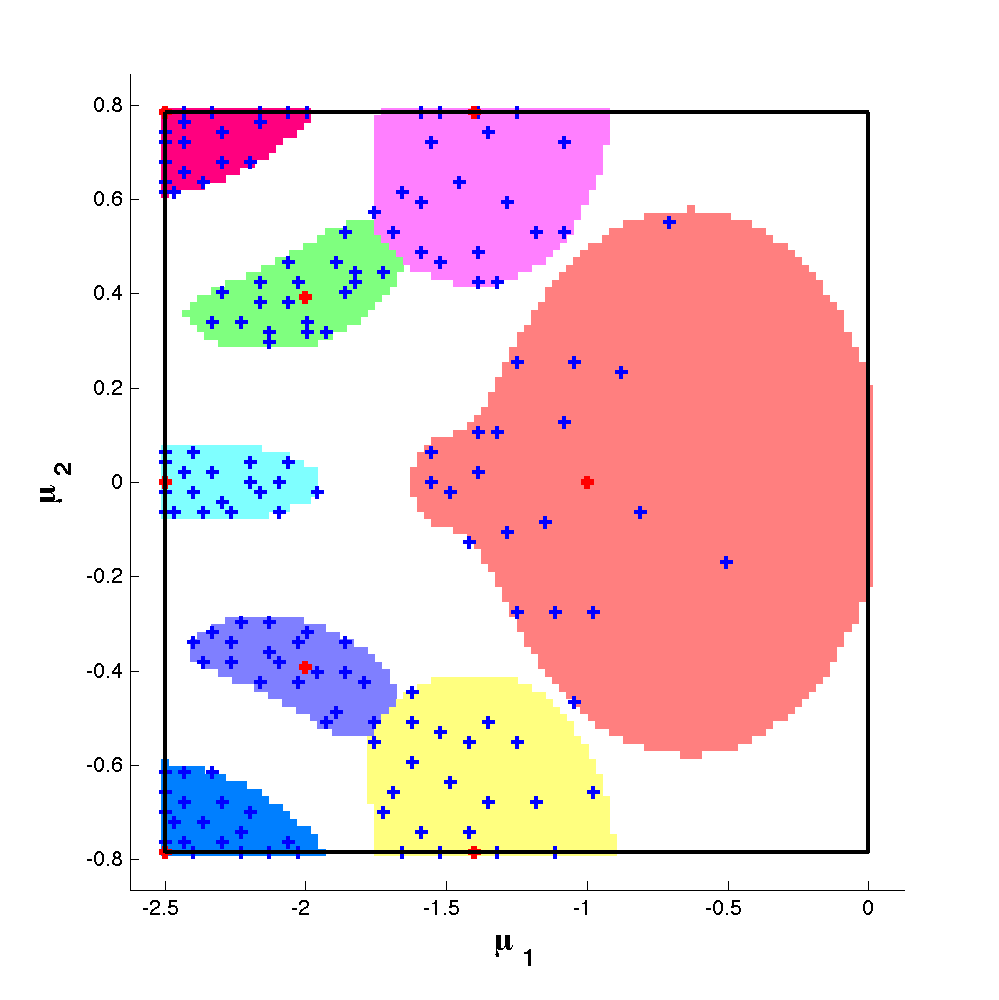}
\end{minipage}
}
\subfloat[Radius.]{
\begin{minipage}[t]{\three}
\includegraphics[Trim=\trimvalr,clip,width=\textwidth]{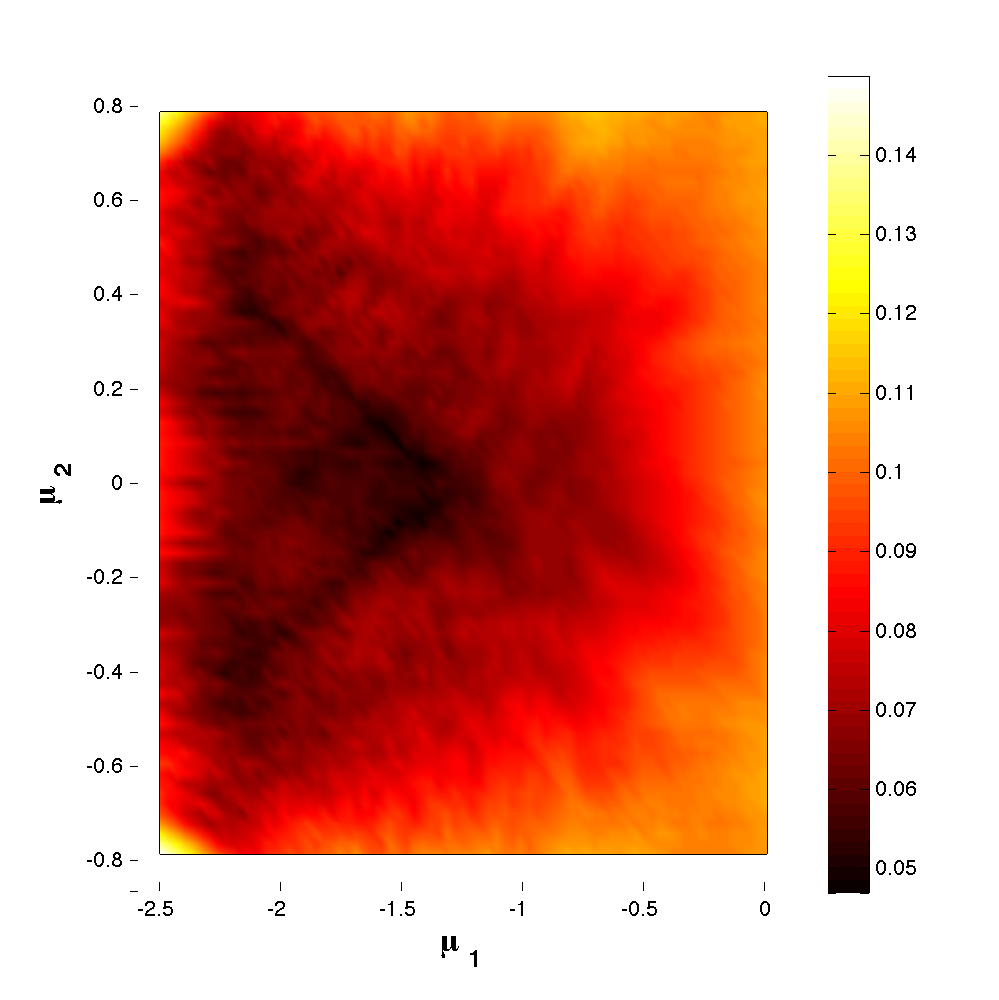}
\end{minipage}
}
\subfloat[Sample points.]{
\begin{minipage}[t]{\three}
\includegraphics[Trim=\trimval,clip,width=\textwidth]{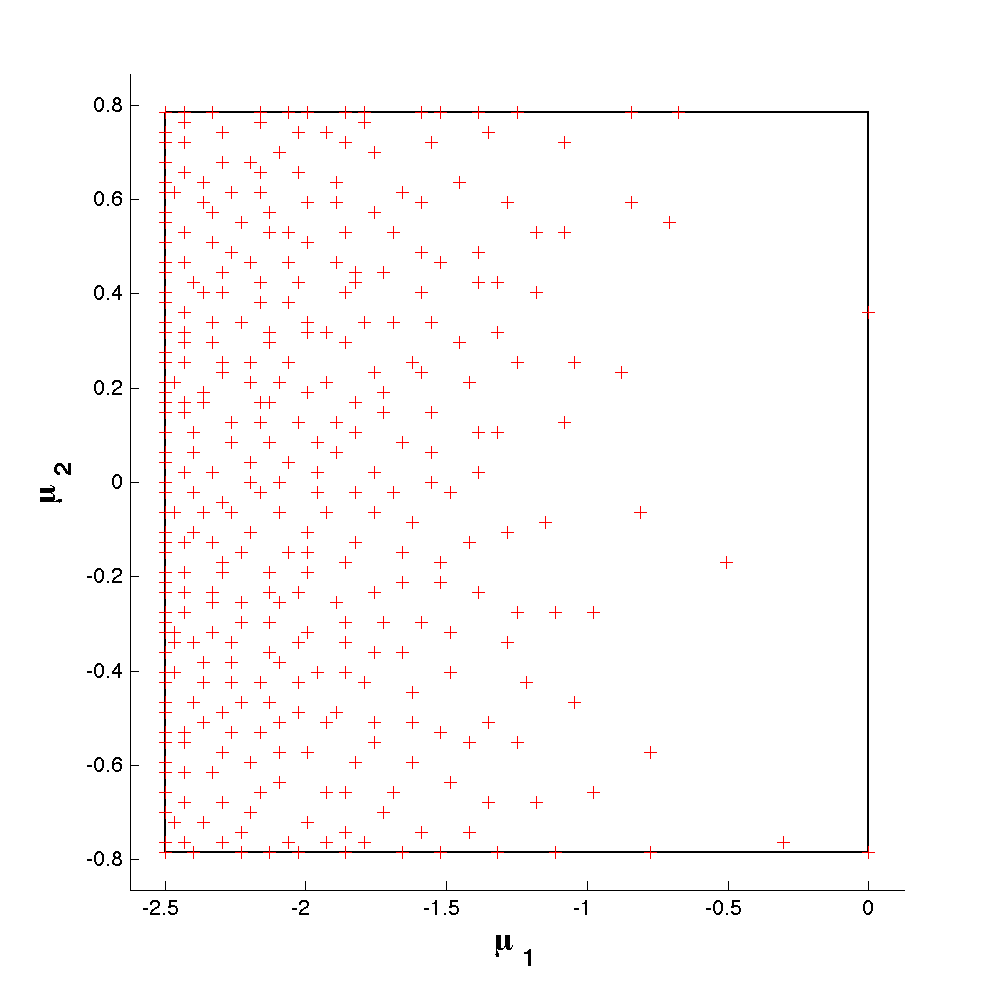}
\end{minipage}
}
\caption{
Test 4:  Local approximation spaces for selected parameter values (left), radius as a function of the parameters (middle) and sample points (right) for $N=20$ for the version using a fixed training set.
}
\label{fig:RB_N20}
\end{figure}

Figure \ref{fig:RB_N20} illustrates the local approximation spaces for some specific parameter values, the radius as a function of the parameter and the chosen sample points for $N=20$ for the version with a fixed training set. We can observe that the symmetry of the problem is reflected in all three quantities and that the shape of the local approximation spaces does vary.
Further, in Figure \ref{fig:RB_N20} we illustrate the training set, the local approximation spaces for some specific parameter values, the radius as a function of the parameter and the chosen sample points for $N=20$ for the version with an adaptive training set.

\begin{figure}[!ht]
  \captionsetup[subfigure]{labelformat=empty,width=\two}
  \centering
\subfloat[Using a fixed training set.]{
\begin{minipage}[t]{\two}
\includegraphics[width=\textwidth]{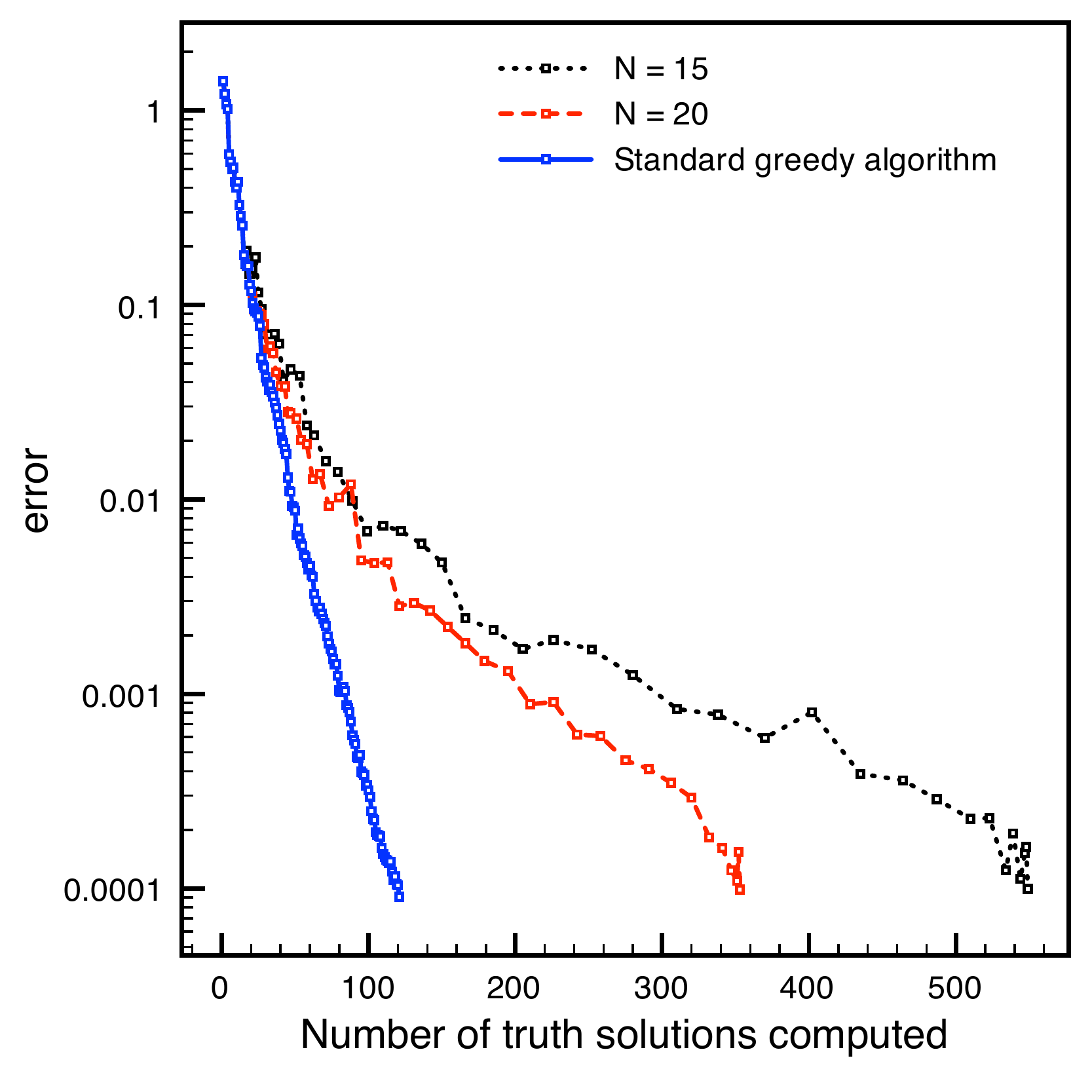}
\end{minipage}
}
\subfloat[Using an adaptive training set.]{
\begin{minipage}[t]{\two}
\includegraphics[width=\textwidth]{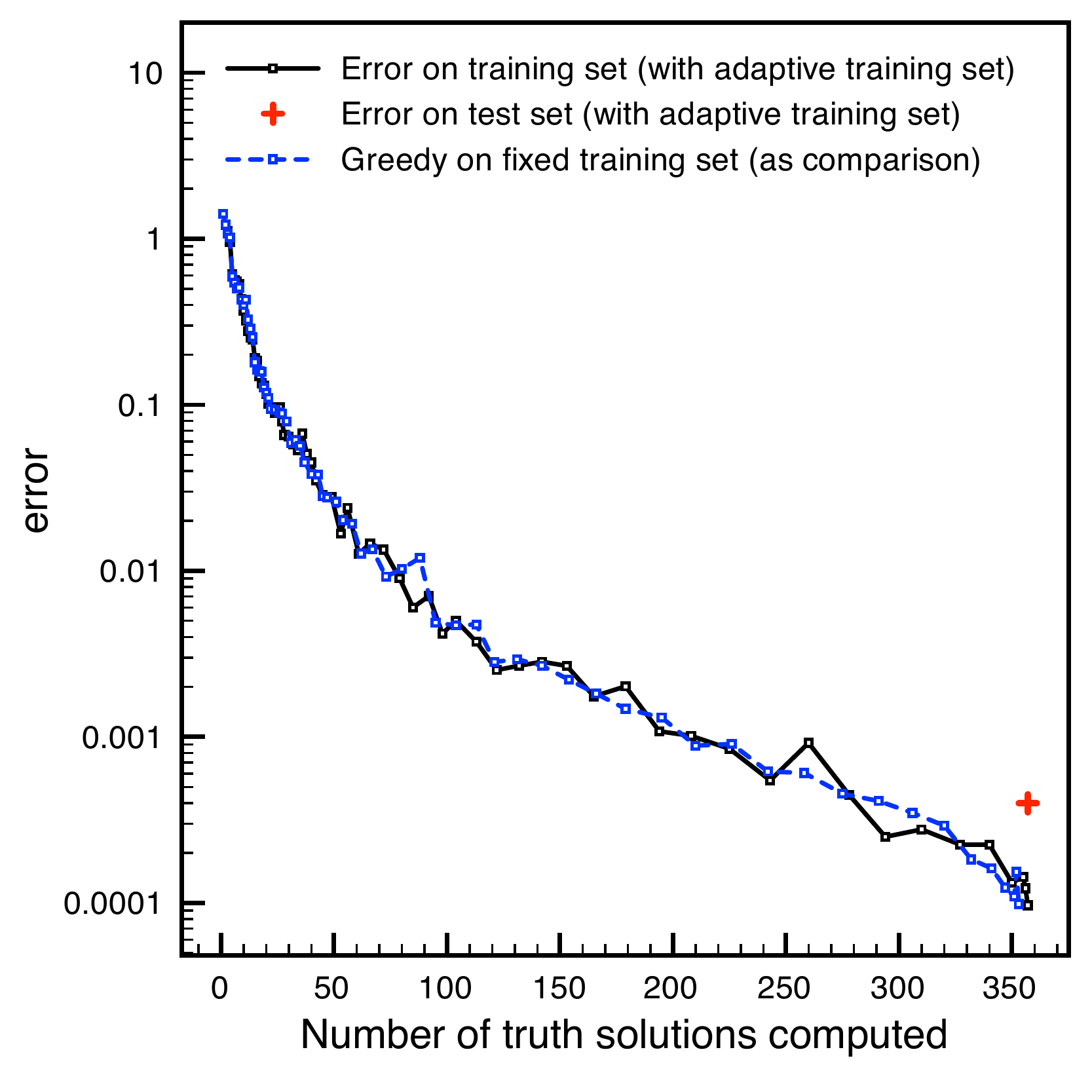}
\end{minipage}
}
\caption{
Test 4: Accuracy with respect to the number of computed truth solutions for the version with fixed training set (left) for $N=15,20$ and the version with an adaptive training set (right) for $N=20$.
}
\label{fig:RB_Conv}
\end{figure}

In Figure \ref{fig:RB_Conv} (left), we present the evolution of the error with respect to the number of computed truth solutions during the locally adaptive greedy algorithm for $N=15,20$ and also plot the convergence of the standard greedy algorithm, which of course converges faster but with a total number of basis functions equal to 121 at the on-line stage.
In Figure \ref{fig:RB_Conv} (right), the evolution of the accuracy on the training set of the version with the fixed and the adaptive training sets are compared. 
The value of about $4\cdot10^{-4}$ of the $L^\infty(\dom)$-norm on the test grid is surprisingly inaccurate.
Further numerical investigations have shown that the violation of the tolerance criterion is mainly focussed on the lower and upper boundary of the parameter space. 
The problem on the lower and upper boundary is caused how our implementation of generating the adaptive training set using FreeFem++ is placing the boundary points.
Indeed, the criterion is violated on only 0.35\% of all points of the test set.
If one studies the error on a test grid of $87\times 99$ points on $[-2.5,0]\times [-0.7,0.7]$ the error of the version with a fixed training set performs with an error of $1.13\cdot 10^{-4}$ and the version with an adaptive training set shows an error of $1.06\cdot 10^{-4}$.


We finally present in Table \ref{tab:times} some computation times of the on-line stage and compare it to a standard greedy-approach in order to highlight possible improvements that have to be made in a future work. 
As can be seen in Algorithm 4, the on-line stage consists of three steps: 1) Finding the $N$ closest sample points, 2) effecting the on-line ortho-normalization and 3) solving the linear system.
In Table \ref{tab:times} the average computing times over a sample of 10'000 random parameter values are given and the computations are done in Matlab restricting the computations to a single thread using a MacBook Pro with an Intel Core i7 processor with 2.66 GHz.
We observe that the proposed locally adaptive greedy algorithm for $N=20$ basis functions is significantly slower by about one order of magnitude than the standard greedy algorithm using $N=121$ basis functions. 
While solving the linear system is of course faster for 20 unknowns than for 121, most time is spent by finding the closest sample points and by effecting the ortho-normalization. 
Certain improvements can certainly be done respecting the implementation of the proposed algorithm such as building a tree-like structure in the parameter space that indicates the considered sample points that might be used in each box. This tree can be constructed off-line. At the on-line stage and for a given parameter value $\bmu$ one has to identify the box that contains $\bmu$ and then searching for the $N$ closest basis functions only among the sample points proposed by the selected box. 
Further, the triangular nature of the matrix $\gamma$ in the on-line ortho-normalization process can be exploited.
In addition, the error estimation/certification is the most expensive part in most cases of the on-line stage of the reduced basis method. 
These computation times have not be accounted for since we worked with the exact error computation. 
Summarizing, these computation times indicate possible bottlenecks of the present approach but do not show the entire truth.

\begin{figure}[!ht]
  \captionsetup[subfigure]{labelformat=empty,width=\four}
  \centering
\subfloat[Training points.]{
\begin{minipage}[t]{\four}
\includegraphics[Trim=\trimval,clip,width=\textwidth]{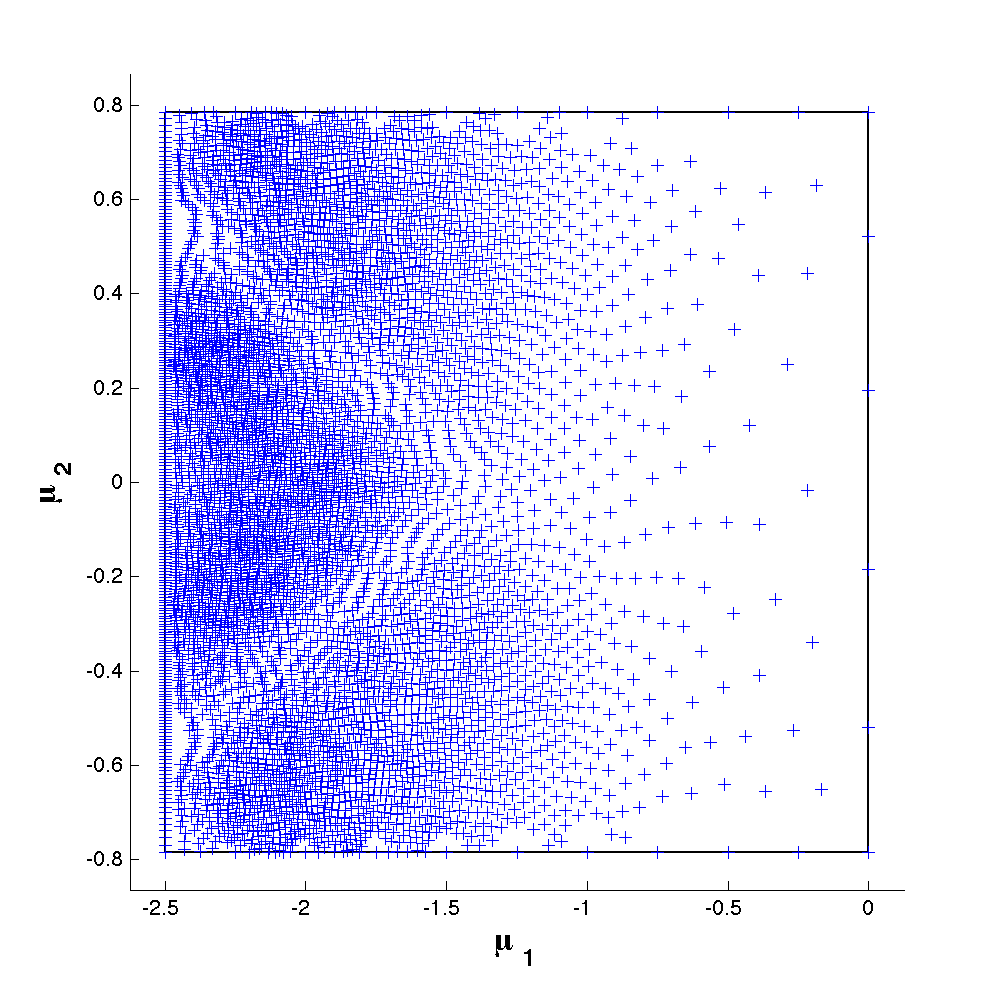}
\end{minipage}
}
\subfloat[Local approximation spaces.]{
\begin{minipage}[t]{\four}
\includegraphics[Trim=\trimval,clip,width=\textwidth]{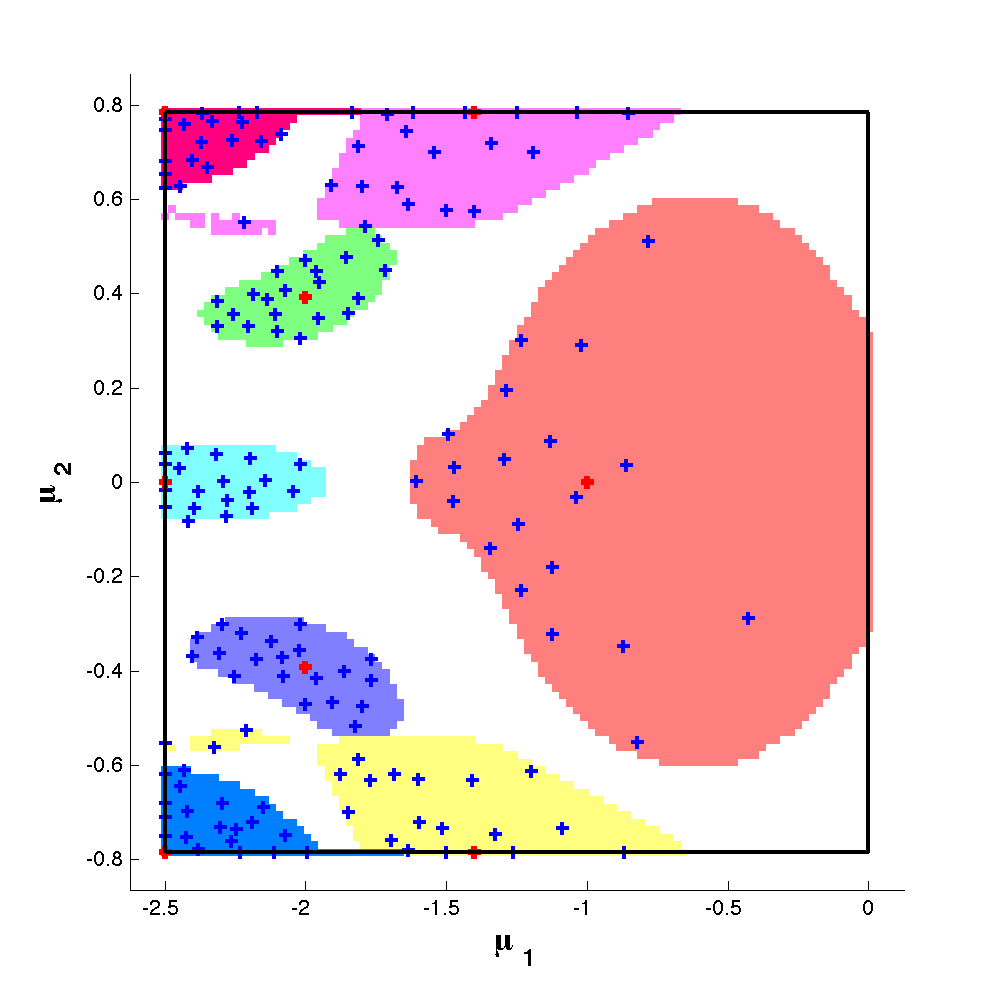}
\end{minipage}
}
\subfloat[Radius.]{
\begin{minipage}[t]{\four}
\includegraphics[Trim=\trimvalr,clip,width=\textwidth]{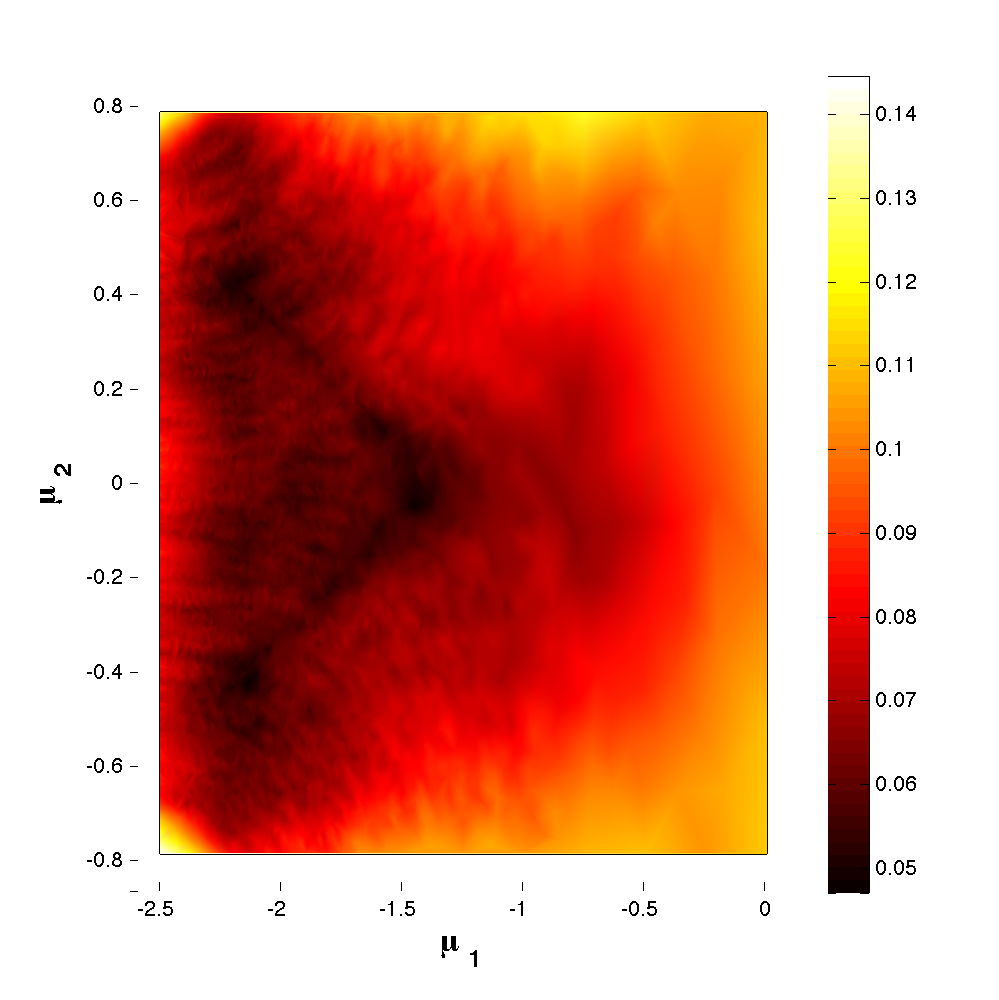}
\end{minipage}
}
\subfloat[Sample points.]{
\begin{minipage}[t]{\four}
\includegraphics[Trim=\trimval,clip,width=\textwidth]{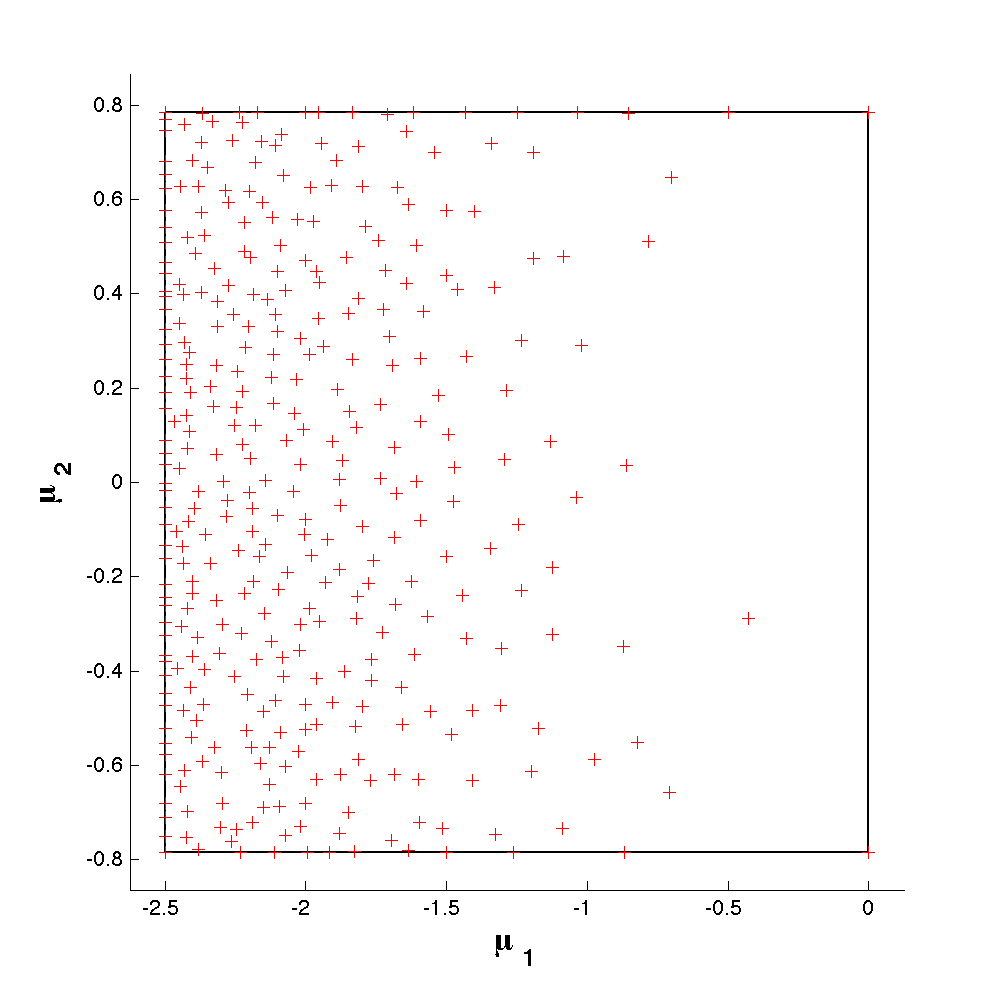}
\end{minipage}
}
\caption{
Test 4: Adaptive training set, local approximation spaces for selected parameter values, radius as a function of the parameters and sample points for $N=20$ using the version using an adaptive training set.
}
\label{fig:RB_N20_ATS}
\end{figure}

\begin{table}
\footnotesize
\begin{center}
\begin{tabular}{ | l | c | c | c | c |}
\hline
\multicolumn{1}{| c |}{[ms]}
& 
\begin{minipage}{2cm}
Total on-line time 
\end{minipage}
& 
\begin{minipage}{2cm}
\vspace{2pt}
Search of $N$-closest sample points
\vspace{2pt}
\end{minipage}
& 
\begin{minipage}{2cm}
Ortho-normalization 
\end{minipage}
& 
\begin{minipage}{2cm}
Solving the linear system 
\end{minipage}
\\
\hline
\begin{minipage}{2cm}
\vspace{2pt}
Fixed training set
\vspace{2pt}
\end{minipage}
&
2.9375
&
2.244  & 0.59266  & 0.10083
\\
\hline
\begin{minipage}{2cm}
\vspace{2pt}
Adaptive training set
\vspace{2pt}
\end{minipage}
&
2.9182
&
2.2285 & 0.58874 & 0.10091
\\
\hline
\begin{minipage}{2cm}
\vspace{2pt}
Standard greedy algorithm
\vspace{2pt}
\end{minipage}
&
0.3977
&
0 & 0 & 0.3977
\\
\hline
\end{tabular}
\end{center}
\label{tab:times}
\caption{Test 4: Average on-line computation times over a random set of 10'000 parameter values.}
\end{table}

%
%

\section{Conclusions}
\label{sec:conc}
We presented a framework that can be seen as a generalization of the classical greedy algorithm that is widely used in reduced basis methods. 
The presented algorithm introduces local approximation spaces (in the parameter space) that also account for local anisotropic behavior instead of a global approach.
The key idea was to consider the $N$ closest basis functions (for a fixed $N$) where the distance is measured in an empirically built distance function which can be constructed on the fly and which is problem-dependent.
In consequence, the number $N$ of basis functions considered for an approximation is user controlled. 
Further, it is proposed to optionally place the training set uniformly distributed in this empirically constructed (approximative) metric space. 
Numerical tests illustrate the characteristics of this approach. 
In our future work, the emphasis will be shed on the implementation of this algorithm for high dimensional parameter spaces, the main challenge in this direction will be to design properly a set of equi-distributed points with an adapted non isotropic metric (in particular our current implementation). 

\section*{Acknowledgement}

The authors want to thank the help and discussions with Frederic Hecht  about the adaptation of FreeFem++ to their purpose. This work was also initiated out of the discussions during the 
workshop LJLL-SMP: ''REDUCED BASIS METHODS IN HIGH DIMENSIONS'' held in June 2011, organized with the support of the Fondation de Sciences Math\'ematiques de Paris.




\bibliographystyle{abbrv}
\bibliography{Biblio}

\end{document}